\documentclass[final,3p,times,onecolumn]{elsarticle}

\usepackage{subfig}
\usepackage{amsmath, bm}
\usepackage{upgreek}
\usepackage{multirow}
\usepackage[export]{adjustbox}
\usepackage[linesnumbered,ruled, noend]{algorithm2e}

\DeclareMathOperator*{\Afem}{A}

\makeatletter
\def\ps@pprintTitle{%
  \let\@oddhead\@empty
  \let\@evenhead\@empty
  \def\@oddfoot{\reset@font\hfil\thepage\hfil}
  \let\@evenfoot\@oddfoot
}
\makeatother

\begin{document}

\begin{frontmatter}



\title{Physically recurrent neural networks for path-dependent heterogeneous materials: embedding constitutive models in a data-driven surrogate}

\author[a]{M. A. Maia}
\author[a]{I. B. C. M. Rocha}
\author[b,c]{P. Kerfriden}
\author[a]{F. P. van der Meer}

\address[a]{{Delft University of Technology, Department of Civil Engineering and Geosciences}, {P.O.Box 5048}, 2600GA, Delft, {The Netherlands}}
\address[b]{{Mines ParisTech (PSL University), Centre des mat\'{e}riaux}, {63-65 Rue Henri-Auguste Desbru\`{e}res BP87}, F-91003, \'{E}vry, France}
\address[c]{{Cardiff University, School of Engineering}, {Queen's Buildings}, The Parade, Cardiff, {CF24 3AA}, {United Kingdom}}

\begin{abstract}

Driven by the need to accelerate numerical simulations, the use of machine learning techniques is rapidly growing in the field of computational solid mechanics. Their application is especially advantageous in concurrent multiscale finite element analysis (FE$^2$) due to the exceedingly high computational costs often associated with it and the high number of similar micromechanical analyses involved. To tackle the issue, using surrogate models to approximate the microscopic behavior and accelerate the simulations is a promising and increasingly popular strategy. However, several challenges related to their data-driven nature compromise the reliability of surrogate models in material modeling. The alternative explored in this work is to reintroduce some of the physics-based knowledge of classical constitutive modeling into a neural network by employing the actual material models used in the full-order micromodel to introduce non-linearity. Thus, path-dependency arises naturally since every material model in the layer keeps track of its own internal variables. For the numerical examples, a composite Representative Volume Element with elastic fibers and elasto-plastic matrix material is used as the microscopic model. The network is tested in a series of challenging scenarios and its performance is compared to that of a state-of-the-art Recurrent Neural Network (RNN). A remarkable outcome of the novel framework is the ability to naturally predict unloading/reloading behavior without ever seeing it during training, a stark contrast with popular but data-hungry models such as RNNs. Finally, the proposed network is applied to FE$^\mathrm{2}$ examples to assess its robustness for application in nonlinear finite element analysis.

\end{abstract}




\end{frontmatter}


\newcommand{\fesqr}{FE$^2$}
\newcommand{\pd}[2]{\frac{\partial{#1}}{\partial{#2}}}
\newcommand{\strainvec}{\bm{\varepsilon}}
\newcommand{\stressvec}{\bm{\sigma}}
\newcommand{\strainstresscurve}{\bm{\varepsilon}\text{-}\bm{\sigma}}
\newcommand{\intvarvec}{\bm{\alpha}}
\newcommand{\macrostrainvec}{\bm{\varepsilon}^{\Omega}}
\newcommand{\macrostressvec}{\bm{\sigma}^{\Omega}}
\newcommand{\microstressvec}{\bm{\sigma}^{\omega}}
\newcommand{\microstrainvec}{\bm{\varepsilon}^{\omega}}
\newcommand{\macrodispfield}{\mathbf{u}^\Omega}
\newcommand{\microdispfield}{\mathbf{u}^\omega}
\newcommand{\weights}{\mathbf{W}}
\newcommand{\biases}{\mathbf{b}}
\newcommand{\hidden}{\mathbf{h}}
\newcommand{\histunits}{\mathbf{h}}
\newcommand{\stiffmatrix}{\mathbf{D}}
\newcommand{\xvec}{\mathbf{x}}
\newcommand{\Bmat}{\mathbf{B}}
\newcommand{\macrodmat}{\mathbf{D}^{\Omega}}

\section{Introduction}
\label{sec:intro}

Driven by the need to accelerate numerical simulations, machine learning (ML) techniques are rapidly growing in the field of computational solid mechanics. The use of surrogate models in particular is especially advantageous in concurrent multiscale finite element analysis (\fesqr). In this approach, each integration point at the macroscale is linked to a microscopic model, allowing complex materials such as composite materials to be explicitly modeled using relatively simple constitutive models. Although appealing, the computational cost associated with the concurrent Finite Element (FE) simulations is often prohibitive, hindering its widespread use in practical engineering-scale applications. In that sense, another layer of approximation, for instance based on ML techniques, that help accelerate \fesqr \ simulations is key. 

However, despite the successful use of this type of technique in several fields (\textit{e.g.}, speech recognition, language processing, and biomedical sciences), the use of surrogate models in material modeling is still filled with challenges and open issues. One of them relates to the data scarcity that exposes how data-driven models do not perform as well in extrapolation as they do within their training range. In general, another shortcoming of black-box models is their lack of interpretability. Although not a limiting aspect by itself, the limited interpretability these models offer can make it difficult to understand whether the fitted model is capable of yielding physically consistent solutions for loading paths different than those seen during training.  

In this scenario, the definition of a Design of Experiments (DoE) strategy plays an important role in the generalization capabilities of the surrogate model. The designer needs to take into account the computational (or experimental) budget available to collect the data and what type of loading scenarios are expected to be experienced by the model. The former is known \textit{a priori}, but the latter is not known exactly until the finite element simulation itself is run. As a consequence, large training and testing datasets are usually employed in an attempt to cover all different loading scenarios for a given strain/stress range. 

When path-dependency is present, conceiving an efficient DoE becomes even more complex since stresses now depend on the history of the material, leading to potentially infinite parameter space. One way to bypass the computationally intensive approach is to run and incorporate new load paths on-the-fly, updating the surrogate model as necessary as in the works of \citet{Gouryetal2016} and \citet{Rocha2021}. In addition to the definition of the DoE, the physics-based assumptions about the material should also be wisely taken into account when choosing the modeling approach. For instance, conventional Artificial Neural Networks (ANN) with strains and stresses taken as input and output, respectively, will fail to capture phenomena such as elastic unloading. This scenario is illustrated by \citet{VLASSIS2021} and occurs due to the unique mapping between inputs and outputs that does not reflect two different stress states for the same strain. One way to differentiate the different paths is by augmenting the feature space of the ANN with extra variables that carry partial information about the history of stresses and/or strains \cite{LEFIK20091785, Huang2020}. 

Another highly popular strategy to handle path-dependency is to use Recurrent Neural Networks (RNNs), a variation of ANNs capable of learning from sequential data. This is done by incorporating the previous state of the network when making predictions for the current state. However, this approach often suffers from short-term memory when dealing with long sequences. To address that issue, more complex architectures with more model parameters and fine control of the flow of information retained through a sequence of information were proposed (\textit{e.g.}, GRU and LSTM). These networks gained popularity and have been used to model a wide variety of materials that show history-dependency \cite{Ghavamian2019, WU2020113234, Mozaffar26414, Chen2021, Koeppe2021, LOGARZO2021}. 

However, RNNs are still severely limited by the curse of dimensionality associated with sampling arbitrarily long strain paths. The overview on the potential of RNNs for modeling path-dependent plasticity presented by \citet{GORJI2020} illustrates the multiple ways in which one can define a DoE and argue that GRUs, in particular, can be used to model the plastic response of materials provided a rich enough training dataset is used. However, the authors do not investigate further how these networks would perform in scenarios slightly different than those seen in training or their efficiency in an \fesqr \ framework.

Recently, a new class of ANNs called Physics-Informed Neural Networks (PINNs) was proposed by \cite{RAISSI2019} to solve any given laws of physics described by general nonlinear partial differential equations. For training PINNs, in addition to the data used to model the governing equation, the loss function is augmented with information on the imposed physical constraints (\textit{e.g.} initial and boundary conditions). In this approach, automatic differentation plays a key role in allowing the outputs of the network to be differentiated with respect to the input. 

Applications of PINNs as surrogate constitutive models are still relatively new, but a few works showcase their potential in modeling elastoplastic \cite{HAGHIGHAT2021, Eghbalianetal2022} and elastic-viscoplastic behaviors \cite{Arora2022}. In \cite{HAGHIGHAT2021} and \cite{Arora2022}, the loss function is also incremented with terms related to the physics constraints assumed by the material models considered (\textit{e.g.} the Karush-Kuhn-Tucker conditions on the yield function). While both works approximate the displacements and the stress based on the position of the material point at the macroscale and on the material properties, \citet{HAGHIGHAT2021} further investigates the applicability of their surrogate model in a sensitivity analysis study. Finally, \citet{Eghbalianetal2022} proposed an architecture where the nonlinear incremental elasticity and the strain decomposition are hardwired into the networks's architecture and the loss function. Thus leading to improved capability in predicting unseen loading scenarios over conventional ANNs.

Inspired by PINNs, \citet{Masi2020} tailored a network architecture capable of ensuring thermodynamically consistent predictions by evaluating the numerical derivatives of the network with respect to its inputs. Later, \citet{Masi2021} proposed a key extension to the approach, in which the identification of the internal variables is no longer user-dependent. For that purpose, an encoder-decoder architecture is incorporated to identify the minimum number of internal variables from the set of internal coordinates of the system (\textit{e.g.}, displacement fields, internal forces, etc.)  in an unsupervised manner. The feature allows the recovery of the full-field information of the microstructure by the decoder.
Another strategy that takes into account the derivatives of the approximated functions to impose thermodynamically consistency is proposed by \citet{VLASSIS2021}. 
The authors treat the yield function as a signed distance function level set and formulate a supervised learning task that deduces the learned evolving yield function against a monotonically increasing accumulated plastic strain. As a result, cyclic loading paths can be predicted based only on monotonic data.
For a recent and comprehensive review on the state-of-the-art literature of ANNs in the constitutive modeling of composite materials, the reader is referred to \citet{LIU2021109152}. 

Different from regular ANNs or RNNs, the Deep Material Network (DMN) proposed by \citet{Liu2019} assigns physical meaning to the hyperparameters. The DMN learns the hidden topology representation of the micromodel (RVE) based on the stiffness matrices and residual stresses of the different material phases. A major advantage of DMNs is the ability to extrapolate to nonlinear behavior based only in elastic snapshots. %
Several follow-up publications and improvements on DMN can be found in the literature. The author's latest work \cite{Liu2021} is dedicated to tackling multiscale failure analysis, a far less explored avenue by the community. Other succesful applications in this field can be found in the works of \citet{KERFRIDEN2013}, \citet{BESSA2017} and \citet{OLIVER2017}, where different reduced-order modeling strategies were employed. 

In an alternative approach that incorporates some level of physics in the surrogate model, \citet{Fuhg2021} combines a physical model that accounts partially for the constitutive behavior of the material and a Gaussian Process (GP) correction to locally improve the baseline physical model. Although this boosts accuracy inside the training space, in practice, it does not translate to benefits in terms of extrapolation. Another alternative that retains some of the physics embedded in the full-order solution is model order reduction. These methods are frequently employed to alleviate the computational cost of \fesqr \ \cite{GHAVAMIAN2017, Rocha2018} and usually offer better generalization properties to unseen points than common surrogate models. The drawback is that they are inherently slower \cite{Rocha2020}. Such techniques can also be coupled with RNNs \cite{vijayaraghavan2021neuralnetwork}, ANNs \cite{Huang2020} and GPs \cite{GUO2021} to quickly infer the coefficients of the reduced basis approximation for arbitrary parameter values. 

Although a large body of literature has been devoted to applying ML techniques in the solid mechanics field, the gaps and issues left by data-driven surrogate models are still an open issue, compromising their reliable and widespread use in practical applications. In this work, a new network design is proposed specifically for the modeling of path-dependent materials and to accelerate concurrent finite element simulations. In Section \ref{sec:fe2} the \fesqr \ method is presented, while in Section \ref{sec:rnn} one of the most popular approaches to tackle the computational bottleneck originated from it is briefly discussed. In Section \ref{sec:mnn}, the main features of the novel neural network are described. In Section \ref{sec:doe}, the Design of Experiments and methodology adopted for the comparative study shown in Section \ref{sec:results} is described. In this study, the performance of the proposed network is compared to a RNN for a single-scale problem. In Section \ref{sec:numex} the novel approach is integrated into an \fesqr \ framework and tested in two applications for robustness and accuracy. Finally, conclusions are presented in Section \ref{sec:conclusion}.

\section{Concurrent multiscale analysis}
\label{sec:fe2}

\begin{figure}[!ht]
\centering
\subfloat[Schematic representation of \fesqr]{\label{fig:fe2}\includegraphics[width=0.35\textwidth]{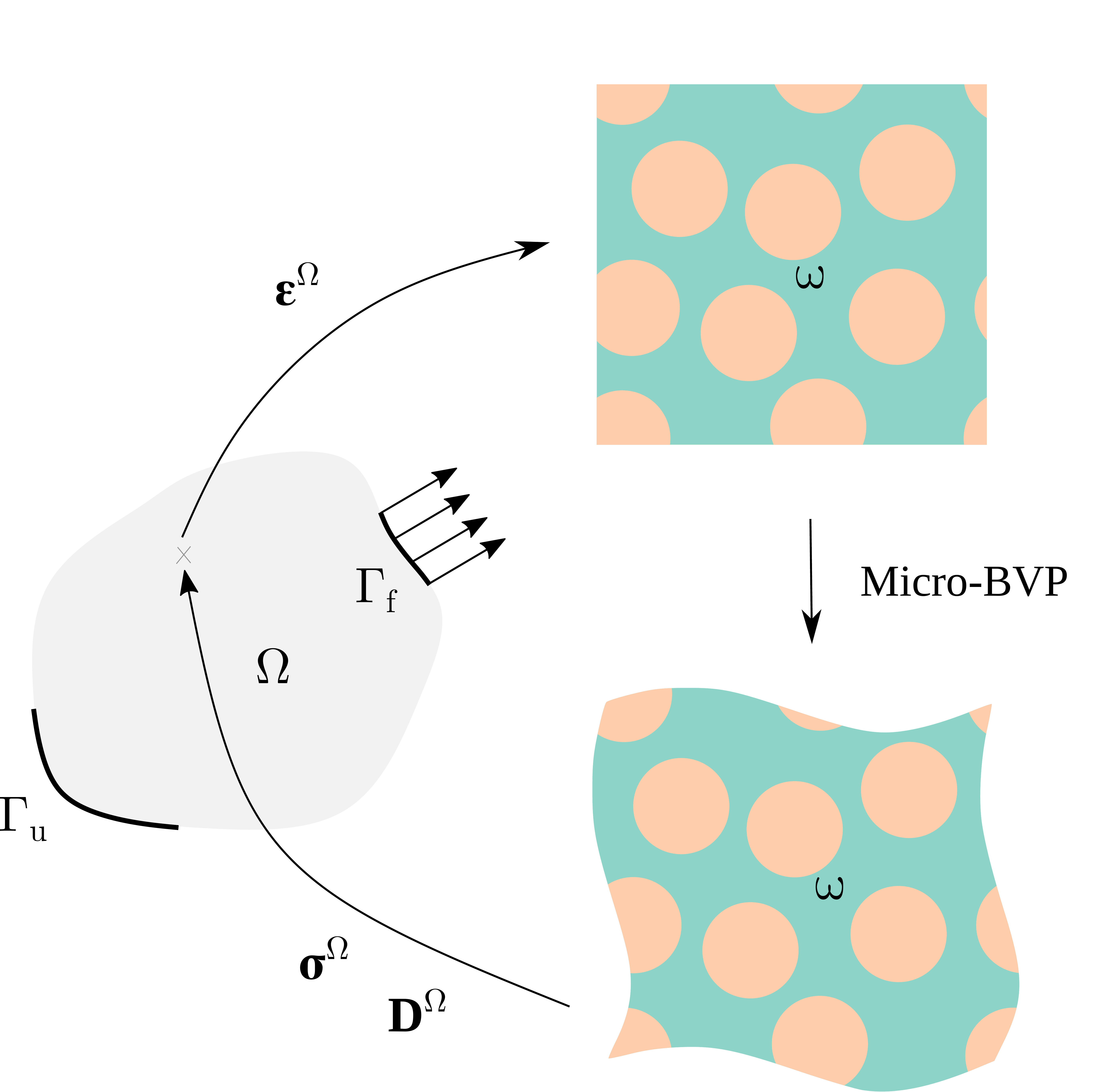}}
\hfill
\subfloat[Controlling nodes of RVE]{\label{fig:rve} \includegraphics[width=0.45\textwidth]{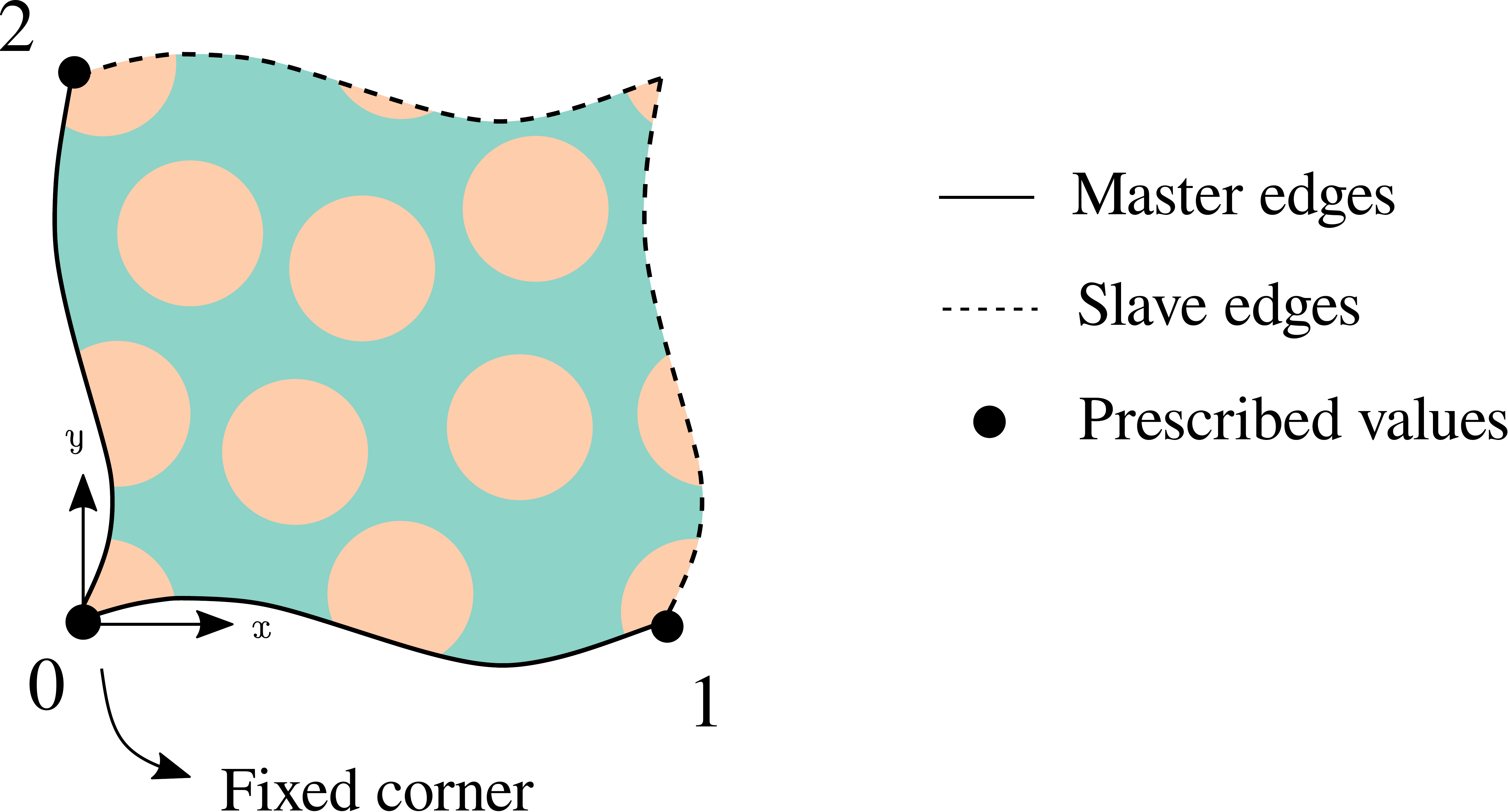}}
\caption{Scheme of \fesqr \ framework and definition of the boundary value problem on RVE}
\end{figure}

Let $\Omega$ define the macroscopic domain being modeled. To find the internal stresses and displacement field of such body in absence of body forces, a boundary value problem that satisfies the following equilibrium equations is defined:
\begin{equation}
\label{eq:dirichlet}
\begin{aligned}
\nabla \macrostressvec = \mathbf{0} 
\end{aligned}
\end{equation}
where $\nabla$ is the divergence operator and $\macrostressvec$ is the macroscopic stress, which depends on the macroscopic displacement field $\macrodispfield$ (for simplicity, this dependence is omitted). The governing equations are subjected to the boundary conditions: 
\begin{equation}
\label{eq:constraints}
\begin{aligned}
\macrostressvec \ \mathbf{n} = \mathbf{t}^{\Gamma_f} \ \ \ \text{on} \ \Gamma_f \qquad
\macrodispfield = \mathbf{u}^{\Gamma_u} \ \ \ \text{on} \ \Gamma_u  
\end{aligned}
\end{equation}
where $\mathbf{n}$ is the normal to the surface $\Gamma_f$ and $\mathbf{u}^{\Gamma_u}$ and $\mathbf{t}^{\Gamma_f}$ represent a set of Neumann and Dirichlet boundary conditions acting on the body surface such that $\Gamma_u \cap \Gamma_f = \emptyset$ as illustrated in Fig. \ref{fig:fe2}. To relate strains and stresses, a constitutive model $\mathcal{D}$ is required:
\begin{equation}
\label{eq:constmodel}
\macrostressvec = \mathcal{D}^{\Omega} \ (\macrostrainvec, \intvarvec^{\Omega}) 
\end{equation}
where $\intvarvec^{\Omega}$ are history variables that account for path-dependency and $\macrostrainvec$ is the macroscopic strain defined under small displacement assumptions as:
\begin{equation}
\label{eq:strainvec}
\macrostrainvec = \frac{1}{2} \Big( \nabla \macrodispfield + (\nabla \macrodispfield)^\text{T}  \Big)
\end{equation}

In the concurrent multiscale approach, the model $\mathcal{D}^{\Omega}$ is not directly formulated but is instead obtained by nesting a lower scale finite element model to each integration point. In that scale, the microscopic stucture of complex materials can be explicitly modeled using simpler constitutive models for each of the components. Further discussion on how to solve the microscopic problem and link both scales is shown in Section \ref{subsec:micro} and \ref{subsec:homo}.

To solve the boundary value problem at the macroscale, the Finite Element (FE) method is employed to discretize the domain $\Omega$ into a number of elements connected by nodes with $N$ degrees of freedom. The global equilibrium is solved iteratively in its discretised weak form:
\begin{equation}
\label{eq:resfe2}
\mathbf{r} = \mathbf{f}^\Gamma - \mathbf{f}^\Omega (\macrodispfield) = \mathbf{0}
\end{equation}
where $\mathbf{r} \in \mathbb{R}^N$ is a residual vector that goes to zero when equilibrium is reached, $\mathbf{f}^\Gamma \in \mathbb{R}^N$ is the global external vector that represents the Neumann boundary conditions and $\mathbf{f}^\Omega \in \mathbb{R}^N$ is the global internal force vector given by a volume integral:
\begin{equation}
\label{eq:intforcesupperscale}
\mathbf{f}^{\Omega} = \Afem_{e=1}^{n_e} \int_{\Omega_e} \Bmat_{e}^{\text{T}} \ \macrostressvec (\macrodispfield_e) \ \text{d} \Omega 
\end{equation}
where $\Afem$ is an assembly operator that takes into account the connectivities between the elements and the global system and $\Bmat$ is a matrix with the spatial derivatives of the shape functions used to interpolate nodal displacements. 
Finally, an iterative procedure is adopted to solve Eq. (\ref{eq:resfe2}) for the macroscopic displacement field:
\begin{equation}
\label{eq:iterativedisp}
\Delta \macrodispfield = \macrodispfield_n - \macrodispfield_o = - \mathbf{K}_o^{-1} \mathbf{r}_o
\end{equation}
where the subscripts $o$ and $n$ refer to old and new analysis increments, respectively and $\mathbf{K} \in \mathbb{R}^{N \times N}$ is the global tangent stiffness matrix given by:
\begin{equation}
\label{eq:tangentupperscale}
\mathbf{K} = \Afem_{e=1}^{n_e} \int_{\Omega_e} \Bmat_{e}^{\text{T}} \ \mathbf{D}_{e}^{\Omega} (\macrodispfield_e) \ \Bmat_{e} \ \text{d} \Omega 
\end{equation}
and $\macrodmat$ is the constitutive tangent matrix, discussed in Section \ref{subsec:homo}. 

The key difference to a classical FE simulation lies in the embedding of another FE model in the macroscopic integration points. Here, to obtain the internal forces in Eq. (\ref{eq:intforcesupperscale}) and the tangent stiffness matrix in Eq. (\ref{eq:tangentupperscale}) for a single integration point of the macroscale, one needs to run an entire FE model instead of a single evaluation of a homogenous material model. This is the most computationally expensive part of the framework and is where the proposed network aims to tackle. The approach here is to replace the solution to the microscopic problem (discussed in Section \ref{subsec:micro}) with a surrogate model, specifically a neural network. The homogenization procedure required to upscale the responses to the macroscale is discussed in Section \ref{subsec:homo}. 

\subsection{Microscopic scale}
\label{subsec:micro}

Let $\omega$ be a Representative Volume Element (RVE) of the microscopic material features whose behavior is to be upscaled. Assuming that the principle of separation of scales (\textit{i.e.} $\Omega \gg \omega$) holds, the two scales can be linked by enforcing: 
\begin{equation}
\label{eq:microdispfield}
\microdispfield = \macrostrainvec \xvec^\omega + \tilde{\mathbf{u}}
\end{equation}
where the linear displacement field is the result of the imposed macroscopic strains $\macrostrainvec$ and the fluctuation field $\mathbf{\tilde{u}}$ is the result of microscopic inhomogeneities.
The principle of separation of scales implies the stain averaging theorem that states that the macroscopic strains are considered uniform over the RVE domain:
\begin{equation}
\label{eq:strainavtheorem}
\macrostrainvec (\xvec^\Omega) = \frac{1}{\omega} \int_{\omega} \microstrainvec (\xvec_\omega) \ \text{d}\omega 
\end{equation}
where $\microstrainvec$ is the microscopic strain tensor. Therefore, the microscopic displacement field in Eq. (\ref{eq:microdispfield}) can only satisfy (\ref{eq:strainavtheorem}) if the fluctuation displacement field vanishes at the RVE boundary when upscaling quantities. An additional requirement on the fluctation field having zero resultant work at the boundaries arises from the Hill-Mandel principle. Both requirements are met using Periodic Boundary Conditions (PBC) to represent the behavior of a macroscopic bulk material point. Fig. \ref{fig:rve} illustrates the node groups and boundary edges needed to implement the PBC. In Section \ref{sec:doe}, the generation of $\strainstresscurve$ paths for the training of the surrogate models is obtained by setting a user-defined function to set the prescribed displacements $\mathbf{u}_1$ and $\mathbf{u}_2$.  

Finally, keeping the hypothesis of small strains, the stress equilibrium problem is described as:
\begin{equation}
\label{eq:microproblem}
\begin{aligned}
\nabla \microstressvec = \mathbf{0} \qquad 
\microstressvec = \mathcal{D}^{\omega} \ (\microstrainvec, \intvarvec^{\omega})  \qquad
\microstrainvec = \frac{1}{2} \Big( \nabla \microdispfield + (\nabla \microdispfield)^\text{T}  \Big)
\end{aligned}
\end{equation}
where $\microdispfield$ is the microscopic displacement field and $\microstressvec$ and $\microstrainvec$ are the microscopic stress and strain tensors, respectively. An analoguous procedure to the one detailed in Section \ref{sec:fe2} is used to find the microscopic displacement field (subjected to the periodic boundary conditions). Note that at this scale, regular physics-based material models $\mathcal{D}^\omega$ (\textit{e.g.} elastoplasticity, viscoelasticity, etc.) are employed to represent the constitutive behavior of the homogeneous material of the discretized elements. 

\subsection{Homogenization procedure}
\label{subsec:homo}
After convergence of the microscopic displacement field $\microdispfield$, the upscaling procedure is performed based on the Hill-Mandel principle. The principle postulates that the macroscopic stress power must equal the volume average of the microscopic power over the RVE. Considering the definition in Eq. (\ref{eq:microdispfield}), an expression similar to Eq. (\ref{eq:strainavtheorem}) is obtained for the homogenized stresses:
\begin{equation}
\label{eq:homogenizedexp}
\begin{aligned}
\macrostressvec =  \frac{1}{\omega} \int_{\omega}\microstressvec \text{d} \omega \qquad
\end{aligned}
\end{equation}

As for the macroscopic constitutive tangent stiffness $\macrodmat$, a probing operator $\mathcal{P}$ is applied on the global microscopic tangent stiffness matrix $\mathbf{K}^\omega$ without the need to invert it as proposed by \citet{NGUYEN2012}.

\section{Recurrent Neural Networks}
\label{sec:rnn}

In this section, a brief overview of the working mechanisms of Recurrent Neural Networks is presented. Although part of the paper is dedicated to comparing them with the novel approach, here the idea is to use well-known concepts from ANNs and RNNs to illustrate features of the proposed network in the following sections. 

As a starting point, consider a conventional feed-forward neural network to surrogate the nonlinear constitutive relationship of a path-independent material given by the following parametric regression model: 
\begin{equation}
    \label{eq:pred}
    \widehat{\boldsymbol{\sigma}}^{\Omega} = \mathcal{NN} \ (\boldsymbol{\varepsilon}^{\Omega}, \weights, \biases)
\end{equation}
where $\weights$ and $\biases$ are weights and biases calibrated through a fitting procedure based on observations of the actual microscopic model. During training, the strains are fed to the first neural layer (input layer) and values are propagated until the final layer (output layer) to give the predicted stresses $\widehat{\stressvec}^{\Omega}$. These are in turn compared to the ground truth value according to a loss function. Based on that, the model parameters are adjusted so that the error between the predicted stresses and the actual stresses is minimized:
\begin{equation}
    \label{eq:nnparameters}
    \weights, \biases = \textrm{argmin} \ \overline{\weights} \ \overline{\biases} \sum_{i \ \in \ \mathbf{X}} \| \ \boldsymbol{\sigma}_{i}(\macrostrainvec_i) - \widehat{\boldsymbol{\sigma}}_{i} (\macrostrainvec_i, \overline{\weights}, \overline{\biases}) \ \|^2
\end{equation}
where $\mathbf{X} \in \mathbb{R}^{n_{\varepsilon} \times N}$ is a snapshot matrix with $N$ pairs of $\macrostrainvec$-$\macrostressvec$ obtained from microscopic simulations. This setting is the most straight-forward way to map pairs of macroscopic strains and stresses but does not offer good generalization properties once path-dependency is introduced. In that case, one way to overcome the lack of history information is to extend their feature space with \textit{e.g.} previous (incremental) strains and/or stresses \cite{LEFIK20091785, Huang2020}.  

\begin{figure}[h]
\centering
\subfloat[Classical Recurrent Unit]{\label{fig:rnncell}\includegraphics[width=0.30\textwidth]{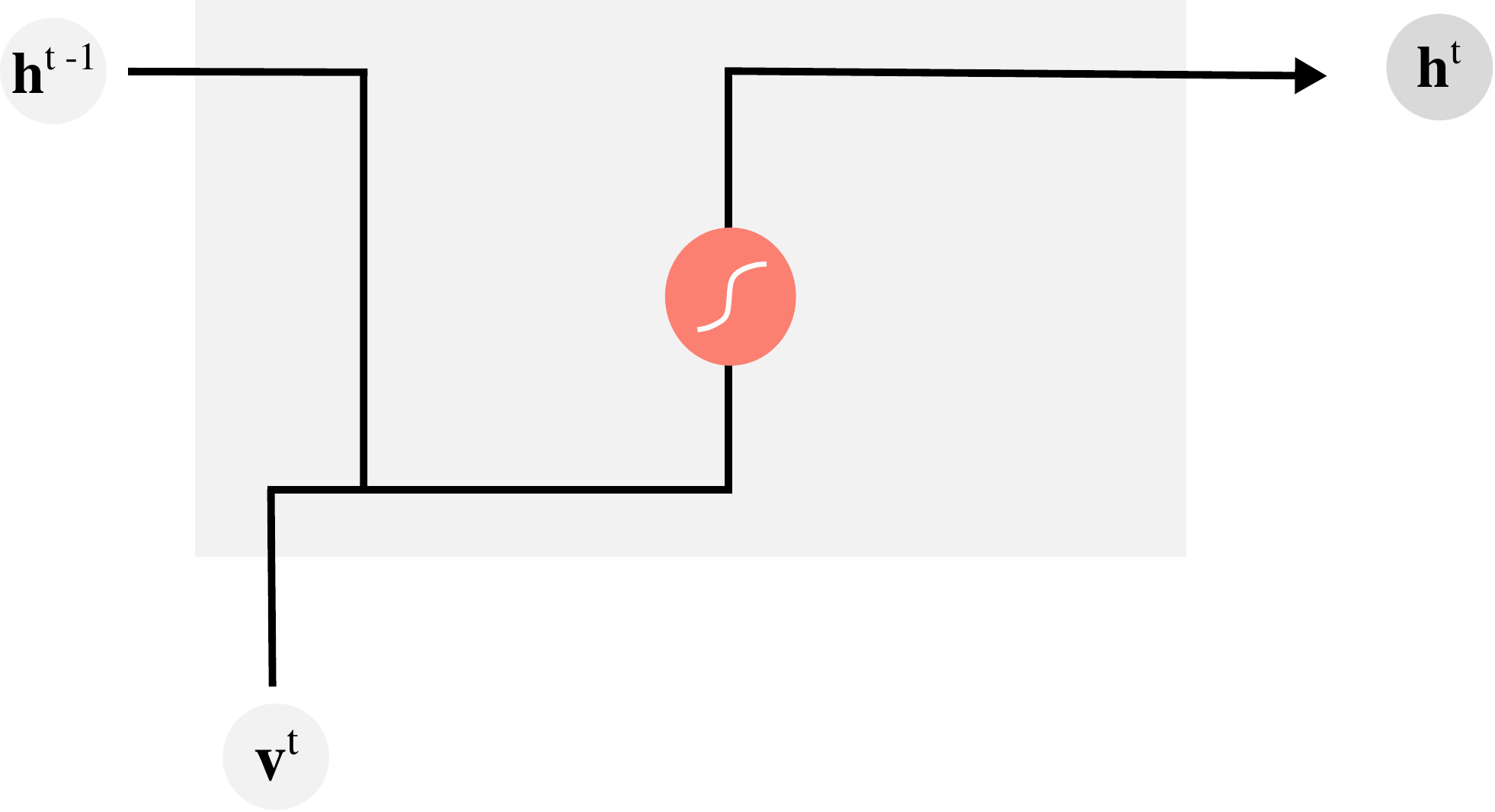}}
\hfill
\subfloat[Gated Recurrent Unit (GRU)]{\label{fig:grucell}\includegraphics[width=0.30\textwidth]{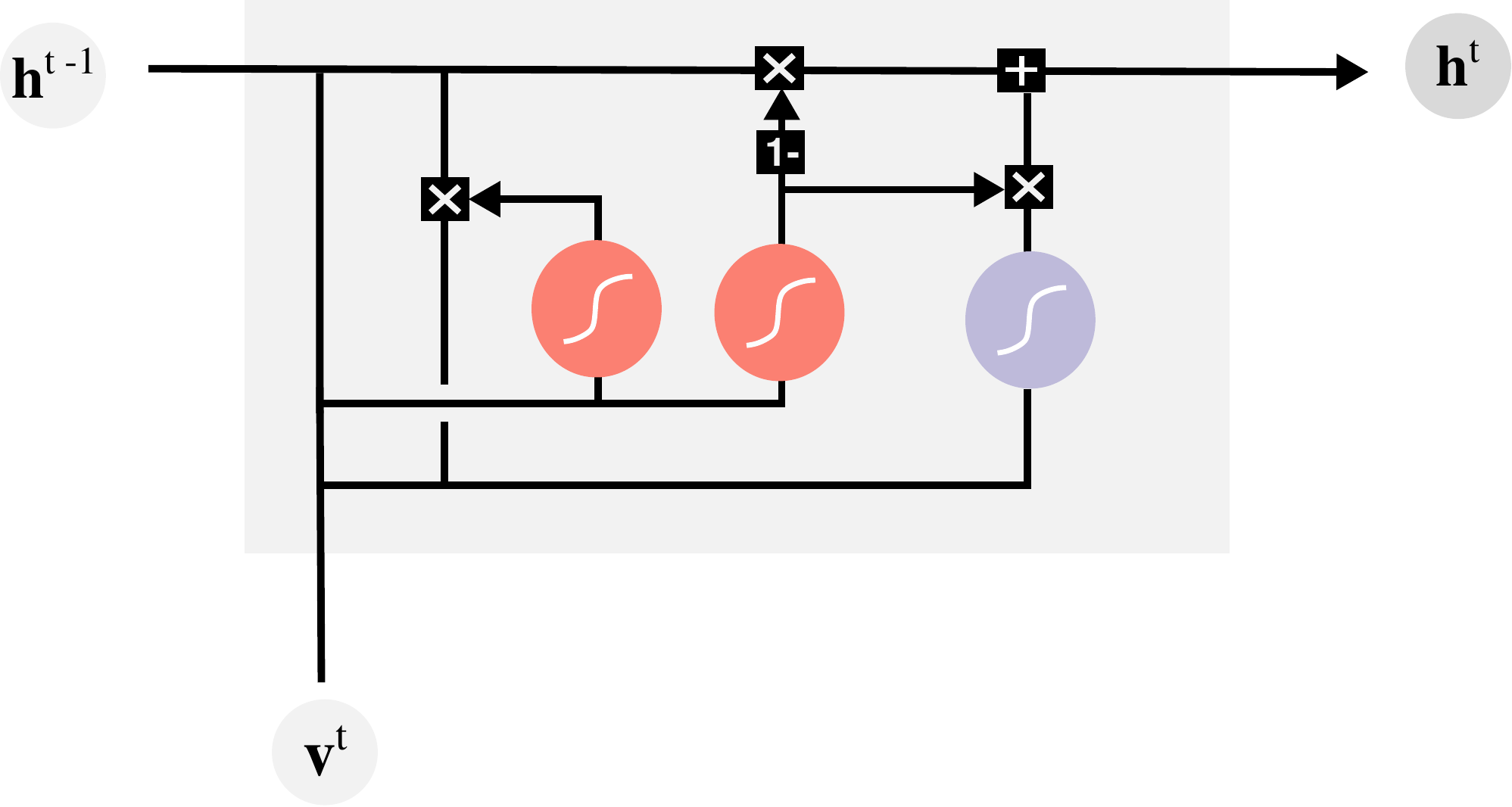}}
\hfill
\subfloat[Long-Short Term Memory Unit (LSTM)]{\label{fig:lstmcell}\includegraphics[width=0.30\textwidth]{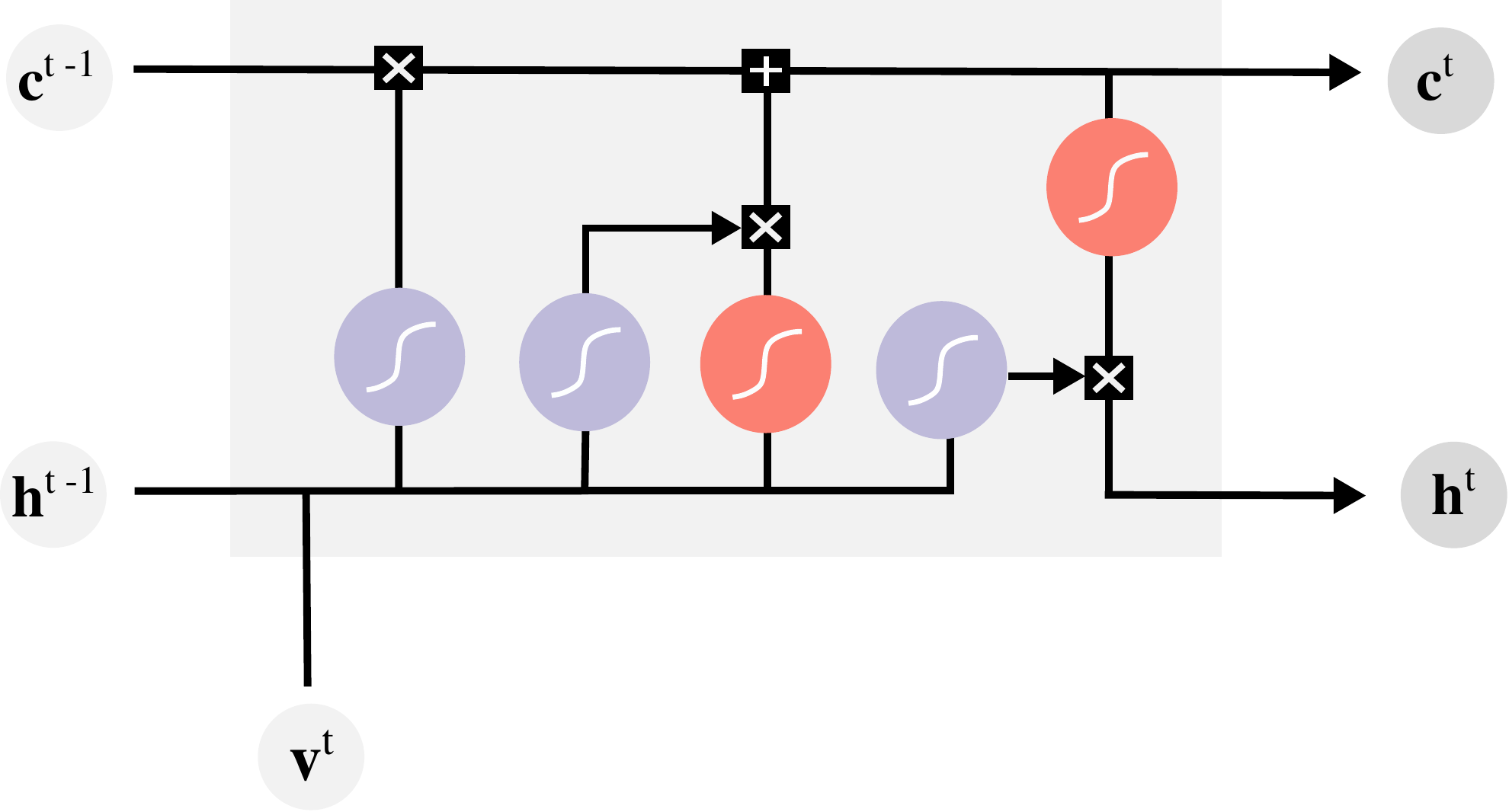}}
\hfill
\caption{Different architectures for recurrent units/cells: red circles correspond to the element-wise application of the sigmoid function, while purple circles correspond to the tanh function.}
\end{figure}

As an alternative, RNNs offer additional parameters (\textit{i.e.} the hidden state) and mechanisms (\textit{i.e.} the gates that control the flow of information being propagated) to learn history information from sequential data in an implicit way. These parameters describe the evolution of the so-called hidden state and can encapsulate information from previous iterations without the need to introduce history variables in the feature space. In a regular RNN, the outputs and hidden state are given by:
\begin{equation}
\label{eq:rnn}
\begin{aligned}
\hidden^{t}  = \upphi ( \weights_1 \mathbf{v}^{t} + \weights_s \hidden^{t-1} + \biases_s) \\
\widehat{\stressvec}^t = \upphi ( \weights_2 \hidden^{t} + \biases_2) 
\end{aligned}
\end{equation}
where $\upphi(\cdot)$ is an activation function, $\weights_s$ and $\biases_s$ are the additional model parameters (compared to conventional feed-forward neural networks), $\mathbf{v}^t$ are the current neuron values coming from the last layer and $\hidden^t$ and $\hidden^{t-1}$ are current and previous states, respectively. This arrangement allows the network to learn how stress evolves for a sequence of strains instead of bulding a regression model from independent stress-strain pairs and is illustrated in Fig. \ref{fig:rnncell}. However, in practice, the efficiency of RNNs are impeded by vanishing gradient problems and are not suitable for long-term history dependent problems. 

To overcome that, more sophisticated architectures (more popularly known as cells) have been proposed. Among the most popular ones are the Gated Recurrent Unit (GRU) and the Long-Short Term Memory (LSTM), illustrated in Fig. \ref{fig:grucell} and \ref{fig:lstmcell}, respectively. The internal mechanisms, also known as gates, used to control the flow of information passing from one state to another are represented by the colored circles. For each gate, additional parameters need to be learned by the network in a way that the element-wise application of the sigmoid (red circles) or tanh (purple circles) functions can wisely retain what should be preserved and what can be forgotten in a long sequence.

Despite best efforts, these architectures are still vulnerable to overfitting, compromising their ability to generalize well to new data. Potential solutions to prevent this phenomenon include regularization techniques such as L2 penalty, early stopping and dropout. In this work, a special type of dropout proposed by \citet{Kingmaetal2015} is used in combination with a GRU architecture is considered to perform the comparison with the proposed network. In this Bayesian GRU, the regular dropout with continuous noise (\textit{i.e.} Gaussian dropout) is reintepreted as a variational method that allows optimal dropout rates to be inferred from the data as opposed to it being fixed and defined in advance as usual. This circumvents the need for a validation set during model selection. For more details, the interested reader is referred to \cite{Kingmaetal2015}. 

\section{Physically Recurrent Neural Network}
\label{sec:mnn}

\begin{figure}[h!]
\centering
\includegraphics[width=0.66\textwidth]{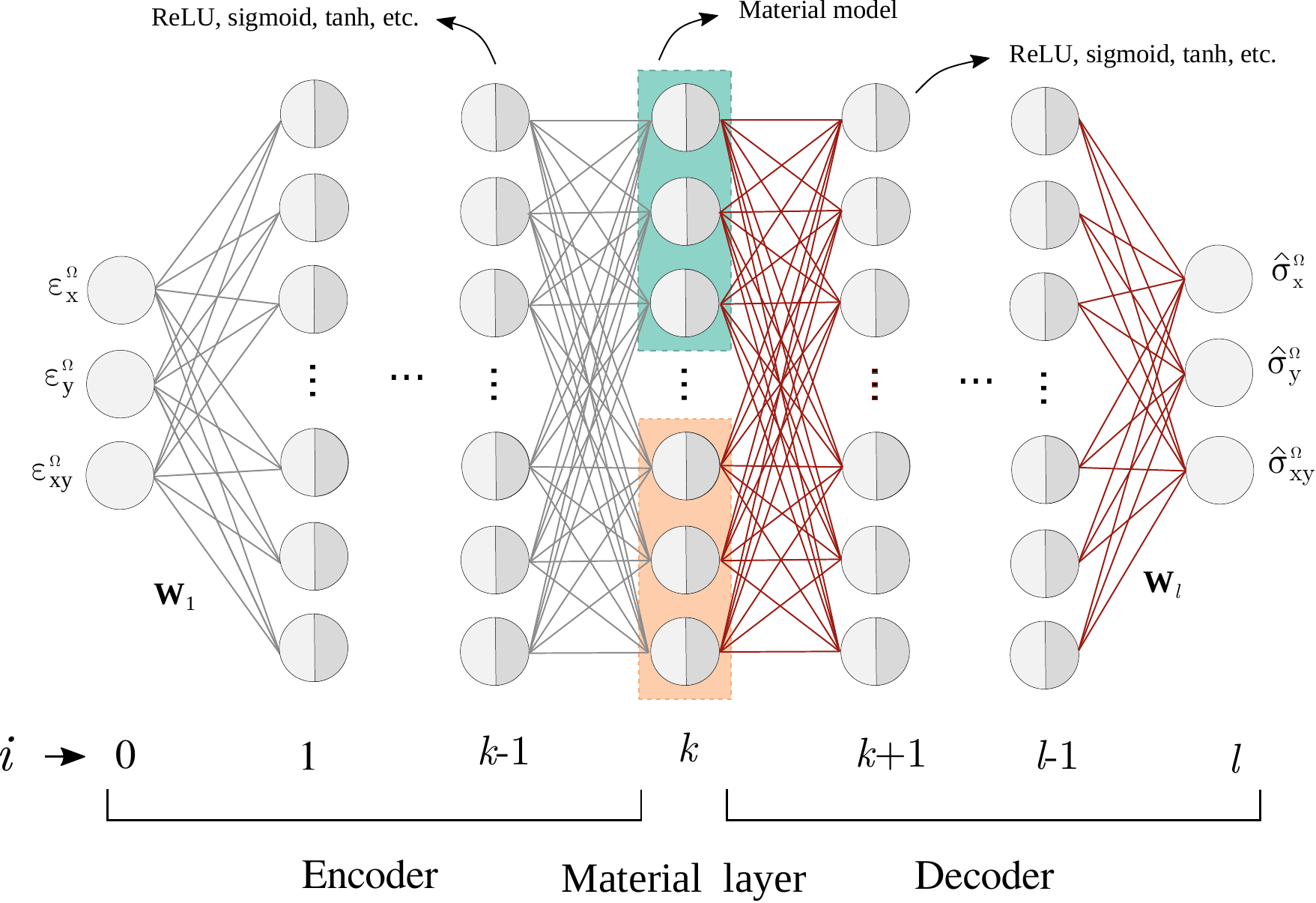}
\caption{The proposed architecture}
\label{fig:mnn}
\end{figure}

This section presents the neural network proposed to capture path-dependent behavior of heterogeneous microscopic models. The core task of the network is to learn how the macroscopic strain $\macrostrainvec$ can be dehomogenized into a small set of representative material points and how their responses can be combined to obtain the homogenized macroscopic stresses $\macrostressvec$. For this, the parametric regression model $\mathcal{R}$ illustrated in Fig. \ref{fig:mnn} is proposed: a combination of a data-driven encoder, a material layer with embedded physics-based material models and a data-driven decoder. Each of these components are discussed in detail in Sections \ref{subsec:encoder}, \ref{subsec:matlayer} and \ref{subsec:decoder}, respectively. Finally, the training process is described in Section \ref{subsec:training} and the use of this network as constitutive model in \fesqr \ frameworks is discussed in Section \ref{subsec:constmodel}. 

\subsection{Encoder}
\label{subsec:encoder}

\begin{figure}[h!]
\centering
\includegraphics[width=0.66\textwidth]{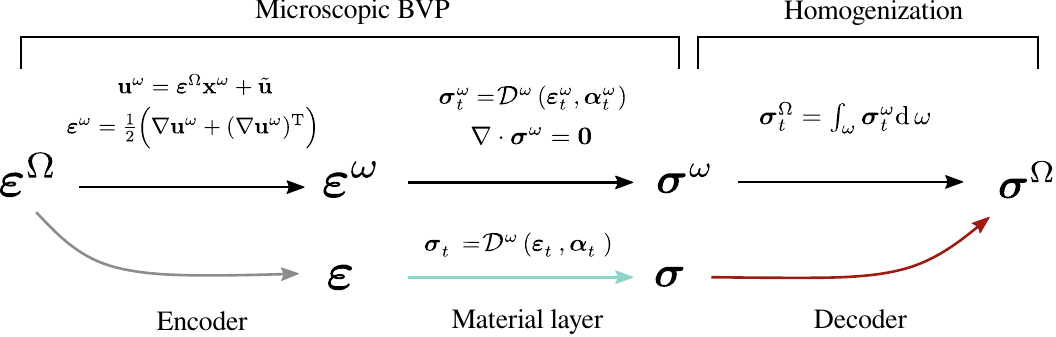}
\caption{Interpretation of the proposed network with respect to a full-order solution}
\label{fig:fe2analogies}
\end{figure}

The encoder consists of all parameters that convert the macroscopic strain from the input layer to the values used as input of the material layer, which corresponds to the grey lines in Fig. \ref{fig:mnn}. Since these values are the inputs of actual material models that are later combined by the decoder into the prediction of the macroscopic stresses, we interpret the role of the encoder as being the microscopic periodic boundary-value problem (BVP) solved with FE, only at a much lower computational cost. On the other hand, no information on the displacement field of the micromodel is retrieved by the network as the encoder is learned based exclusively on snapshots of macroscopic stresses.  
This understanding is depicted in Fig. \ref{fig:fe2analogies} by the grey curved line linking the strains from the macroscopic scale and the fictitious strains seen by the material points in the network. 

As for the architecture, an arbitrary number of layers and units (with conventional activation functions as illustrated in Fig. \ref{fig:mnn}) can be used. In case regular dense layers are employed, the neuron states ($\mathbf{a}_{i-1}$) from the previous layer $i$-1 are propagated to the following layer $i$ according to: 
\begin{equation}
    \begin{aligned}
    \mathbf{v}_i = \weights_i \mathbf{a}_{i-1} + \biases_i
    \qquad
    \mathbf{a}_i = \upphi (\mathbf{v}_i)
    \end{aligned}
\end{equation}
where $\weights_i \in \mathbb{R}^{n_i \times n_{i-1}}$ is a weight matrix and $\biases$ is a bias term with $n_i$ being the number of neurons of layer $i$, and $\upphi$ is an activation function applied in an element-wise manner to the neuron values of $i$ to introduce nonlinearity into the network. In the particular case where the dense layer is either the input or the output layers, no activation function is applied and $\mathbf{a}_i = \mathbf{v}_i$. Popular activation functions include the sigmoid, tanh and ReLU. In the present investigation, a single dense layer is considered before the material layer, which results in a linear relationship between macroscopic strains and local strains in the material layer.  
 
\subsection{Material layer}
\label{subsec:matlayer}

\begin{figure}[!ht]
\centering
\subfloat[Fictitious material point]{\label{fig:intpoint}
       \includegraphics[width=0.40\textwidth, valign = c]{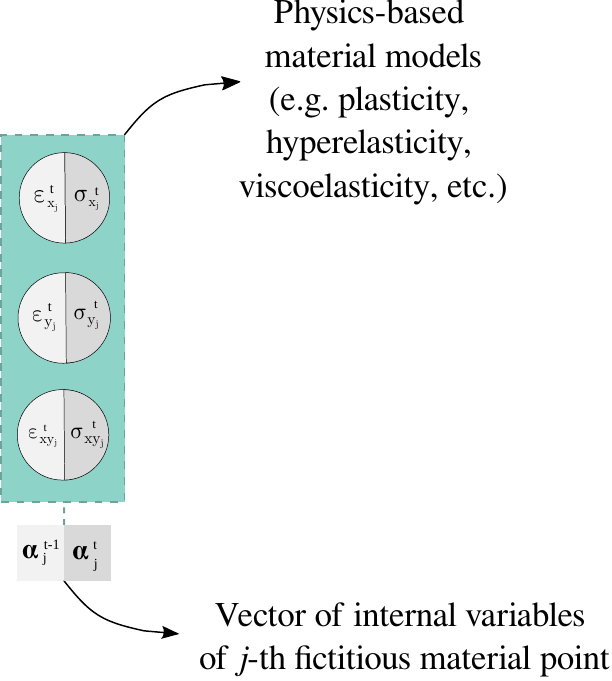}
     }
     \hfill
     \subfloat[Material cell]{\label{fig:mnncell}\includegraphics[width=0.40\textwidth, valign=c]{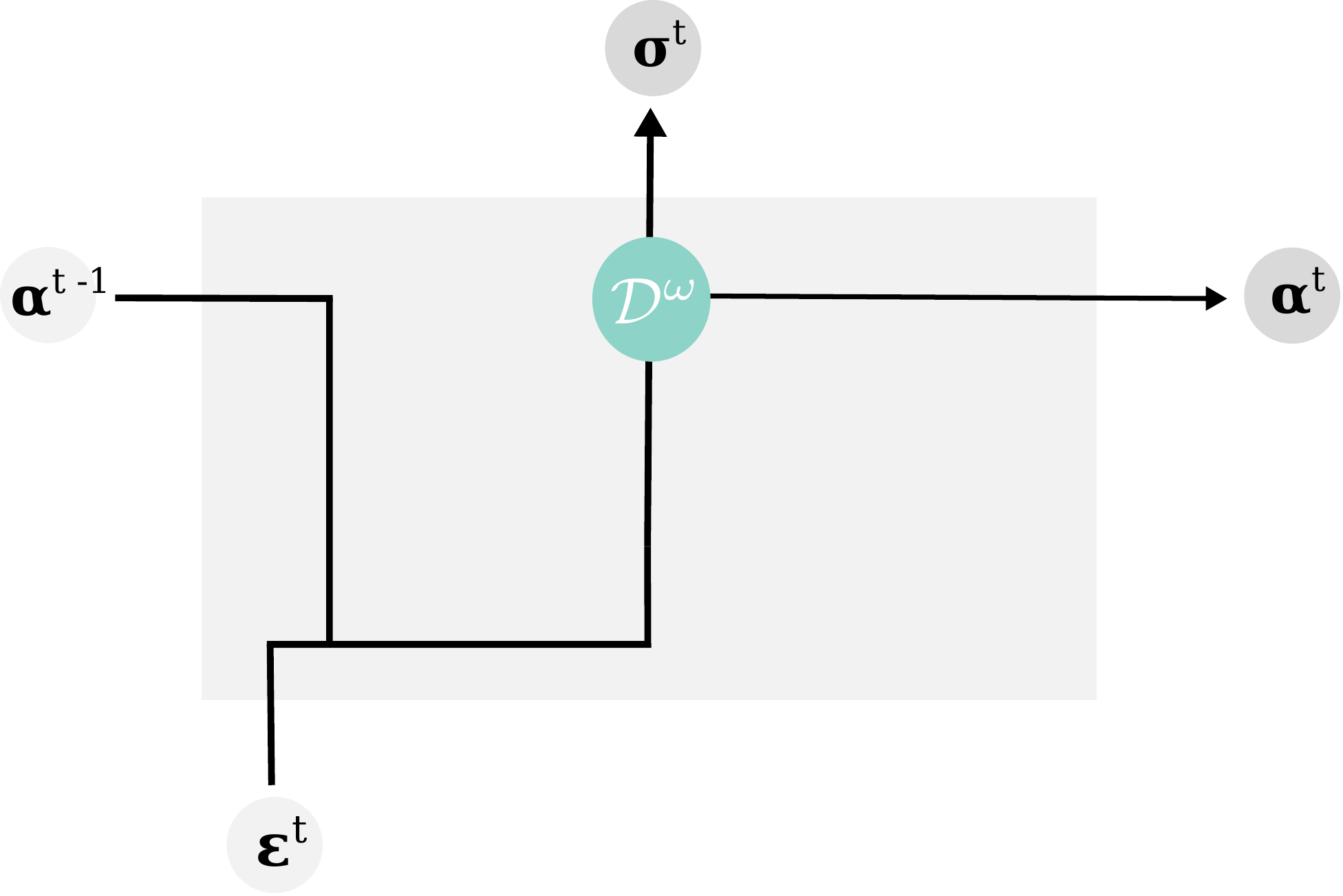}\vphantom{\includegraphics[width=0.40\textwidth,valign=c]{pictures/intpoint.pdf}}}
    \caption{Schemes of (a) fictitious material point and (b) material layer as a cell}
\label{fig:ficintpoint}
\end{figure}

The material layer is responsible for introducing explicitly the same physics-based material models used in the RVE that the network will be a surrogate for. 
To properly incorporate them and take full advantage of its outputs, important changes on how neurons are evaluated compared to regular dense layers are proposed. First, instead of introducing nonlinearity in a element-wise manner with a scalar-to-scalar activation function, neurons are grouped in $m$ sets of the size of the input layer (see colored boxes in Fig. \ref{fig:mnn}) and then evaluated as a subgroup. Each subgroup is referred to as a \textit{fictitious material point} and its size is equal to the length of the strain vector (\textit{i.e.} length 3 for the present investigation in two dimensions). In this arrangement, each neuron of the subgroup $j$ represents one component of the strain vector $\strainvec_j$, as illustrated in Fig. \ref{fig:intpoint}. 

Supposing the micromodel contains $n$ material models $\mathcal{D}^{\omega}_1$, $\dots  \mathcal{D}^{\omega}_n$, several combinations of them can be employed in the material layer. The choice on which material models should be used to evaluate the fictitious material points depends on the types of non-linearity embedded in each of these models. For instance, in this work, the micromodel consists of a composite microstructure with two material models: a linear elastic model $\mathcal{D}^{\omega}_1$ to describe the fibers and an elastoplastic model with isotropic hardening $\mathcal{D}^{\omega}_2$ to describe the matrix. However, in the material layer, only model $\mathcal{D}^{\omega}_2$ is employed in all fictitious material points. This is possible because the network still can make any of the subgroups to behave linear elastically by passing small strains and scaling the stresses to give a significant elastic contribution.  

However, if instead of a linear elastic model, a nonlinear hyperelastic model was used to describe the fibers the network would not perform optimally. In that case, although the elastoplastic model does introduce nonlinearity to the network, the (nonlinear) contribution from the fibers is no longer embedded in that model. Furthermore, if the nonlinear hyperelastic model was the one chosen to evaluate all subgroups, the network would essentially become a feed-forward one with no history information taken into account (explicitly or implicitly), losing the ability to predict elastic unloading. For such a micromodel, both material models would need to be considered.

Another interesting hypothetical case is that of a micromodel with two elastoplastic phases with different material properties. This time, depending on the contrast of the material properties, a single material model with a fixed set of properties coming from one of the two phases might be enough to reproduce the homogenized response of the micromodel. In that case, making the material properties of each fictitious material point a trainable feature would also be advantageous. This could also be a valuable feature when dealing with experimental data or with a micromodel with continuously varying material properties. Naturally, far more complex arrangements than those discussed here are found in practice. Addressing these extensions is object of ongoing research. A general guideline is to employ all different sources of nonlinearity with their respective known material properties in the material layer. 

For simplicity, we choose to illustrate the network with a general material model $\mathcal{D}^{\omega}$ for all fictitious material points. Such model can take the form of any of the material models with known material properties used in the micromodel. Here, $\mathcal{D}^{\omega}$ takes as input the current strain $\strainvec^{t} \in \mathbb{R}^{n_\varepsilon}$ and the internal variables from previous time step $\intvarvec^{t-1} \in \mathbb{R}^{n_{IntVar}}$, where $n_{IntVar}$ is the number of internal variables of the material model. With that, the model is used to evaluate current stress state $\stressvec^{t} \in \mathbb{R}^{n_{\epsilon}}$ and updated internal variables $\intvarvec^{t}$. These quantities motivated the tailor-made architecture of the proposed layer. 

To store the internal variables used as input/output of the material model, an auxiliary vector $\histunits_j \in \mathbb{R}^{n_{IntVar}}$ referred as history vector is defined. In the particular case of a subgroup with a material model with no internal variables (\textit{e.g.} linear elastic model), $\histunits_j$ does not exist since $n_{IntVar} = 0$. For the first time step, the history vector is initialized as zero for all $m$ subgroups. As information reaches the material layer and the material model is called, three outputs are made available: the stresses $\stressvec^t_j$, the updated internal variables $\intvarvec^t_j$ and the tangent stiffness matrix $\stiffmatrix^t_j \in \mathbb{R}^{n_\varepsilon\times\varepsilon}$. In this layer, only the stresses are propagated forward. To do this, each stress component is associated to a unit of the subgroup, as illustrated in Fig. \ref{fig:intpoint}. Then, the updated internal variables $\intvarvec^t$ are stored in $\histunits^t_j$ so that when new strains $\strainvec^{t+1}_j$ are fed to the fictitious material point, the material model is aware of its own history so far, making the $\strainstresscurve$ path of each subgroup unique. This architecture is illustrated in Fig. \ref{fig:mnncell}.

Note that $\histunits_j$ has no weighted connections with the neurons in the material layer (or any other layer). It is merely an object to store information about the history of the fictitious material point. It is also worth mentioning that since no data from the microscale has been collected and imposed in the network, the paths seen by the fictitious material points do not need to hold any similarity with actual integration points of the microscopic model. 

Using standard machine learning notation, the material layer propagates previous neuron states ($\mathbf{a}_{k-1}$) and applies the material model $\mathcal{D}^{\omega}$ as follows:
\begin{equation}
    \label{eq:inputtomaterial}
    \mathbf{v}_k = \weights_k \mathbf{a}_{k-1} + \biases_k \Rightarrow \mathbf{a}_k, \mathbf{h} = \mathcal{D}^{\omega} (\mathbf{v}_k, \histunits^{t-1})
\end{equation}
where $\weights_k \in \mathbb{R}^{n_k \times n_{k-1}}$ is the weight matrix connecting layers $k$-1 and $k$, $\biases_k \in \mathbb{R}^{n_k}$ is a bias term. In addition, $\mathbf{v}_k$ are the neuron values (correspond to the concatenated vector of all microscopic strains $\strainvec_j$), $\mathbf{h}^{t-1}$ and $\mathbf{h}^{t}$ are history-related term (correspond to concatenated vector of all internal variables $\intvarvec_j^{t-1}$) resulting from the material models with path-dependent behavior from past and current time step and $\mathbf{a}_k$ are the current neuron states (correspond to the concatenated vector of all microscopic stresses $\stressvec_j$). 

\subsection{Decoder}
\label{subsec:decoder}

The decoder consists of all parameters that convert the outputs from the material layer to the predicted macroscopic stress $\widehat{\mathbf{\sigma}}^{\Omega}$ in the output layer, which corresponds to the brown lines in Fig. \ref{fig:mnn}. Similar to the encoder, an arbitrary number of conventional layers and units can be employed. In the full-order solution, after convergence of the microscopic BVP, the macroscopic stresses are obtained by the volume average of microscopic stresses over the entire RVE. In the network, since the solution of the microscopic BVP is replaced by the encoder and the microscopic material points in the RVE are replaced by the few fictitious material points, the decoder is then analogous to the homogenization operator that transforms local stresses to macroscopic stresses, as illustrated by the brown curved line in Fig. \ref{fig:fe2analogies}.  

In the present work, we use a single dense layer (output) with linear activation and physics-motivated modifications to perform the task. With this, all the nonlinearity of the network arises from the models in the material layer. As discussed previously, the decoder can be understood as the averaging operator in a multiscale approach and with the chosen architecture (dense-material-dense), the weights of the output layer can be seen as the relative contribution of each fictitious material point to the macroscopic stress. Based on that, a constraint on the positivity of the weights of the output layer is considered. For that, a softplus function $\rho(\cdot)$ is applied element-wise on the weights matrix before computing the neuron values of the last layer:  
\begin{equation}
\label{eq:activateweights}
 \mathbf{v}_l = \rho (\weights_{l}) \ \mathbf{a}_{l-1} + \biases_l
\end{equation}
where $\biases_l$ is set to zero and $\mathbf{a}_{l-1}$ corresponds to the stresses coming from the material layer. This procedure guarantees that, after the tranformation, weights will always be positive. 

\subsection{Training}
\label{subsec:training}
The goal of the training phase is to minimize a loss function given by:
\begin{equation}
    \label{eq:lossfunc}
    L = \frac{1}{N} \ \sum_{i=1}^{N} \frac{1}{2} \ \| \ \macrostressvec(\boldsymbol{\varepsilon}^{\Omega}_i) - \widehat{\boldsymbol{\sigma}}^{\Omega}(\boldsymbol{\varepsilon}^{\Omega}_i) \ \|^2   
\end{equation}
where $N$ is the number of snapshots. Based on it, a Stochastic Gradient Descent (SGD) optimization algorithm is used to update the trainable parameters $\mathbf{W}$ and $\mathbf{b}$:
\begin{equation}
    \label{eq:gradients}
    \begin{aligned}
    \mathbf{W}^n = \mathbf{W}^o - \mathcal{A} \Big( \frac{1}{B} \sum_{i=1}^{B} \pd{L_i}{\mathbf{W}} \Big) \\ \qquad
    \mathbf{b}^n = \mathbf{b}^o - \mathcal{A} \Big( \frac{1}{B} \sum_{i=1}^{B} \pd{L_i}{\mathbf{b}} \Big)
    \end{aligned}
\end{equation}
where $L_i$ is the loss of the $i$-th sample, $o$ indicates current values, $n$ indicates updated values and $B$ is the size of the sample mini-batch used in the update. Finally the $\mathcal{A}$ operator depends on the solver. In this work, the Adam optimizer \cite{Adametal2014} is used. 

To compute the gradients appearing in Eq. (\ref{eq:gradients}), backpropagation in time is employed in a similar fashion as done to RNNs: based on the network state ($\mathbf{v}$ and $\mathbf{a}$) after computing each training curve with $n$ pairs of $\stressvec-\strainvec$, the chain rule is used to propagate the derivative of the loss functions starting from the output layer and progressively moving back through the network and through time. Commonly, this process is dealt with by automatic differentiation, but we present the expressions to allow for integrating the network directly into existing FE software. For this, two auxiliary quantities are defined. The first is defined for each layer and helps propagating the error through the network $\mathbf{d}_i \in \mathbb{R}^{n_i}$. Starting from the output layer $l$, it is defined as:
\begin{equation}
    \label{eq:lastlayer}
    \mathbf{d}_l = \pd{L}{\mathbf{a}_l} = \widehat{\boldsymbol{\sigma}}^{\Omega} - \macrostressvec
\end{equation}

Next, the effect of the activation function is taken into account as:
\begin{equation}
    \label{eq:lastlayerderiv}
    \bar{\mathbf{d}_i} = \mathbf{d}_i \odot \pd{\upphi(\mathbf{v}_i)}{\mathbf{v}_i} 
\end{equation}
after which it is possible to compute the gradients of the trainable parameters:
\begin{equation}
\label{eq:derivatives}
\begin{aligned}
\pd{L}{\mathbf{W}_i} = \bar{\mathbf{d}}_i \mathbf{a}_i^T \qquad
\pd{L}{\mathbf{b}_i} = \bar{\mathbf{d}}_i
\end{aligned}
\end{equation}

Finally, the values $\mathbf{d}$ of the previous layer (the next layer to be backpropagated) can be computed as:
\begin{equation}
    \label{eq:backnextlayer}
    \mathbf{d}_{i-1} = \mathbf{W}_i^T \bar{\mathbf{d}}_i 
\end{equation}
and the algorithm moves to Eq. (\ref{eq:lastlayerderiv}) for layer $i$-1.

When reaching the material layer, recall that the internal variables are stored in $\histunits$ and used for keeping track of the evolution of the internal variable through time. For that reason, a second auxiliary quantity is introduced and Eq. (\ref{eq:lastlayerderiv}) is replaced by:
\begin{equation}
    \label{eq:derivmatlayer}
    \bar{\mathbf{d}_i} = \mathbf{d}_i \odot \pd{\mathbf{a}_i}{\mathbf{v}_i} + \mathbf{d}_h^{t+1} \odot \pd{\mathbf{h}}{\mathbf{v}_i} 
\end{equation}
where the first term concerns the derivatives of stresses with respect to strains, the second term concerns the derivatives of the current internal variables with respect to strains and $\mathbf{d}_h \in \mathbb{R}^{n_i}$ is given by:
\begin{equation}
    \label{eq:intvarvecder}
    \mathbf{d}_h = \mathbf{d}_i \odot \pd{\mathbf{a}_i}{\mathbf{h}}  +  \mathbf{d}_h^{t+1} \odot \pd{\mathbf{h}}{\mathbf{h}^{t-1}}
\end{equation}

Note that the derivatives of the stresses with respect to the strains of material point $j$ are an output of the material model: the tangent stiffness matrix $\mathbf{D}_j$.  The remaining derivatives in Eqs. (\ref{eq:derivmatlayer}) and (\ref{eq:intvarvecder}) are evaluated using central finite differences. Naturally, computing gradients with other methods would also be possible. For instance, if the material model used in the network supports automatic differentation, storing the internal variables in $\mathbf{h}$ for backpropagation can be bypassed as the derivatives are automatically obtained in this approach. 

Finally, to obtain the gradients of the trainable parameters including the history-dependence coming from the material layer and compute the values $\mathbf{d}$ of the previous layer, we consider Eq. (\ref{eq:derivmatlayer}) instead of (\ref{eq:lastlayerderiv}) in the expressions shown in Eqs. (\ref{eq:derivatives}) and (\ref{eq:backnextlayer}), respectively.

\subsection{Use as constitutive model}
\label{subsec:constmodel}
To make new stress predictions, the macroscopic strain $\macrostrainvec$ is fed to the input layer and a complete forward pass is performed. The final activated neuron values of the output layer give the predicted stress. To obtain the macroscopic consistent tangent stiffness matrix $\mathbf{D}^{\Omega}$, a complete backward pass is required:
\begin{equation}
    \label{eq:stiffmatrix}
    \mathbf{D}^{\Omega} = \pd{\widehat{\stressvec}^{\Omega}}{\strainvec^{\Omega}} = \pd{\mathbf{a}_l}{\mathbf{v}_0} = \mathbf{J}
\end{equation}
which is obtained with a backward pass through the network:
\begin{equation}
    \label{eq:tangent}
    \mathbf{J}_i = \mathbf{J}_{i+1} \mathbf{I}_i^{\upphi} \weights_i \ \ \ \textrm{with} \ \ \ \mathbf{J}_{l+1} = \mathbf{I}
\end{equation}
where $\mathbf{I}_i^{\upphi}$ is a matrix whose diagonal contains the derivatives of the activation function with respect to the neuron values $\mathbf{v}$:
\begin{equation}
    \label{eq:diagmatrix}
    \mathbf{I}_i^{\upphi} = \textrm{diag} ( \pd{\upphi(\mathbf{v}_i)}{v} )
\end{equation}
except for the material layer. In that case, such matrix is full and consists of the concatenation of the tangent stiffness matrix of all fictitious material points. It is worth mentioning that despite the linear dependency on the tangent stiffness matrices of the material models, the Jacobian matrix of the network does not inherit their spectral properties.

\subsection{Analogies to other methods}
\label{subsec:analogies}

In this section, the parallels between features of the proposed network and related works in the literature are briefly discussed. One possible analogy comes from hyper-reduced-order models \cite{HERNANDEZ2017}. With the architecture chosen for the present investigation, both methods work on a reduced number of material points with modified (integration) weights. However, in the network, these points are only fictitious and learned by the encoder based on snapshots of the homogenized stresses. Moreover, each stress component is associated with a different weight. By contrast, the material points in the hyper-reduction approach exist in the microscopic model and a single modified integration weight of each material point selected is used to compute all its stress/internal force components.

Following the discussion on the encoder, it is worth highlighting how this feature would be framed with respect to asymptotic homogenization schemes such as Mori-Tanaka \cite{MORI1973}. In this type of solution, the microscopic problem is also not solved explictly and only average fields are calculated. Relying on the equivalent inclusion idea and on Eshelby's solution \cite{Eshelby1957}, the strain concentration tensor is obtained analytically and yields the full solution of the microcopisc model as it correlates the average field of the phases in the micromodel with its average field. In our network, although the macroscopic stresses are also obtained by relating macroscopic and (fictitious) microscopic strains through an encoder, here no average field is calculated for each of the phases. Indeed, not every phase needs to be included in the material layer and multiple strain paths for the same phase are considered. Furthermore, while Mori-Tanaka is accurate for moderate volume fractions of the inclusions, such restriction is not present in our method. 

Compared to PINNs, in which physical constraints are explicitly included in the loss function, here, most physical constraints are naturally taken care of by the physics-based material models directly embedded in the material layer. The proposed approach is also more general as it can be directly used for arbitrary material models and is not particularly tailored to a single type of model (\textit{e.g.} elastoplastic behavior \cite{HAGHIGHAT2021, Eghbalianetal2022}). As an added benefit, our model selection procedure makes physical sense: we add more material points or material models to the network.

Another noteworthy strategy with relevant analogies to our method is the DMN \cite{Liu2019}. In this approach, the contribution of a few material points evaluated using the classical constitutive models in the RVE is also employed to make predictions in the online phase. On the same reasoning as discussed in Section \ref{subsec:matlayer}, since the inputs come directly from actual material models, path-dependency is captured naturally. However, the main concept and architecture of DMNs are different from the ones explored here. In the offline phase, the goal of the DMN is to find a topological representation of the RVE with fewer degrees of freedom (\textit{i.e.} material points) based only on the elastic stiffness matrices of the different material phases that compose the original micromodel. For the online phase, the feature space is increased to include residual stresses of the micromodel components, and an iterative procedure is implemented. The authors compare the incremental strains of the material points at the beginning of the iteration with the one obtained by a de-homogenization process that backpropagates the macroscopic incremental strain from the output layer to the bottom layer (\textit{i.e.} input layer). Upon convergence, the set of internal variables of each material point at the bottom layer is therefore updated.

In the present work, the feature space is the same in both phases and no iterative procedure is employed in the forward pass, which simplifies implementation and reduces even further the number of material model calls. Here, the strain path each fictitious material point follows is simply described by the encoder and not all phases need to be included in the network. The homogenized stresses and tangent stiffness are obtained in a single forward and backward pass, respectively. Furthermore, the backpropagation in our approach is considerably simpler than the DMN. Although the use of homogenization (and de-homogenization) operations in the DMN assigns physical interpretation to the model, it also makes training a rather intricate process.

Finally, to draw a parallel with LSTMs, one might understand $\histunits$ as the cell state $\mathbf{c}$, but instead of using bijective and smooth functions such as the sigmoid and tanh functions to describe the evolution of the material response, the material model itself is directly employed. This bypasses the need to learn new parameters to regulate the flow of information kept or forgotten throughout time (see Fig. \ref{fig:grucell}) and has important implications for the training process. The most important one is the ability to mirror physical behaviors such as elastic unloading/reloading without ever seeing the pattern during training, a stark contrast with LSTMs and GRUs that usually require extensive training sets with multiple cycles of loading and reloading at different strain levels with different step sizes. The physical interpretation of the nonlinearity is directly embedded in the network. In the numerical examples of this work, the nonlinearity is due to plasticity, but other effects such as hyperelasticity, visco-plasticity, stiffness degradation, or any combination thereof, could be embedded by adapting the constitutive model that is used in the material layer.

\section{Design of Experiments}
\label{sec:doe}

One critical aspect of the training and testing of surrogate models is the formulation of a sampling plan. Typically, a uniform distribution of the sampling points is desirable, but that task becomes more complex when path-dependent behavior is present. In this case, pairs of strains and stresses are collected and processed as sequences, which leads to potentially infinite-dimensional parameter spaces. 

In this work, two strategies are considered. In the first approach, proportional loading paths are generated, which means that the stress ratio between the components is constant. Here, the sequence of strains is created based on two features: the loading function $\lambda (\Delta\varepsilon, t)$ and the loading direction given by the unit vector $\mathbf{n}$, where $\Delta\varepsilon$ is the step size and $t$ is the current time step. For each time step, $\mathbf{n}$ is multiplied by the scalar-valued loading function $\lambda$ creating a new set of strains, which is in turn applied at the controlling nodes of the microscopic model. 

\begin{figure}[h!]
\centering
\subfloat[Types of loading used in this study]{\label{fig:doenonproprandom}\includegraphics[width=0.45\textwidth]{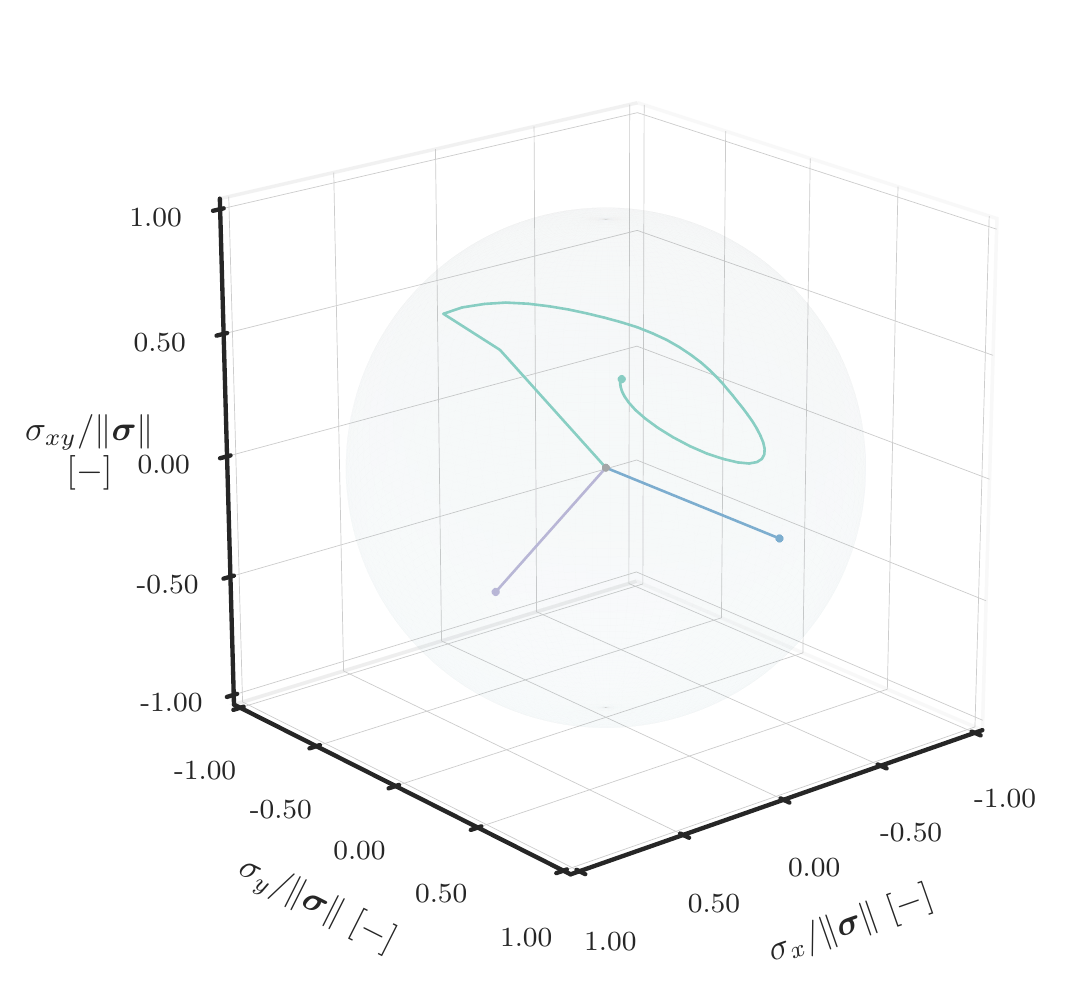}}
\hfill
\subfloat[Proportional loading with \textit{a priori} known directions]{\label{fig:doepropbasic}\includegraphics[width=0.45\textwidth]{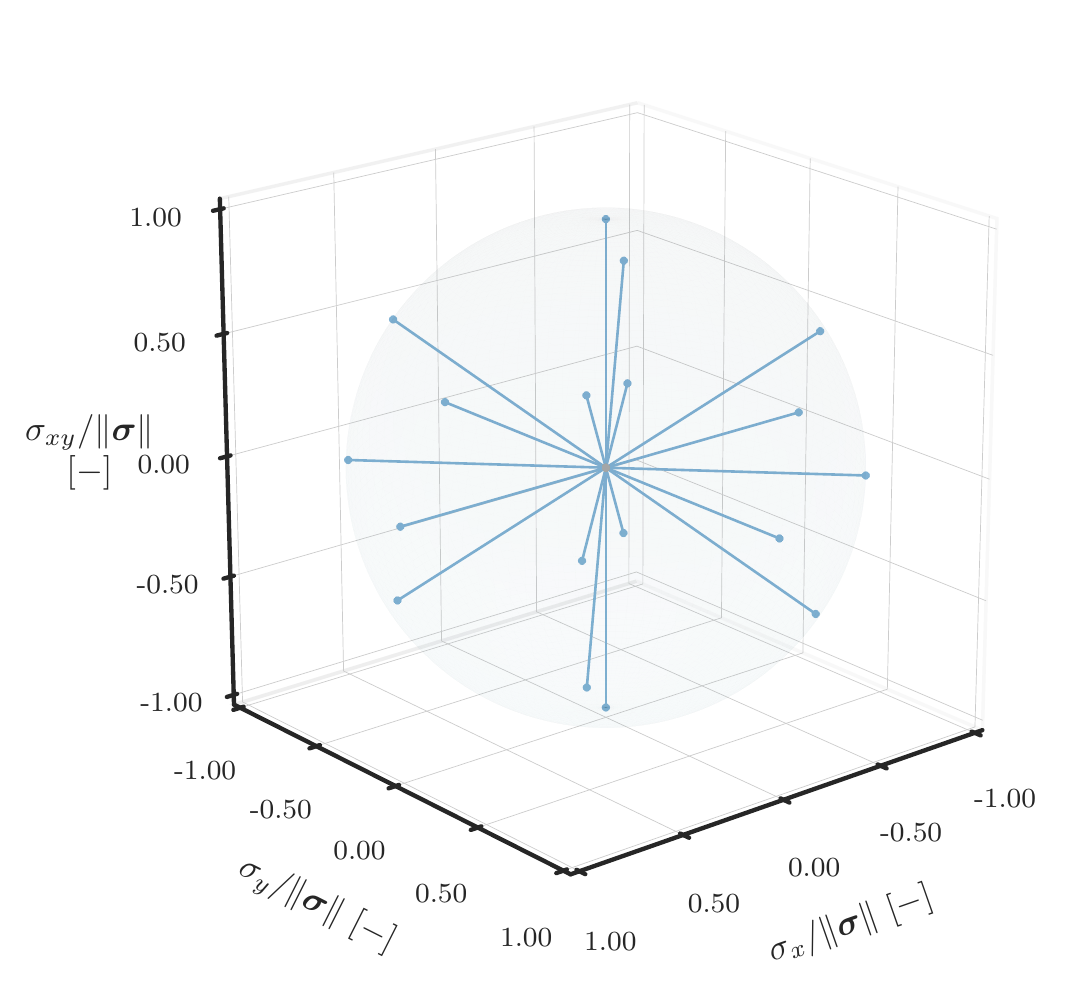}}
\caption{Different design of experiments strategies}
\label{fig:doe}
\end{figure}

\begin{figure}
\centering
\subfloat[Monotonic loading (Types I and II)]{\label{fig:doetype1}
\includegraphics[width=0.33\textwidth]{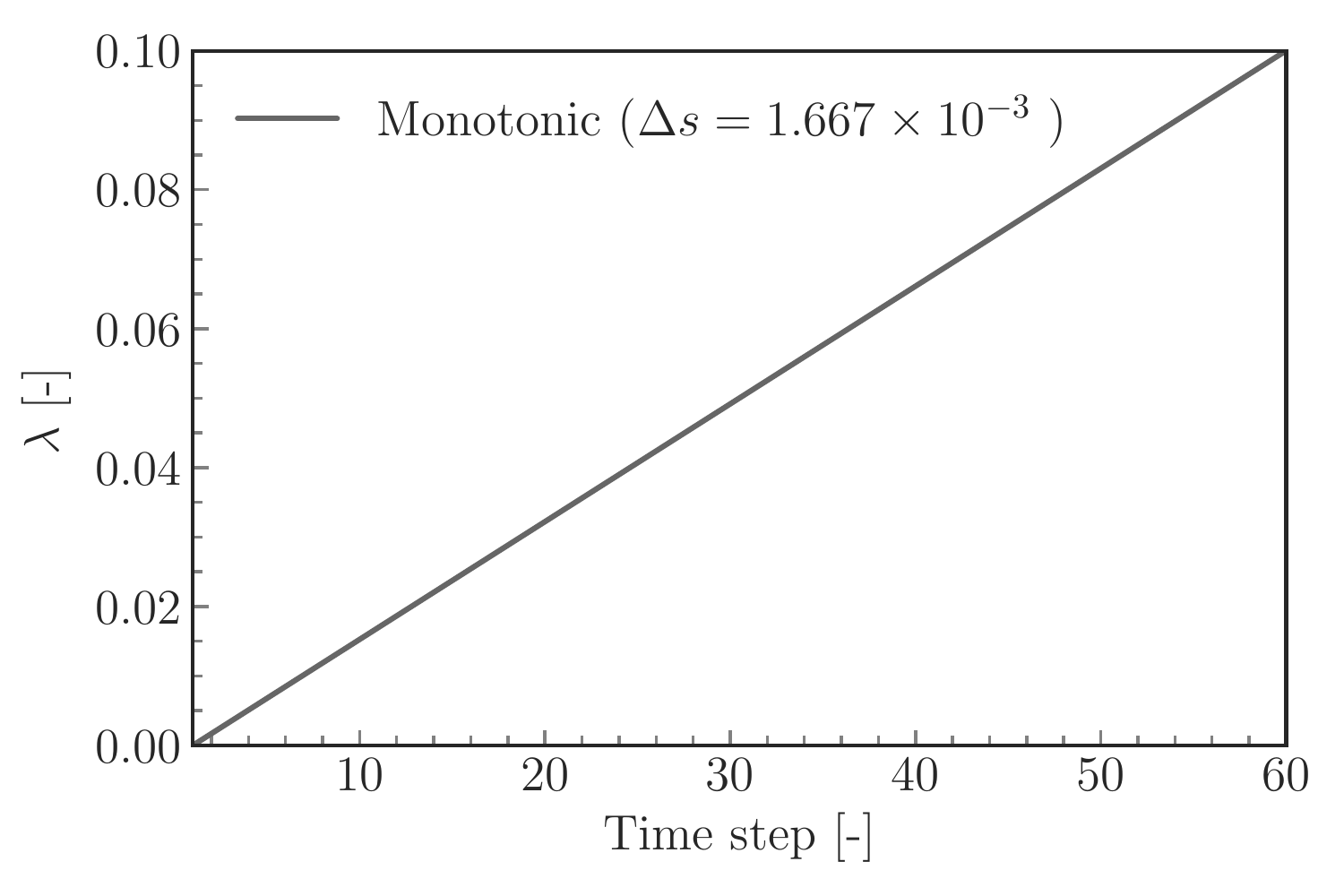}}
\subfloat[Non-monotonic loading (Type III)]{\label{fig:doetype3}
\includegraphics[width=0.33\textwidth]{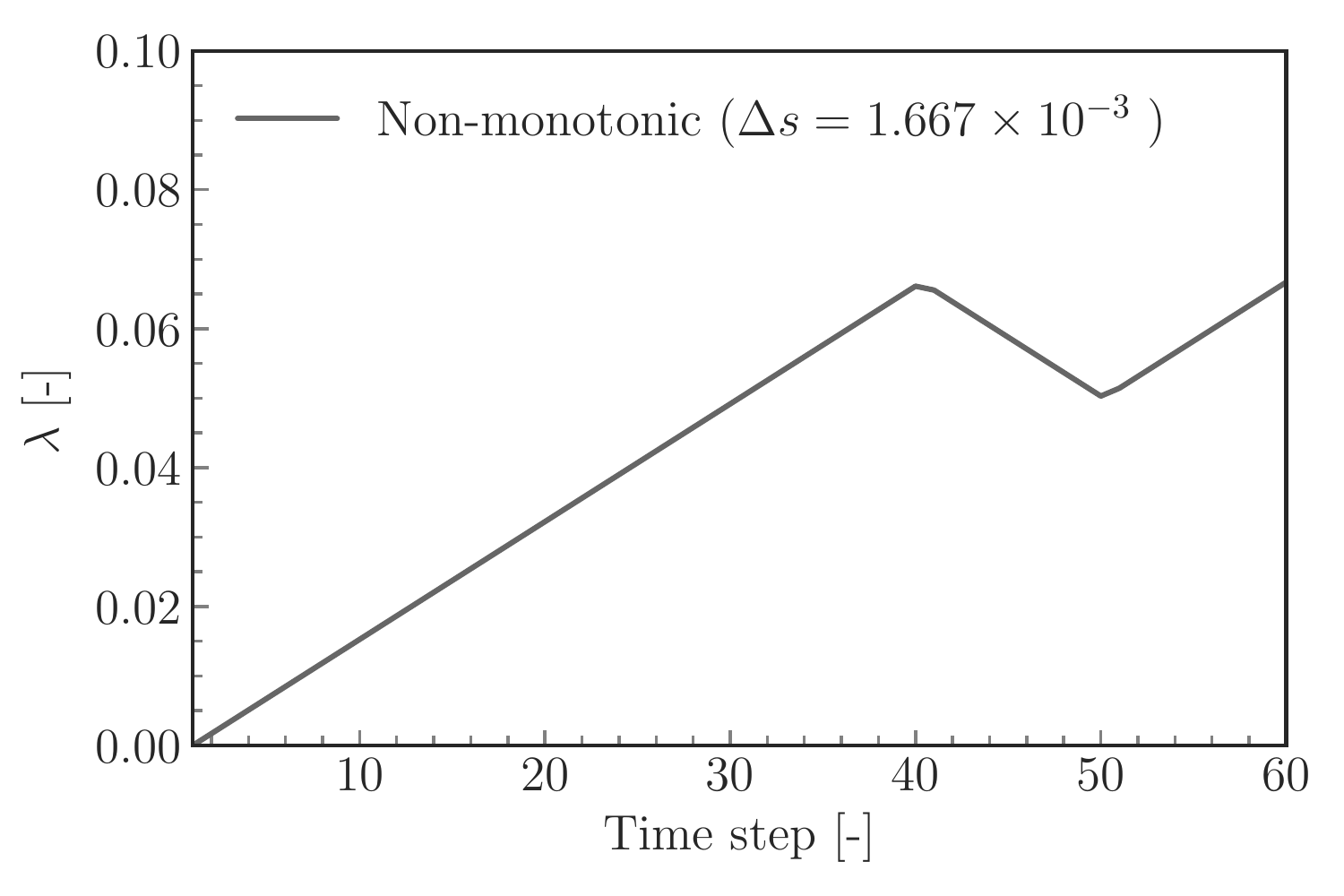}}
\subfloat[Different unloading/reloading (Type IVa)]{\label{fig:doetype3longunl}
\includegraphics[width=0.33\textwidth]{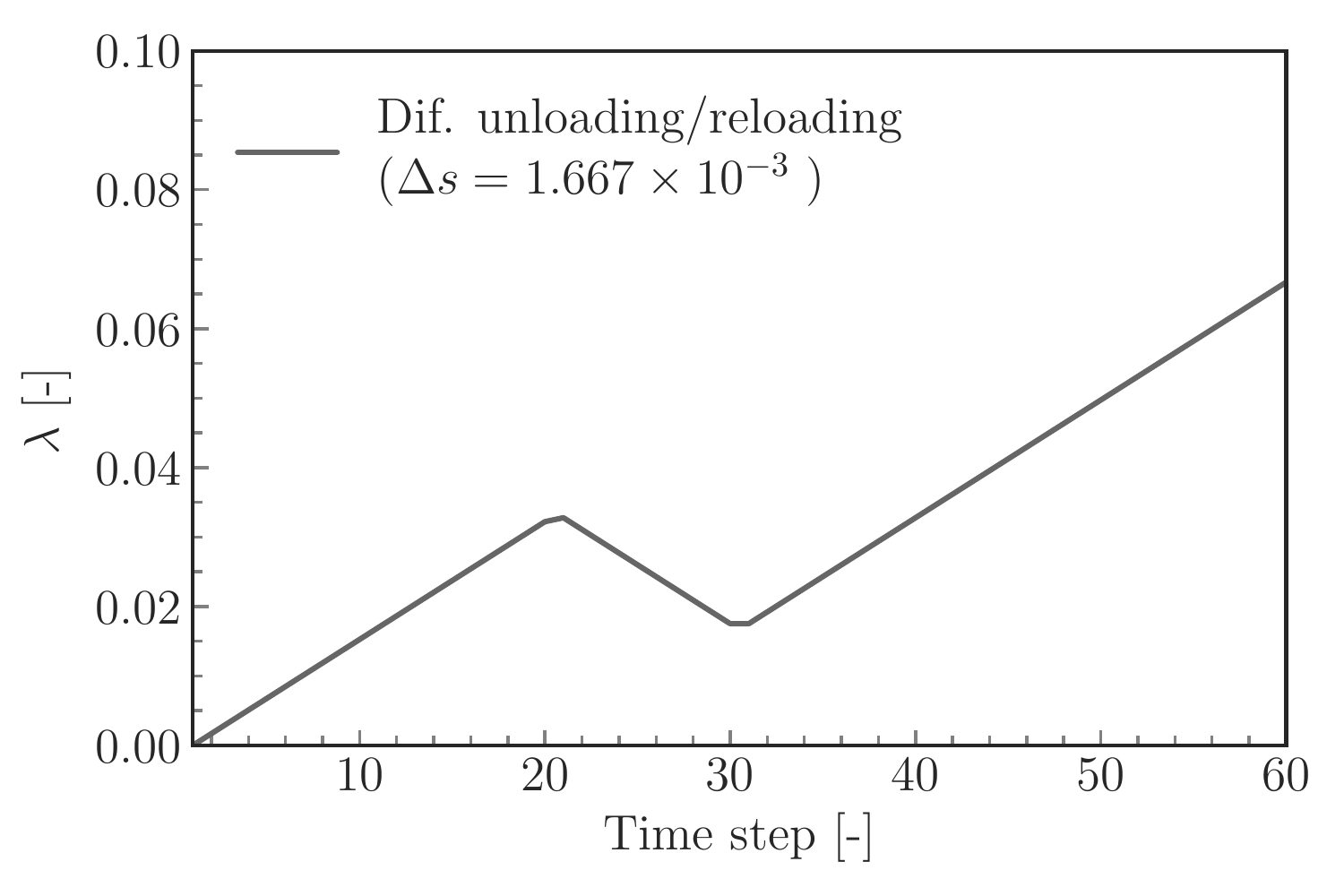}}
\hfill
\subfloat[Non-monotonic with smaller step size (Type IVb)]{\label{fig:doetype2step}
\includegraphics[width=0.33\textwidth]{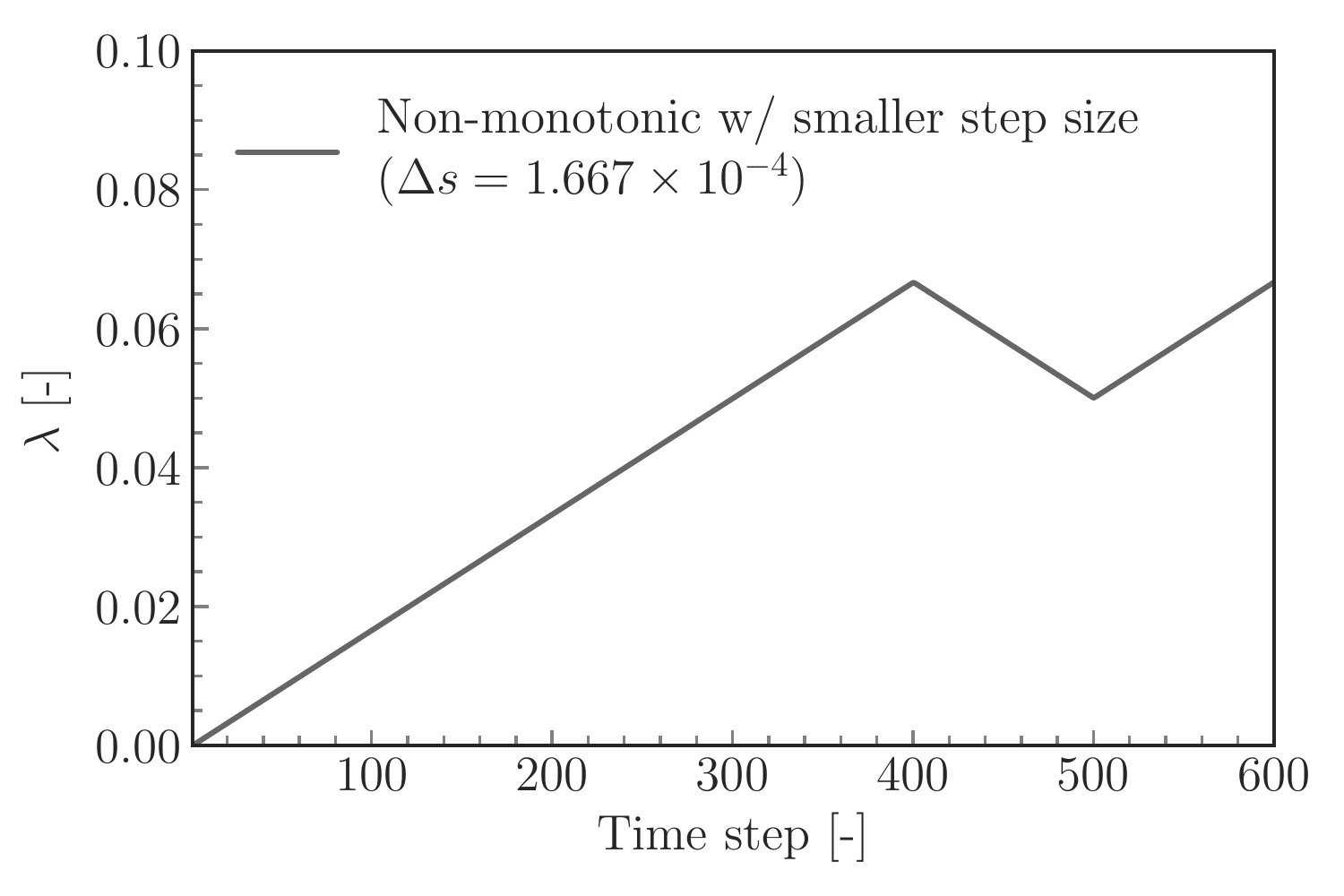}}
\subfloat[Random non-proportional and non-monotonic with $\ell = 20$ (Type V)]{\label{fig:doetype5}
\includegraphics[width=0.33\textwidth]{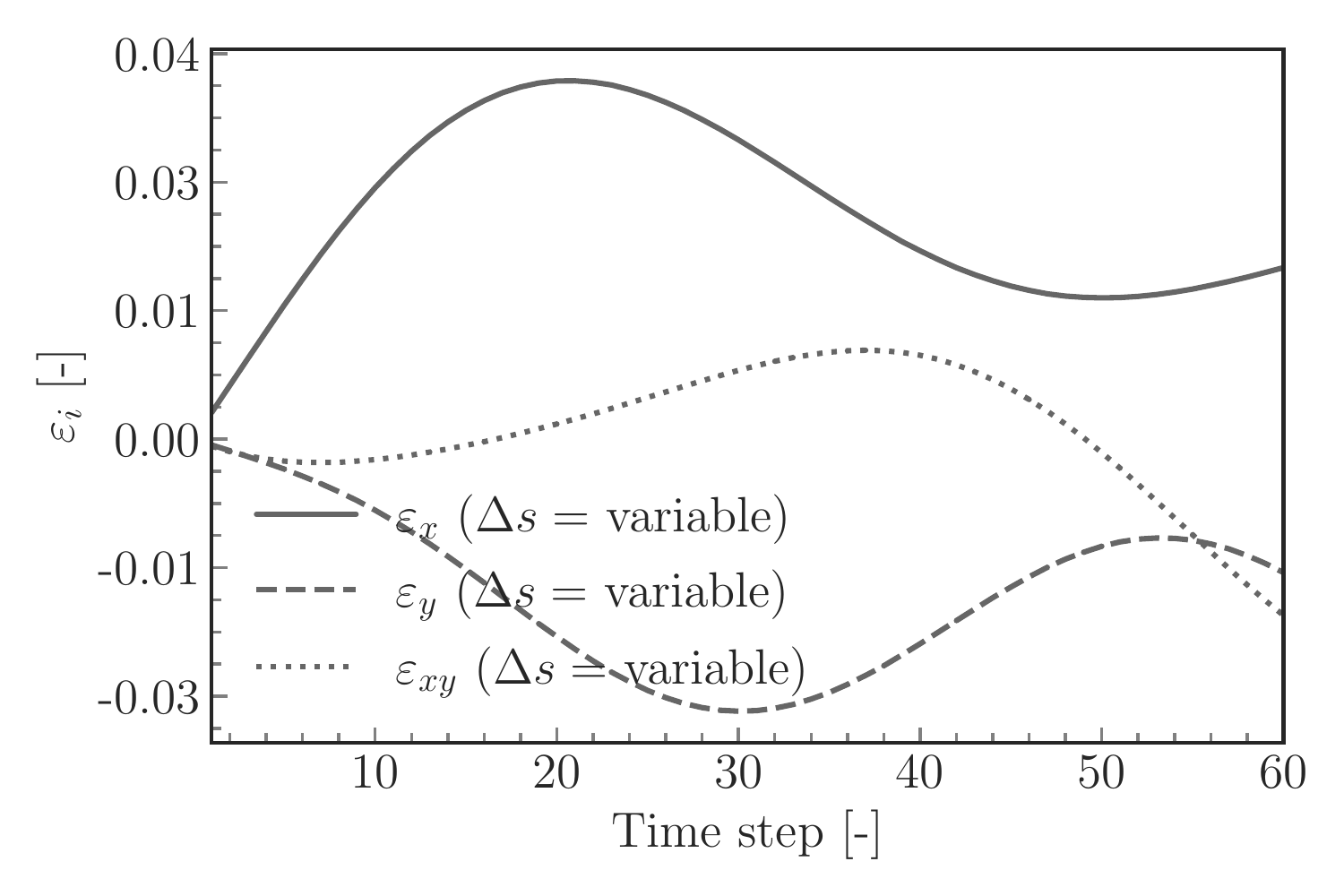}}
\caption{Proportional and non-proportional loading functions}
\label{fig:shapefunctions1and2}
\end{figure}

For monotonic loading, the loading function is as depicted in Fig. \ref{fig:doetype1}. The values in the unit vector can come from prior knowledge of the material as illustrated in Fig. \ref{fig:doepropbasic}, in which only fundamental cases such as uniaxial strain, pure shear, and biaxial cases are considered, or from random distributions as represented by the purple line Fig. \ref{fig:doenonproprandom}. In the present work, the random directions are obtained by sampling values from $n_\varepsilon$ independent Gaussian distributions (X $\sim \mathcal{N}(0, 1)$) and subsequently normalizing the vectors.

Despite the simplicity in creating such paths, RNNs trained exclusively on monotonic cannot predict cylic responses. Thus, to create non-monotonic sequences, a linear piecewise function as the one depicted in Fig \ref{fig:doetype3} is used. Note that even though the loading function is changed, the unit vector is kept constant for the entire strain sequence, yielding proportional loading. However, to cover the entire space of possible cyclic responses, a large (and \textit{a priori} unknown) number of curves comprehending different unloading points with different duration of unloading/reloading and step sizes is necessary. In this work, that matter is first handled in a simplified way by only sampling two different cycles of unloading/reloading. 

Finally, in a more general approach, a second strategy to create the $\strainstresscurve$ paths is considered: the \textit{random walks}. These are typically defined by sampling random strain increments with random loading directions for each time step, resulting in non-proportional loading. In this work, random walks are created by associating the prescribed strains to independent Gaussian Processes (GPs) with X $\sim \mathcal{N}(\mu,\,\sigma^{2})$ and covariance function given by: 
\begin{equation}
\label{eq:gpcurves}
k(\xvec_p, \xvec_q) = \sigma_f^2 \exp \Big( - \frac{1}{2 \ell^2} \|  \xvec_p - \xvec_q \|^2 \Big)
\end{equation}
where $\sigma_f^2$ is the variance and $\ell$ is a length scale. In this setting, the lengthscale controls the smoothness of the strain path and the variance controls how large the step size can be for each prescribed degree of freedom of the controlling nodes. Similar approaches were employed by \citet{Mozaffar26414} and \citet{LOGARZO2021}. 
\begin{algorithm}
    \SetKwInOut{Input}{Input}
    \SetKwInOut{Output}{Output}
    \Input{lengthscale $\ell$, variance $\sigma_f^2$, number of strain components $N_{\textrm{components}}$, number of time steps $N_{\textrm{steps}}$}
    \Output{macroscopic strains $\mathfrak{D}_{\varepsilon}$ and macroscopic stresses $\mathfrak{D}_{\sigma}$}
    initialize datasets: $\mathfrak{D}_{\varepsilon} \leftarrow \varnothing$, $\mathfrak{D}_{\sigma}  \leftarrow \varnothing$ \\
    \For{$i \in [1, 2, ..., N_{\text{components}}]$}{
    initialize input and output datasets for $\textrm{GP}_i$: $\mathbf{X}_{\textrm{GP}_i} \leftarrow  0 $, $\mathbf{Y}_{\textrm{GP}_i} \leftarrow 0$ \\ 
    initialize $\textrm{GP}_i$: $\textrm{GP}_i \leftarrow $  \texttt{initGP}( $\mathbf{X}_{\textrm{GP}_i}, \mathbf{Y}_{\textrm{GP}_i}, \ell, \sigma_f^2$) \\    
    }   
    \For{$t \in [1, 2, ..., N_{\text{steps}}]$}{    
        initialize current macroscopic strain: $\boldsymbol{\varepsilon}_{\textrm{current}} \leftarrow \varnothing$ \\
    \For{$i \in [1, 2, ..., N_{\text{components}}]$}{
    sample from posterior distribution: $\varepsilon_i \leftarrow$ GP$_{i}$::\texttt{samplePosterior} (  $t$ ) \\
    add value to strain vector: $\boldsymbol{\varepsilon}_{\textrm{current}} \leftarrow \boldsymbol{\varepsilon}_{\textrm{current}} \cup \varepsilon_i$
    }
    solve micromechanical BVP: $\boldsymbol{\sigma}_{\textrm{current}} \leftarrow$ \texttt{fullModel::materialUpdate} ( $\boldsymbol{\varepsilon}_{\textrm{current}}$ ) \\
    \If{convergence}
    {
      store equilibrium solution of micromodel: \texttt{fullModel::storeSolution( )} \\ 
      add macroscopic strains and stresses to dataset: $\mathfrak{D}_{\varepsilon} \leftarrow \mathfrak{D}_{\varepsilon} \cup \boldsymbol{\varepsilon_{\textrm{current}}}$, $\mathfrak{D}_{\sigma} \leftarrow \mathfrak{D}_{\sigma} \cup \boldsymbol{\sigma}_{\textrm{current}}$ \\
    \For{$i \in [1, 2, ..., N_{\text{components}}]$}{   
      {
      add time step and current strain to dataset of GP$_i$: $\mathbf{X}_{\textrm{GP}_i} \leftarrow \mathbf{X}_{\textrm{GP}_i} \cup t$, $\mathbf{Y}_{\mathrm{GP}_i} \leftarrow \mathbf{Y}_{\textrm{GP}_i} \cup \varepsilon_i$ \\
      update GP$_i$ with new data: GP$_i$::\texttt{update} ( $\mathbf{X}_{\textrm{GP}_i}$, $\mathbf{Y}_{\textrm{GP}_i}$  ) 
      }
      }
   }
   }
   \textbf{return} ($\mathfrak{D}_{\varepsilon}$, $\mathfrak{D}_{\sigma}$)
   \caption{Generation of random loading path using GPs}
   \label{alg:samplegps}
\end{algorithm}

The details of the present implementation are given in Algorithm $\ref{alg:samplegps}$. Note that instead of drawing the entire strain sequence for a given component, we sample it step by step and update the GP dataset before sampling again. This strategy results in the same strain sequence given a fixed random seed throughout the steps, but in this way the GPs can also be used in applications where the number of loading steps is changed on-the-fly. Following the work of \citet{LOGARZO2021}, we define the mean of all GPs to be zero and include $t=0$ and $\varepsilon_i = 0$ as a prior. In addition to that, references to \textit{fullModel} (\textit{i.e.} the full-order microscopic model) in Algorithm \ref{alg:samplegps} are kept as minimal and general as possible. One example of loading path resulting from Algorithm \ref{alg:samplegps} is illustrated in Fig. \ref{fig:doenonproprandom}. 

In this paper, both strategies generate $\strainstresscurve$ curves containing 60 time steps, unless stated otherwise. To summarize the types of loading studied in the following sections:
\begin{itemize}
\item Type I: monotononic and proportional loading paths with \textit{a priori} known directions. The 18 directions used to train the proposed network are illustrated in Fig. \ref{fig:doepropbasic} and include uniaxial strains, pure shear, biaxial cases and biaxial with shear cases. 
\item Type II: monotonic and proportional loading paths randomly spread across the design space. The loading directions are generated randomly and the loading function is as shown in Fig. \ref{fig:doetype1}. 
\item Type III: non-monotonic and proportional loading paths randomly spread across the design space. Again, the loading directions are random, but the loading function is now given by Fig. \ref{fig:doetype3} and includes one cycle of unloading;
\item Type IV: Variations to Type III:
\begin{itemize}
\item Type IVa: same loading directions as the test set of Type III, but unloading/reloading takes place at a different point in time as shown in Fig. \ref{fig:doetype3longunl};
\item Type IVb: same loading directions as the test set of Type III, but time step is 10 $\times$ smaller. Thus, to reach the same norm as the original curve in Type II, 600 time steps are evaluated, as depicted in Fig. \ref{fig:doetype2step};
\end{itemize}
\item Type V: non-monotonic and non-proportional loading paths randomly spread across the design space. A GP-based path described by Eq. (\ref{eq:gpcurves}) is illustrated in Fig. \ref{fig:doenonproprandom}. Fig. \ref{fig:doetype5} ilustrates the strain paths of each component using this approach with lengthscale $\ell = 20$ and $\sigma_f = 1.0 \times 10^{-3}$. 
\end{itemize}

\section{Assessing network performance}
\label{sec:results}

In this section, the performance of the proposed network is compared to a state-of-the-art RNN trained on different training dataset sizes and methods to sample the design space. The comparison is done for a single micromodel. Specifically, four scenarios are investigated: (i) predicting unloading/reloading behavior from monotonic data, (ii) predicting unloading/reloading behavior from non-monotonic data, (iii) predicting unseen patterns from non-monotonic data, and (iv) training with non-monotonic and non-proportional loading paths. In the first three scenarios, our network is trained exclusively on the fundamental loading cases of Type I (18 curves), while the training of the RNN is an open question to be addressed in the following sections. 

From now on, the network presented in this work will be referred to as Physically Recurrent Neural Network, or simply PRNN. The PRNN was trained for 80000 epochs, while the RNN was trained for 60000 epochs with an early stopping criterion, which consists of interrupting training if the best training loss so far is not improved over a given period (in this work, 5000 epochs). The Adam optimizer is used in all cases with batch size of 9 and default parameters suggested by \citet{Adametal2014}, with the exception of the learning rate of 0.01 for the RNNs. The layer sizes are chosen through model selection to provide optimal results and fair comparison with our approach to the best of our knowledge. The methodology adopted is briefly described in Section \ref{subsec:modelsel}. 

\begin{figure}[h!]
\centering
\subfloat{\includegraphics[width=0.25\textwidth]{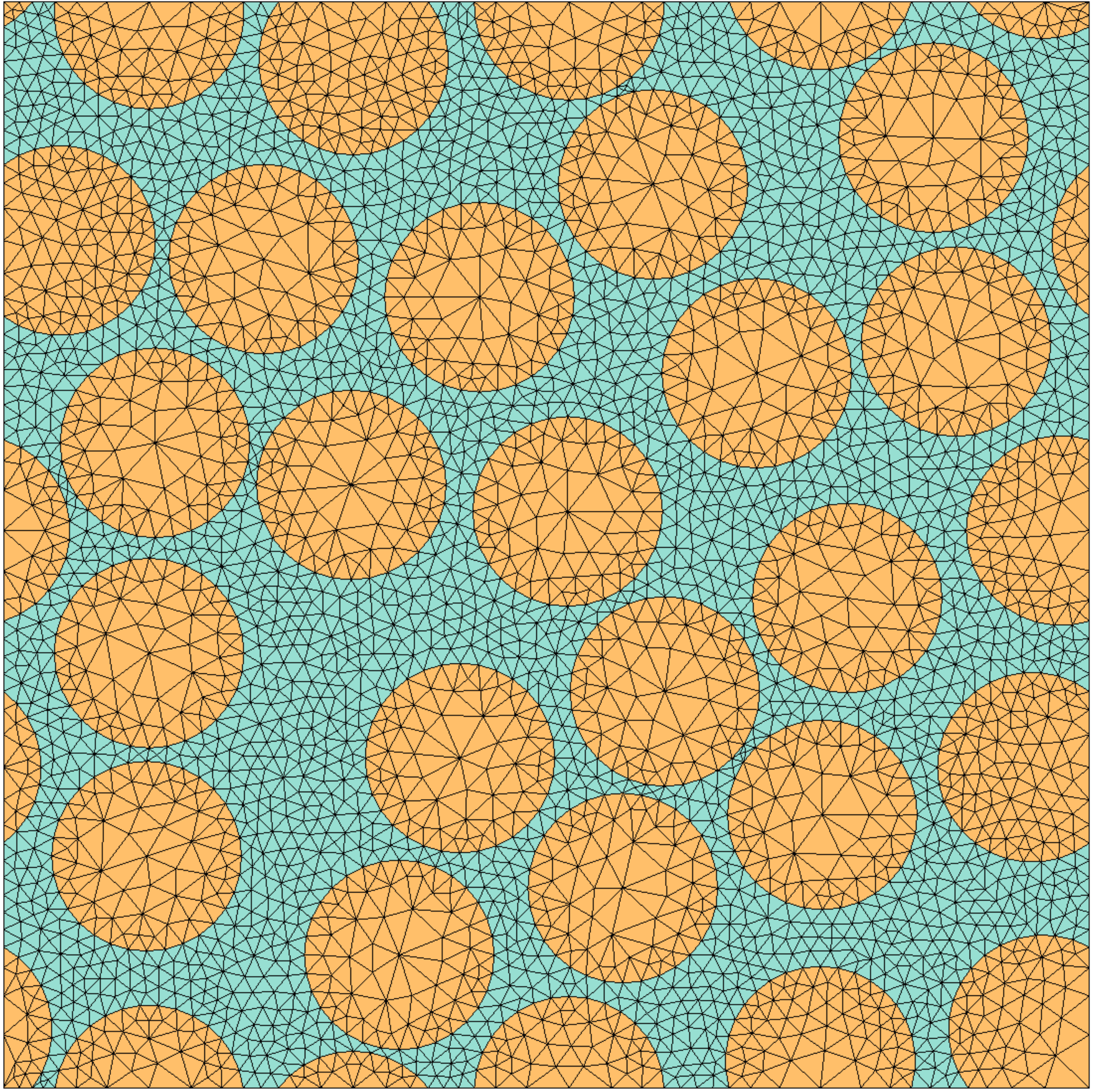}}
\caption{Geometry and mesh discretization of microscopic model adopted in this work.}
\label{fig:rve9}
\end{figure}

The microscopic model consists of an RVE with 36 elastic fibers (volume fraction = 0.6) with properties E = 74000 MPa and $\nu$ = 0.2 embedded in an elastoplastic matrix with isotropic hardening. The geometry and the mesh with 7048 elements are depicted in Fig. \ref{fig:rve9}. The elastoplastic matrix is modeled using the von Mises yield criterion with properties E = 3130 MPa, $\nu$ = 0.3 and yield stress given by:
\begin{equation}
\sigma_y = 64.8 - 33.6 \exp^{-\varepsilon^{p}_{eq}/0.0003407}
\end{equation}
where $\varepsilon^{p}_{eq}$ is the equivalent plastic strain defined as:
\begin{equation}
\label{eq:epeq}
\varepsilon^{p}_{eq} = \sqrt{\frac{2}{3} \ \bm{\varepsilon}^p : \bm{\varepsilon}^p}
\end{equation}
and $\bm{\varepsilon}^p$ is the plastic strain. Plane stress conditions are assumed.

\subsection{Model selection}
\label{subsec:modelsel}

In this section, a grid-search strategy is employed to choose the best architecture for the networks. The goal is to find the optimum layer size before heading to the testing sections. Here, four different layer sizes and 10 different weight initializations are considered to mitigate the effect of randomness. 

\begin{figure}[ht!]
\centering
\subfloat[Validation error for PRNNs trained on $\mathfrak{D}_{\mathrm{PRNN}} = \{$18 Type I curves$\}$ and $\mathcal{V}_{\mathrm{PRNN}} = \{$54 Type II curves$\}$]{\label{fig:paramstudymnn}\includegraphics[width=0.48\textwidth]{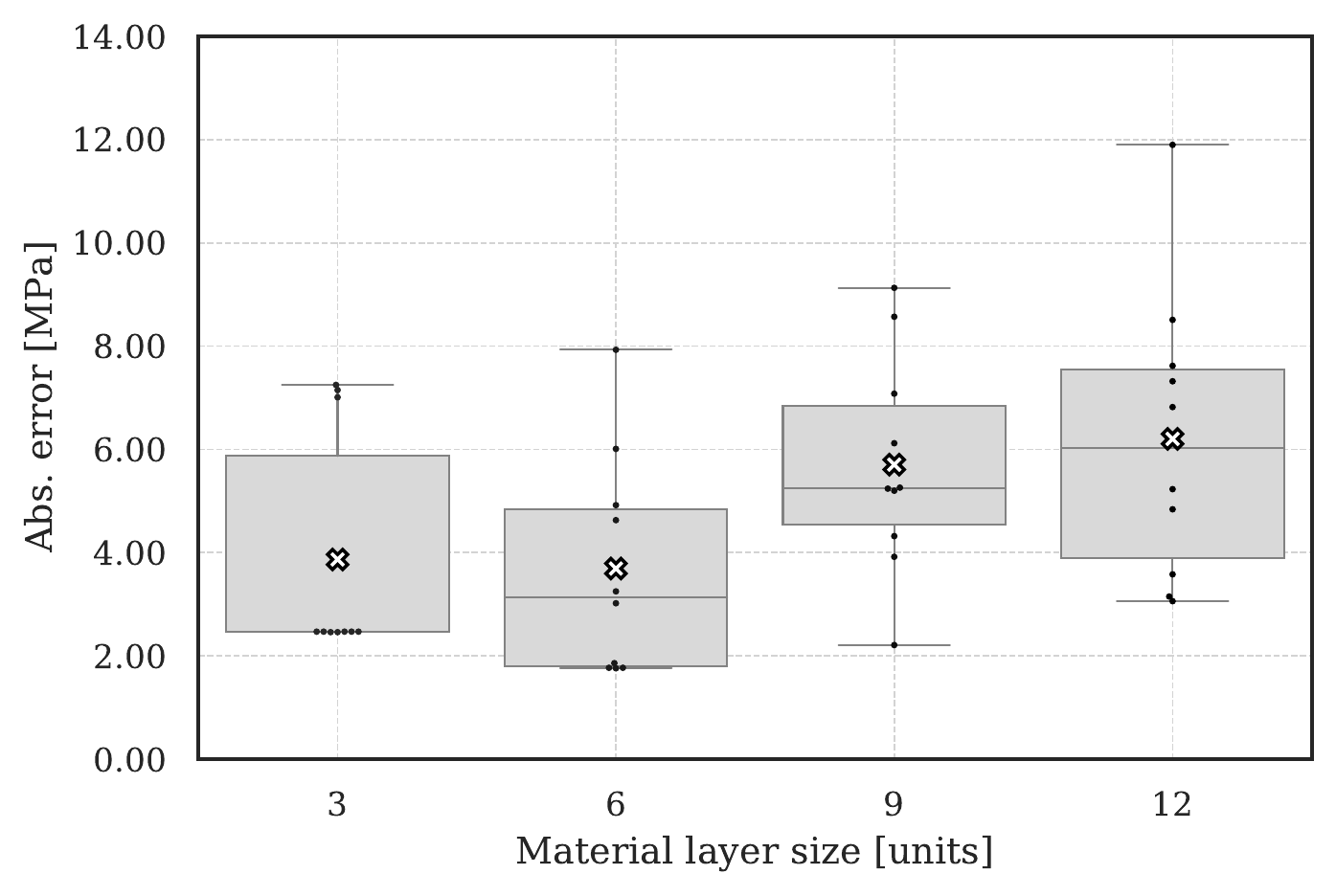}}
\hfill
\subfloat[Training error for RNNs trained on $\mathfrak{D}_{\mathrm{RNN}} = \{$18 Type I curves, 90 Type II and 90 Type III curves$\}$]{\label{fig:paramstudyrnn}\includegraphics[width=0.48\textwidth]{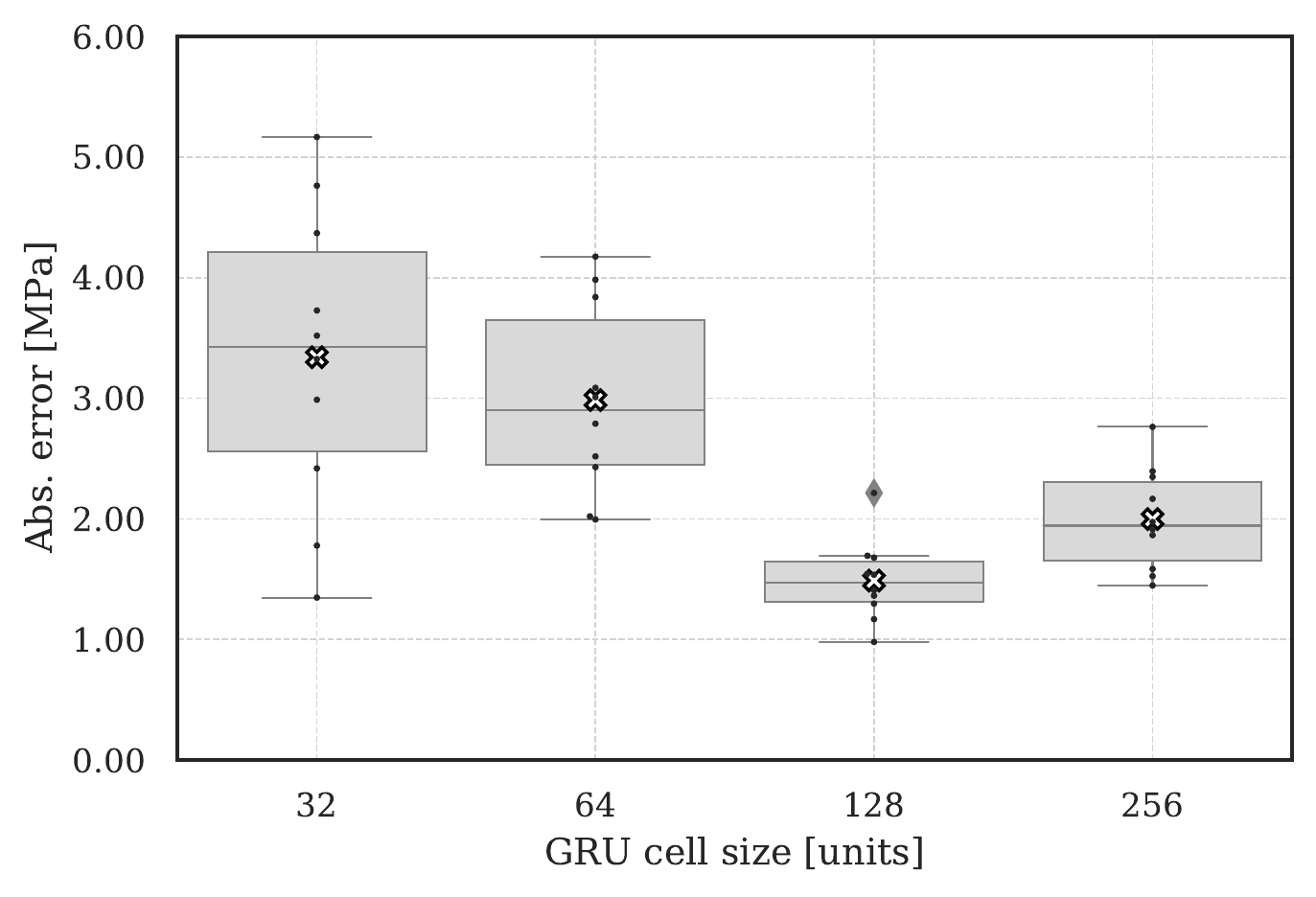}}
\caption{Model selection of PRNN and RNN}
\label{fig:paramstudy}
\end{figure}

For the PRNN, a single material layer without biases and four different sizes are considered. In this case, only the elastoplastic model is used so that the size of the material layer is the only variable in the model selection. Recall that the network still can make a subgroup to behave elastically by passing small strains to the material model. The training and the validation sets consist of 18 Type I curves and 54 Type II curves, respectively. Fig. \ref{fig:paramstudymnn} shows the boxplot with the average validation error of each run alongside the mean error value. The networks with 6 units (or two fictitious material points) performed best.  

For the Bayesian RNN, a single GRU cell is employed and in each initialization, the training set consists of the fixed set of 18 Type I curves, 90 Type II curves randomly chosen from a pool of 1800 curves, and 90 Type III curves also randomly chosen from a pool of 1800 curves, amounting to 198 loading paths. That way, all types of curves used for training in the following sections are covered. Fig. \ref{fig:paramstudyrnn} shows the boxplot with the average training error of each run alongside the mean error value of all runs represented by the x marker. In this study, it is found that the GRU with 128 units performs best. In this case, no validation set is needed to determine the best dropout rate as the type of RNN used in this investigation infers it from the training data by default (see Section \ref{sec:rnn}).

\subsection{Predicting unloading/reloading from monotonic data}
\label{subsec:trlvl2predlvl3}

\begin{figure}[ht!]
\centering
\includegraphics[width=0.54\textwidth]{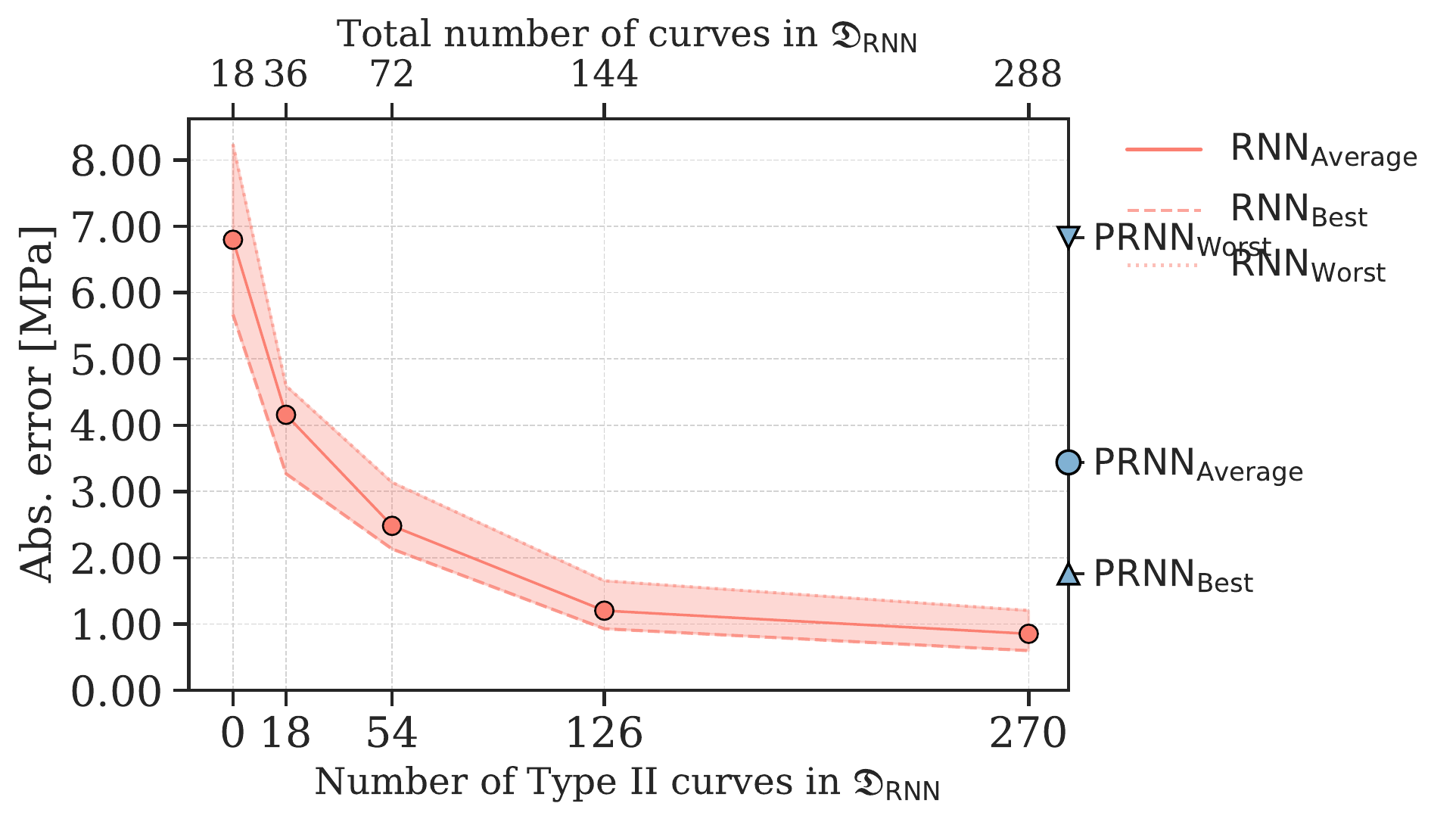}
\caption{Absolute error over monotonic test set $\mathcal{T}_{\mathrm{II}} = \{$100 Type II curves$\}$ for RNNs trained on Types I and II and PRNNs on Type I only}
\label{fig:convtrlvl2predlvl2}
\end{figure}

In this section, the training process of the RNNs on Type II curves (\textit{i.e.} without unloading) is reported. The first test set consists of 100 Type II curves. Fig. \ref{fig:convtrlvl2predlvl2} shows the average error of the RNNs over the test set compared to the best (blue triangle), worst (upside-down blue triangle), and average error (blue circle) found by the PRNNs. Note that the secondary axis starts with 18 curves, this is because the known directions used for training the PRNNs are also a fixed set in the training of the RNNs. The training of the RNNs is stopped with 288 curves. At that stage, a similar level of accuracy between the PRNNs and the RNNs is obtained  (although with a training set 16 times larger) and the addition of new curves only yields a marginal gain in accuracy. 

\begin{figure}[h]
\centering
\subfloat[Absolute error over test set $\mathcal{T}_{\mathrm{III}}= \{$100 Type III curves$\}$ for RNNs trained on Types I and II and PRNNs on Type I only]{\label{fig:convtrlvl2predlvl3}\includegraphics[width=0.54\textwidth, valign = c]{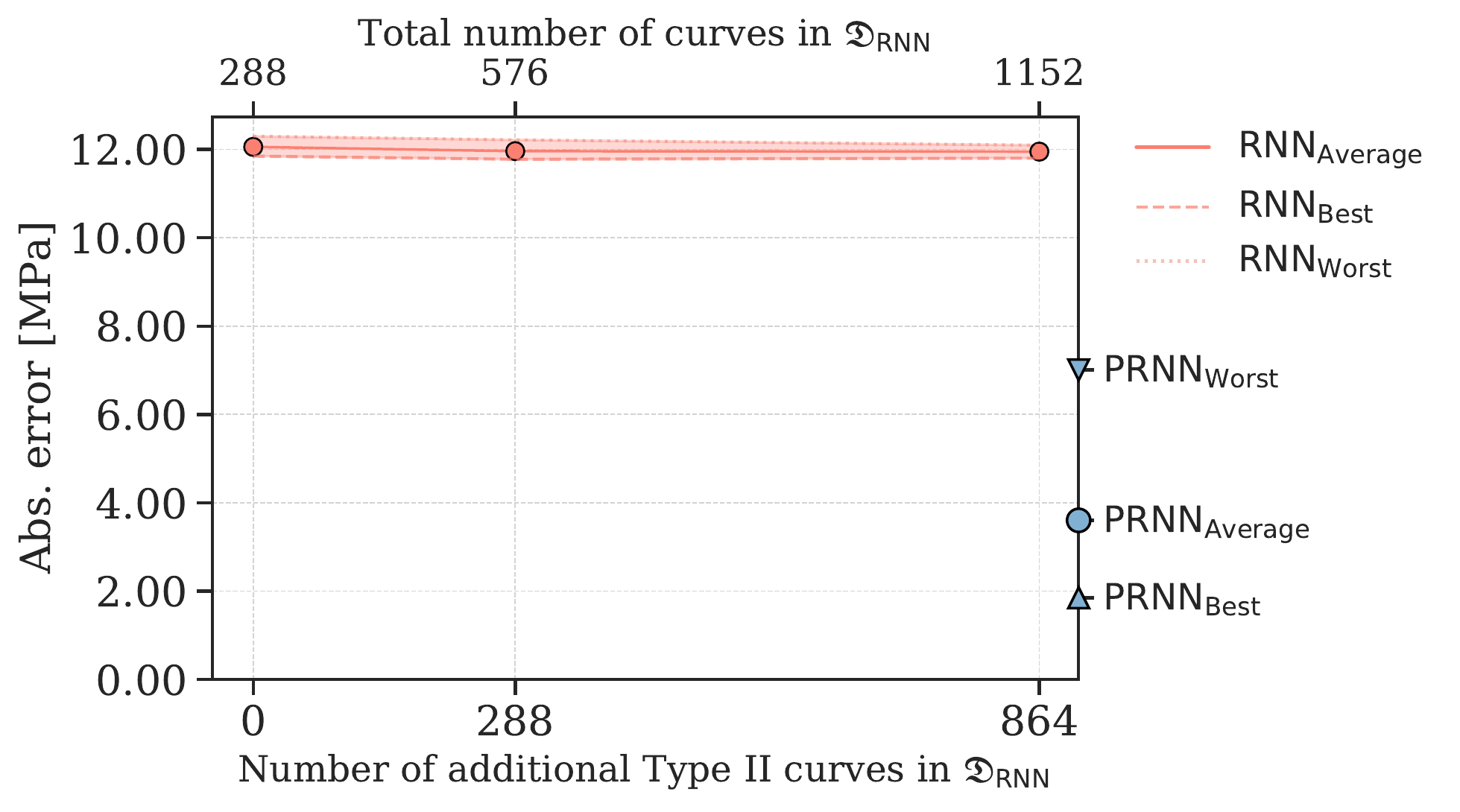}\vphantom{\includegraphics[width=0.45\textwidth,valign=c]{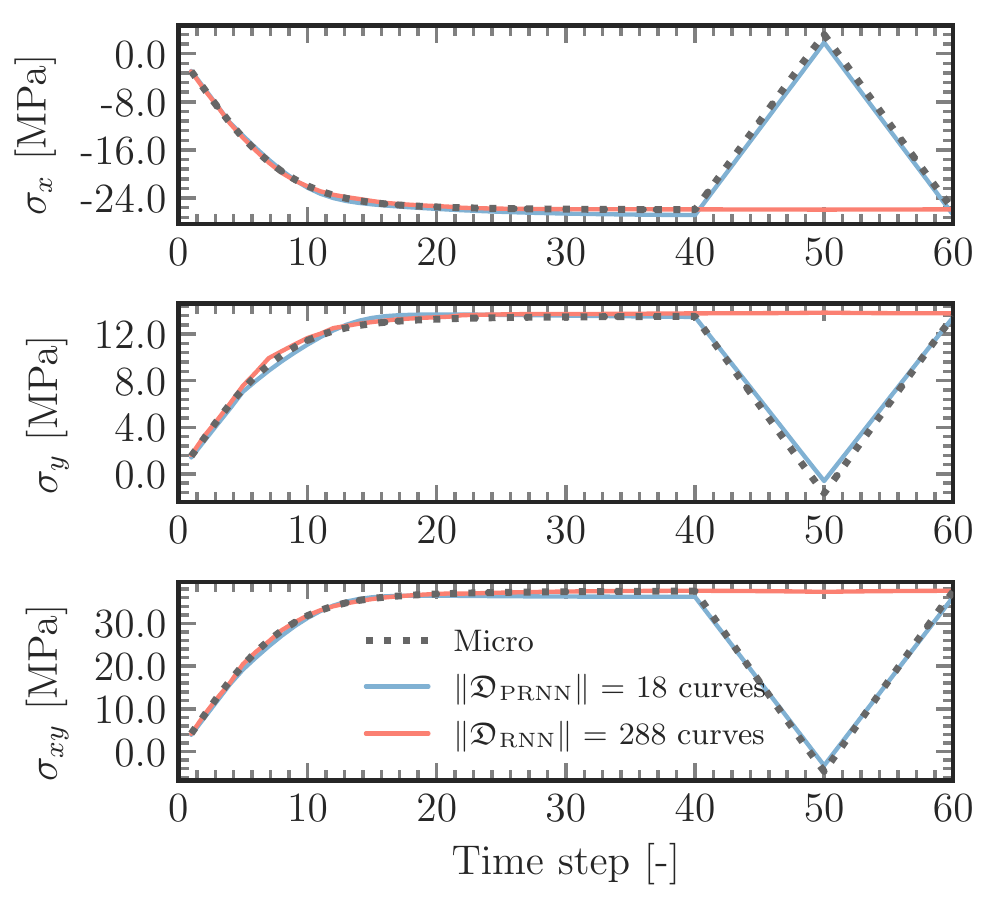}}}
\hfill
\subfloat[Stress-time view of representative case from test set $\mathcal{T}_{\mathrm{III}}$ using the best RNN and PRNN]{\label{fig:samplecurvelevel2onlevel3}\includegraphics[width=0.45\textwidth, valign = c]{pictures/curvestep89.pdf}}
\caption{Absolute error over non-monotonic test set $\mathcal{T}_{\mathrm{III}}$ for RNNs and PRNNs trained only on monotonic data and representative case}
\label{fig:convergencepredlvl2onlvl3}
\end{figure}

Next, a new test set with 100 Type III curves is evaluated by the same networks trained on the 288 monotonic curves. This time, the RNN fails to capture unloading and the addition of more monotonic data is ineffective, as shown in Fig. \ref{fig:convtrlvl2predlvl3}. This outcome is not new to the literature and it is not surprising that RNNs need to see unloading behavior during training in order to be able to describe it. However, in contrast to the RNNs, the PRNNs provide the same level of accuracy for the test sets with and without unloading, even when not exposed to unloading data during training. In Fig. \ref{fig:samplecurvelevel2onlevel3} a single representative case from test set $\mathcal{T}_{\mathrm{III}}$ is plotted using the best RNN and PRNN. Both networks show good agreement with the reference solution until unloading starts (a feature not covered during training), but only the PRNN is capable of capturing the elastic unloading/reloading.

\subsection{Predicting unloading/reloading behavior from non-monotonic data}
\label{subsec:trlvl2and3predlvl3}

\begin{figure}[ht!]
\centering
\includegraphics[width=0.54\textwidth]{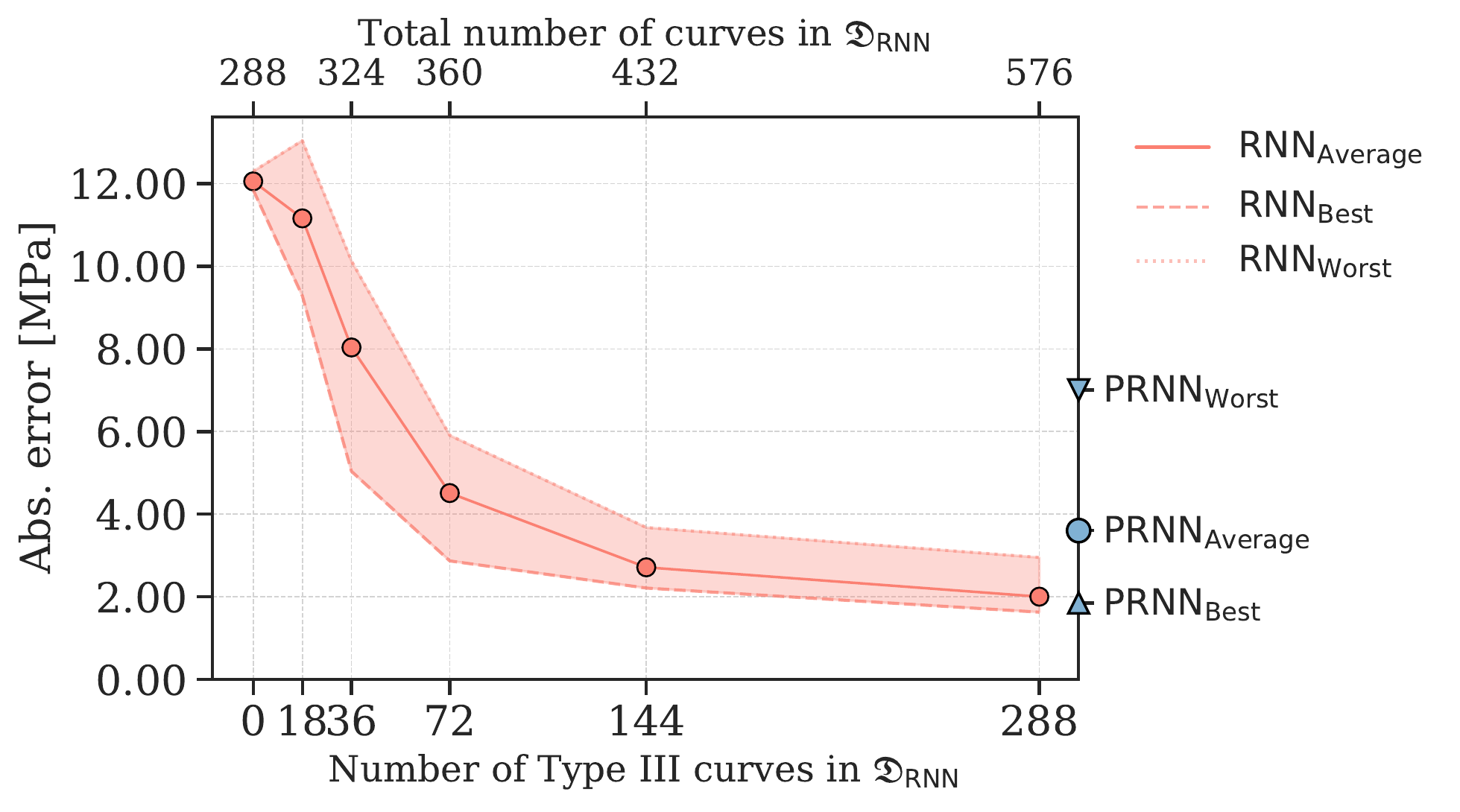}
\caption{Absolute error over non-monotonic test set $\mathcal{T}_{\mathrm{III}} = \{$100 Type III curves$\}$ for RNNs trained on Types I, II and II and PRNNs on Type I only}
\label{fig:convergencepredlvl3onlvl3}
\end{figure}

Following the conclusions of Section \ref{subsec:trlvl2predlvl3}, the training set of the RNN is expanded to include curves with the same unloading behavior as the one observed in the test set $\mathcal{T}_{\mathrm{III}}$. The 288 monotonic curves of Types I and II from the previous section are combined with an increasing number of non-monotonic curves of Type III. This time, with the right features included in the training set, Fig. \ref{fig:convergencepredlvl3onlvl3} shows a monotonic decrease of the average error for the RNN on $\mathcal{T}_{\mathrm{III}}$ curves. However, the performance of the RNN only meets the one obtained by the PRNN with around 32 times more data.   

\subsection{Predicting unseen patterns from non-monotonic data}
\label{subsec:trlvl2and3predunseen}

In this section, three additional test sets are considered for the RNN trained on Types I, II and III and the PRNN trained on Type I only. The goal is to test the ability of the networks to predict the macroscopic stress with patterns different from those seen during training. 

\begin{figure}[ht!]
\centering
\subfloat[Absolute error over test set $\mathcal{T}_{\mathrm{IVa}}= \{$100 Type IVa curves$\}$ for RNNs trained on Types I, II and III and PRNN on Type I only]{\label{fig:convlevel3onfar}\includegraphics[width=0.54\textwidth, valign = c]{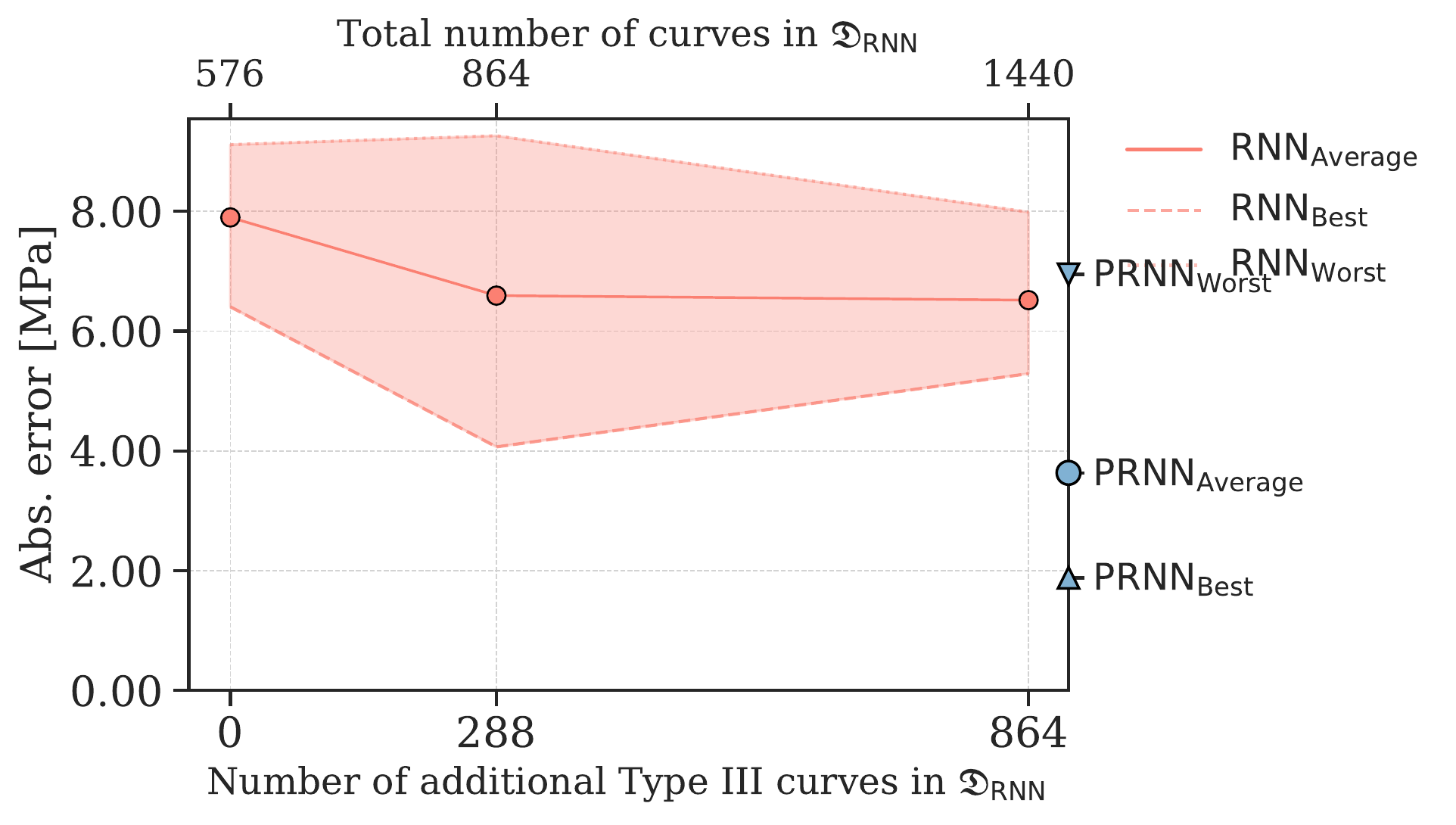}\vphantom{\includegraphics[width=0.45\textwidth,valign=c]{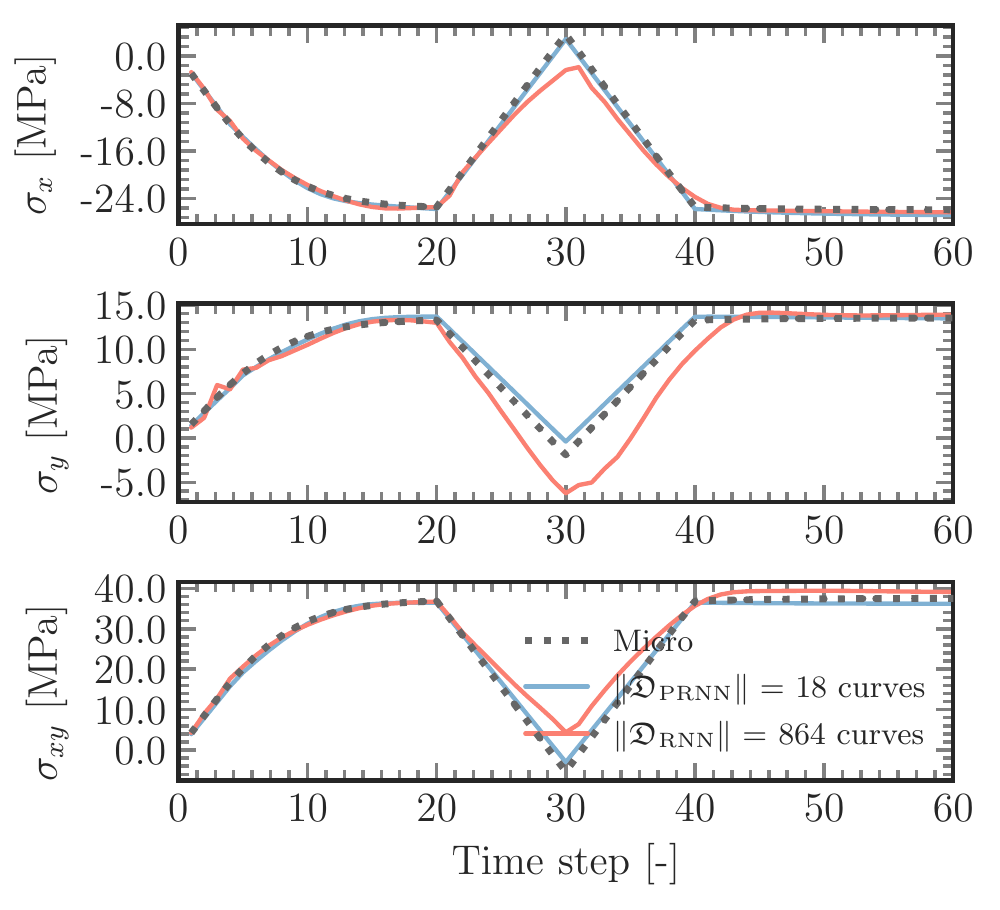}}}
\hfill
\subfloat[Stress-time view of representative case from test set $\mathcal{T}_{\mathrm{IVa}}$ using the best RNN and PRNN]{\label{fig:samplecurvelevel3onlevel3far}\includegraphics[width=0.45\textwidth, valign = c]{pictures/curvestep90_864curves.pdf}}
\caption{Absolute error over non-monotonic test set $\mathcal{T}_{\mathrm{IVa}} = \{$100 Type IVa curves$\}$ for RNNs trained on Types I, II and II and PRNNs on Type I only and representative case}
\label{fig:convergenceunl}
\end{figure}

First, we consider a test set with 100 unseen curves of Type IVa, which consists of proportional curves in random directions with a different predefined unloading/reloading behavior than that of Type III. The average error for that set is shown in Fig. \ref{fig:convlevel3onfar}. Here, the 576 curves from the previous section are no longer enough to provide good accuracy when predicting a different unloading/reloading. By adding more Type III curves to the training set of the RNN, the average error decreases from 6.4 MPa to around 4.0 MPa but no significant gain in the accuracy is observed when the total number of curves is larger than 864 curves. Based on that, a representative case from test set $\mathcal{T}_{\mathrm{IVa}}$ is shown in Fig. \ref{fig:samplecurvelevel3onlevel3far}. Despite the relative low error from both networks, note that the RNN looses performance once unloading starts while the PRNN continues to show good agreement throughout the entire loading path.     

\begin{figure}[ht!]
\centering
\subfloat[Absolute error over test set $\mathcal{T}_{\mathrm{IVb}}= \{$100 Type IVb curves$\}$ for RNNs trained on Types I, II and III and PRNNs on Type I only]{\label{fig:convlevel3onstep}\includegraphics[width=0.54\textwidth, valign = c]{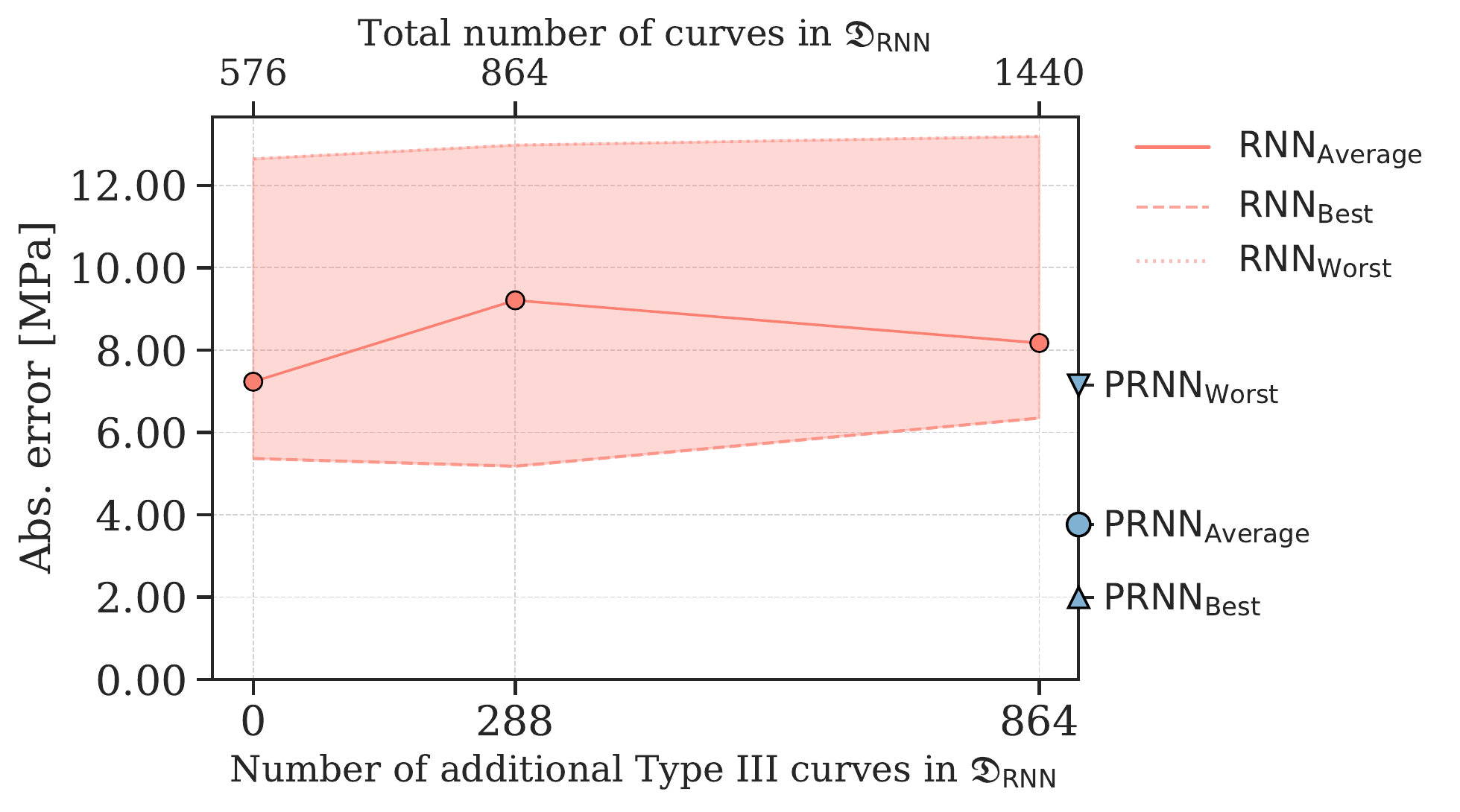}\vphantom{\includegraphics[width=0.45\textwidth,valign=c]{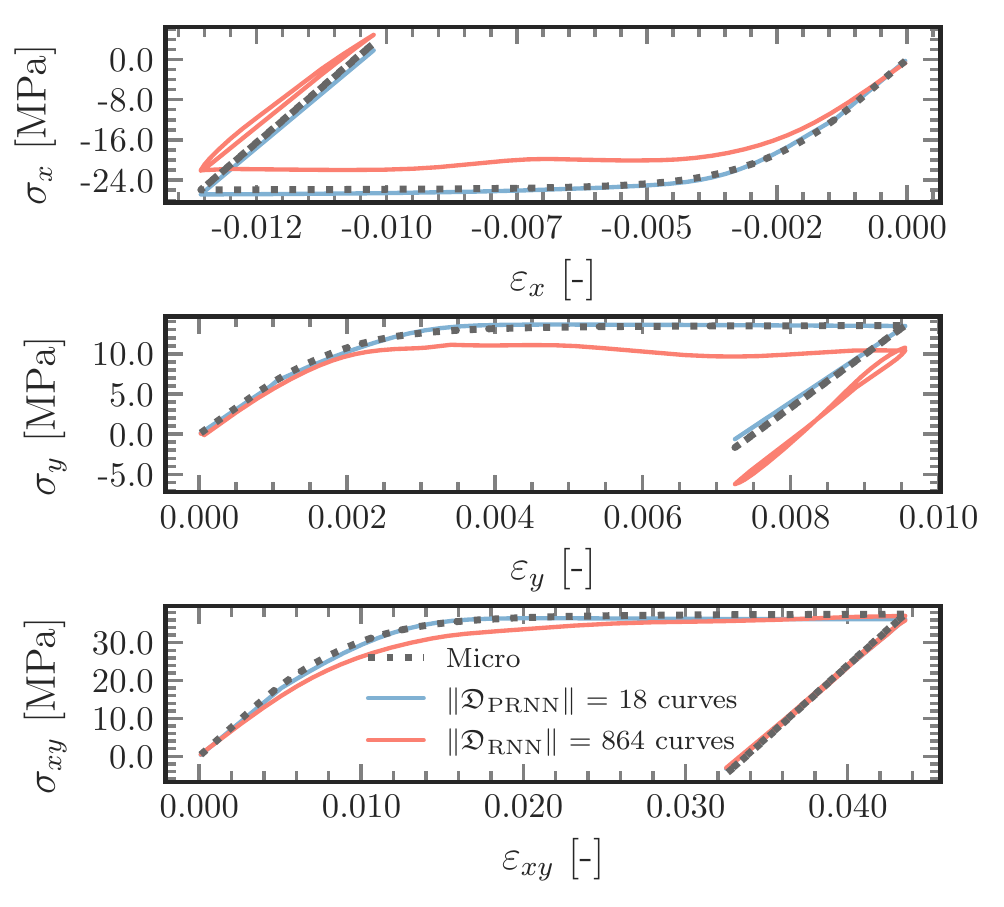}}}
\hfill
\subfloat[Representative $\strainstresscurve$ curve from test set $\mathcal{T}_{\mathrm{IVb}}$ using the best RNN and PRNN]{\label{fig:samplecurvelevel3onstep}\includegraphics[width=0.45\textwidth, valign = c]{pictures/curve91.pdf}}
\caption{Absolute error over different step size test set $\mathcal{T}_{\mathrm{IVb}} = \{$100 Type IVb curves$\}$ for RNNs trained on Types I, II and III and PRNNs on Type I only and representative case}
\label{fig:convergencepredlvl3onstep}
\end{figure}

For the next scenario, a test set with 100 Type IVb curves are considered. These curves have the same unloading/reloading behavior as Type III, but with a 10$\times$ smaller time step. Fig. \ref{fig:convlevel3onstep} illustrates the average error of 10 networks over that test set and again, it is clear that the addition of new curves with patterns different from the exact one being tested is not beneficial to the RNN. Again, the PRNN provides good accuracy. Essentially, the PRNN is only as sensitive to step size as the material models embedded in it. Fig. \ref{fig:samplecurvelevel3onstep} illustrates the networks' predictions for a curve in the test set $\mathcal{T}_{\mathrm{IVb}}$.

\begin{figure}[ht!]
\centering
\subfloat[Absolute error over test set $\mathcal{T}_{\mathrm{V}}= \{$100 Type V curves$\}$ for RNN trained on Types I, II and III and PRNN on Type I only]{\label{fig:convlevel3ongp}\includegraphics[width=0.54\textwidth, valign = c]{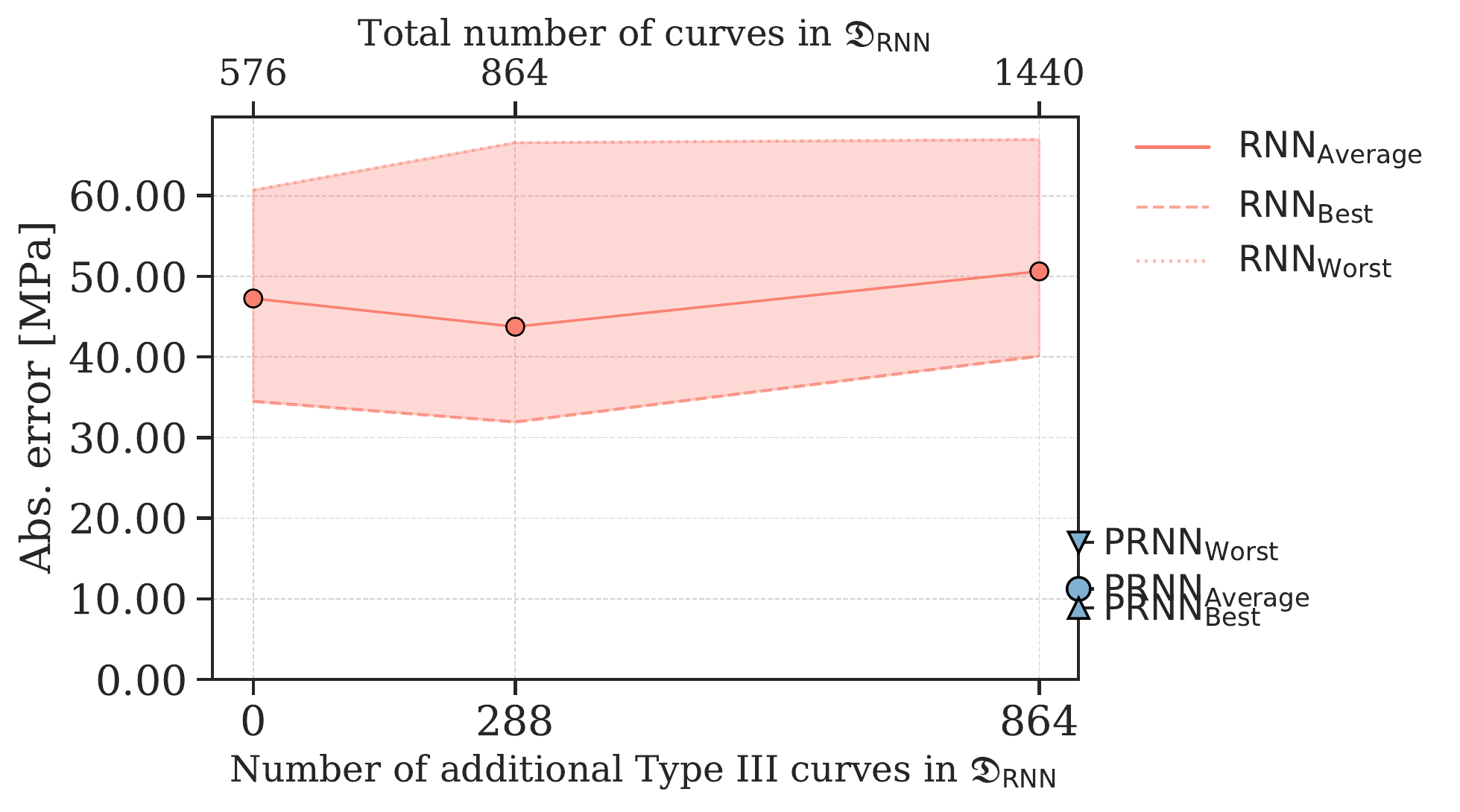}\vphantom{\includegraphics[width=0.45\textwidth,valign=c]{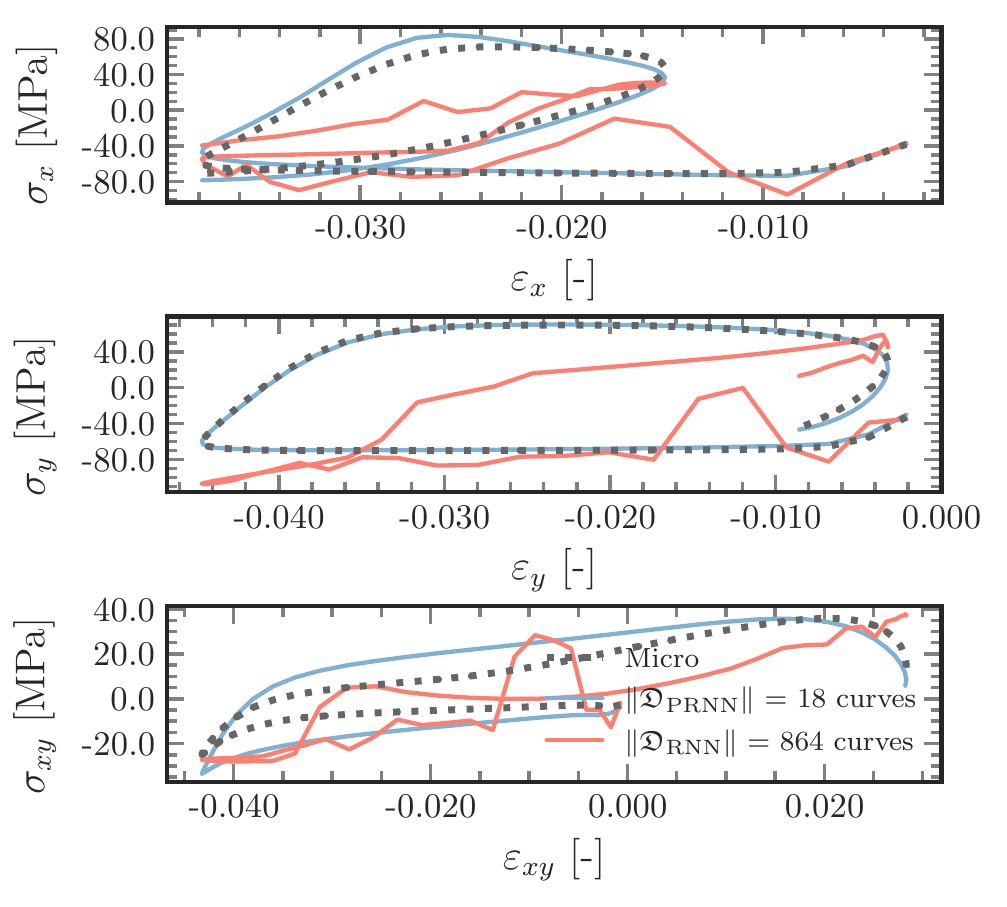}}}
\hfill
\subfloat[Representative $\strainstresscurve$ curve from test set $\mathcal{T}_{\mathrm{V}}$ using the best RNN and PRNN]{\label{fig:samplecurvelevel3ongp}\includegraphics[width=0.45\textwidth, valign = c]{pictures/curve26.pdf}}
\caption{Absolute error over non-proportional and non-monotonic test set $\mathcal{T}_{\mathrm{V}}= \{$100 Type V curves$\}$ for RNNs trained on Types I, II and III and PRNNs on Type I only and representative case}
\label{fig:convergenceunseenpatterns}
\end{figure}

As a final test, a set of 100 Type V curves, which corresponds to non-proportional and non-monotonic paths, is considered. This type of curve combines the two previous features: different step sizes and different unloading/reloading locations. Fig. \ref{fig:convlevel3ongp} shows the average error for additional non-monotonic curves in the training of the RNNs. It is clear that the RNN completely fails to capture non-proportional paths (lowest error around 32 MPa) and that the addition of more data with features different than those being tested is a waste of resources. Although in different levels, a similar trend of loss of accuracy is also observed in the PRNNs, where the best, average and worse performances result in errors around 8.9 MPa, 17.0 MPa and 11.2 MPa respectively. Fig. \ref{fig:samplecurvelevel3ongp} illustrates the different order of error between the PRNN and the RNN on a representative case from test set $\mathcal{T}_{\mathrm{V}}$. 

\subsection{Training both networks on non-monotonic and non-proportional loading}
\label{subsec:trlvlgppredlvlall}

In this section, both networks are trained on the most generic set of curves, \textit{i.e.} random non-monotonic and non-proportional curves of Type V. In addition to that, we trained the PRNNs on the known and proportional loading cases for comparison purposes. In that case, the size of the material layer is kept at 6 units and three training dataset sizes are considered. First, only the pure uniaxial cases are included, which yields 6 loading cases. Then, the 4 biaxial cases are added to the previous training dataset, resulting in 10 loading cases. And finally, we add the 8 cases with biaxial and shear loading, which amounts to the 18 fundamental paths shown in Fig. \ref{fig:doepropbasic}. 

\begin{figure}[ht!]
\centering
\subfloat[Validation error for PRNNs trained on $\mathfrak{D}_{\mathrm{PRNN}} = \{$198 Type V curves$\}$ and $\mathcal{V}_{\mathrm{PRNN}} = \{$54 Type V curves$\}$]{\label{fig:modselmnnrandom}\includegraphics[width=0.48\textwidth]{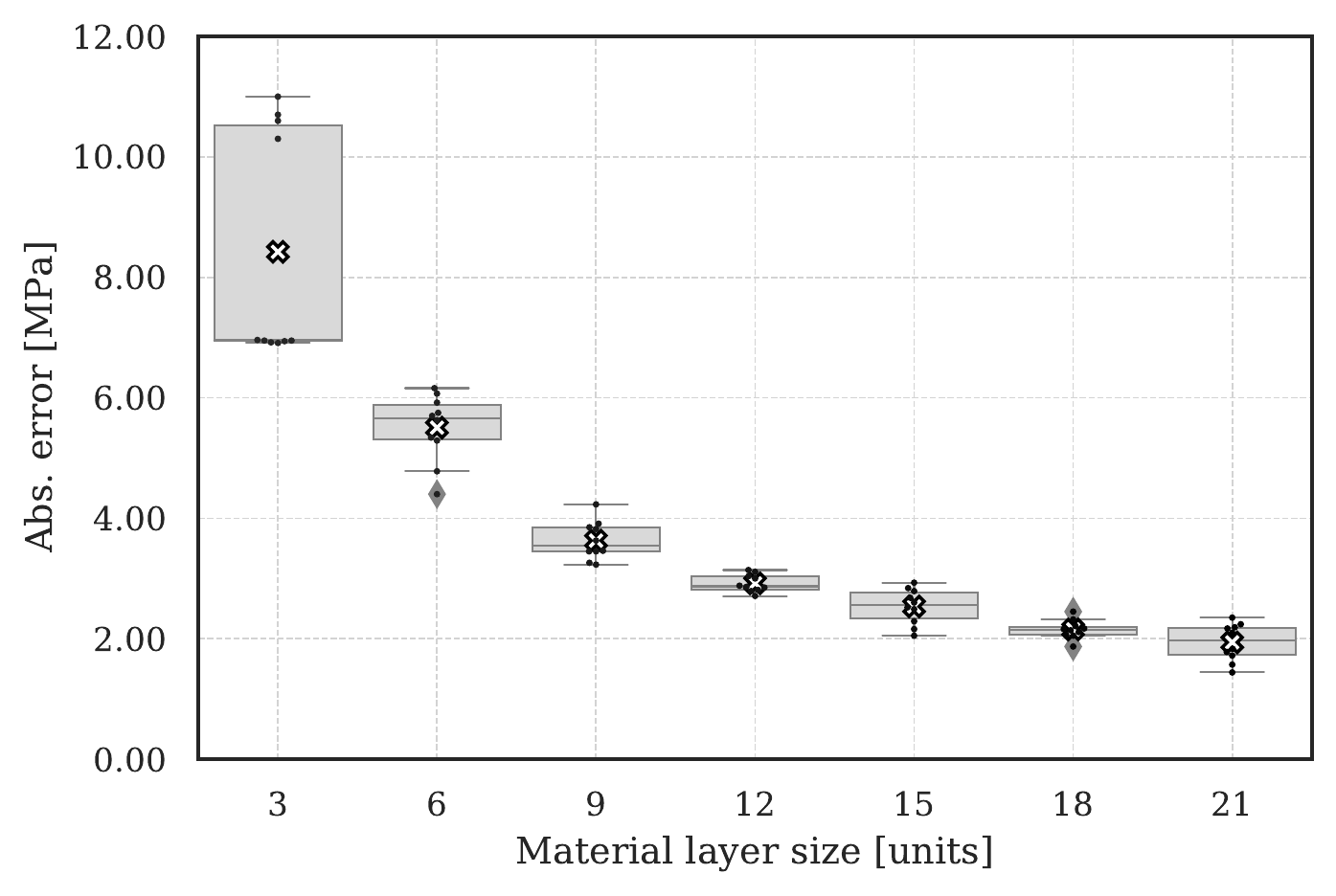}}
\hfill
\subfloat[Training error for RNNs trained on $\mathfrak{D}_{\mathrm{RNN}} = \{$2304 Type V curves$\}$]{\label{fig:modselrnnrandomnonprop}\includegraphics[width=0.48\textwidth]{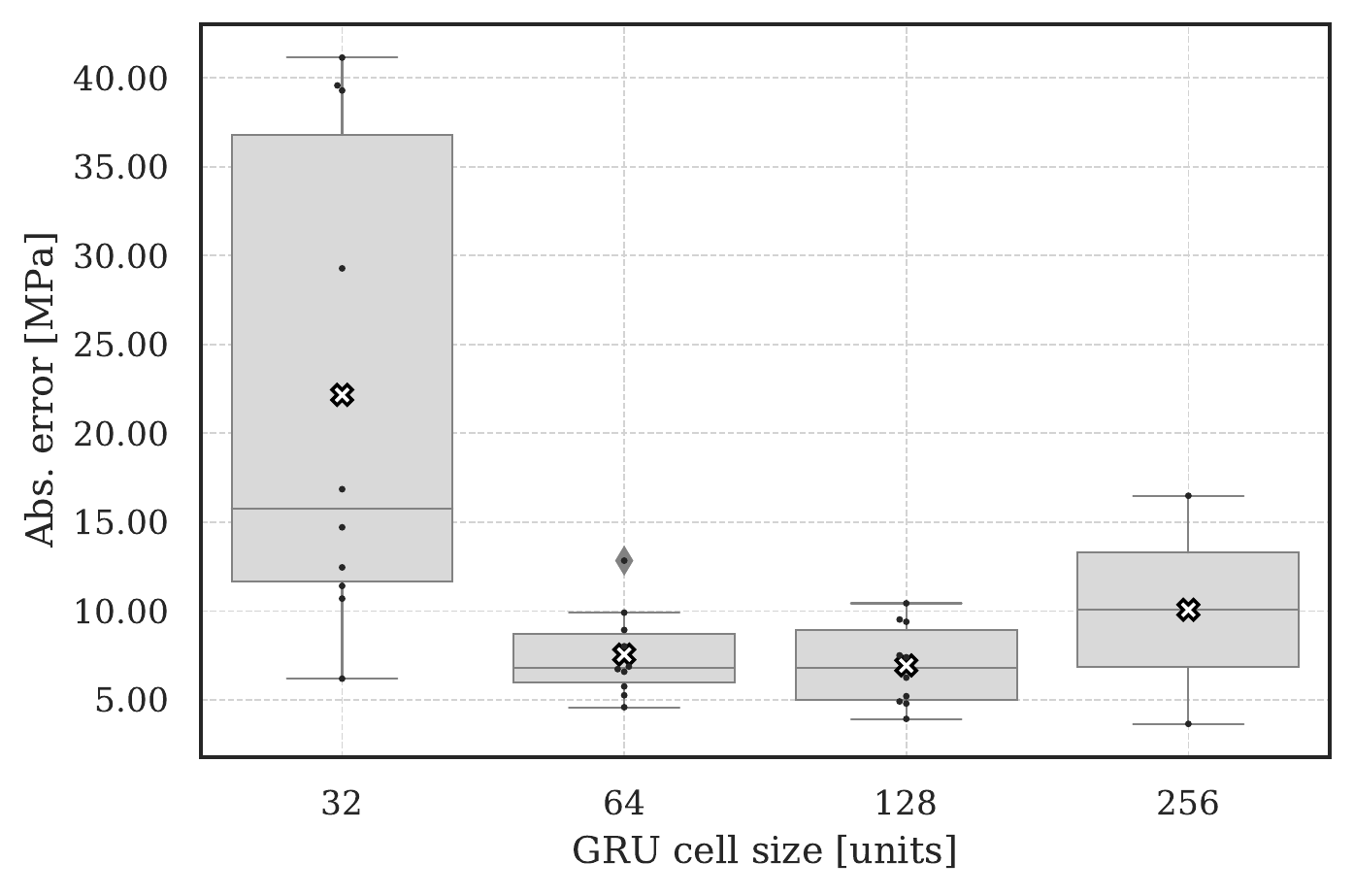}}
\caption{Model selection of PRNN and RNN for non-monotonic and non-proportional loading}
\label{fig:paramstudynonprop}
\end{figure}

When training both networks on non-proportional paths a new model selection procedure was carried out to determine the optimum size of the material layer and the GRU cell, respectively. In this preliminary study, 10 different weights initialization are considered again. For training the RNNs and the PRNNs, 2304 and 198 Type V curves are used, respectively. Figs. \ref{fig:modselmnnrandom} and \ref{fig:modselrnnrandomnonprop} show the boxplot with the average error of each run alongside the mean error value. In this case, the networks with 18 units (which corresponds to six fictitious material points) performed better. For the Bayesian RNN, the GRUs with 128 units continue to provide the best performance. Therefore, this architecture is the one used in the comparison presented below.  

\begin{figure}[h!]
\centering
\subfloat[Test set $\mathcal{T}_{\mathrm{II}}$ = \{100 Type II curves\}]{\label{fig:converrordiftrtestlvl2}\includegraphics[width=0.46\textwidth]{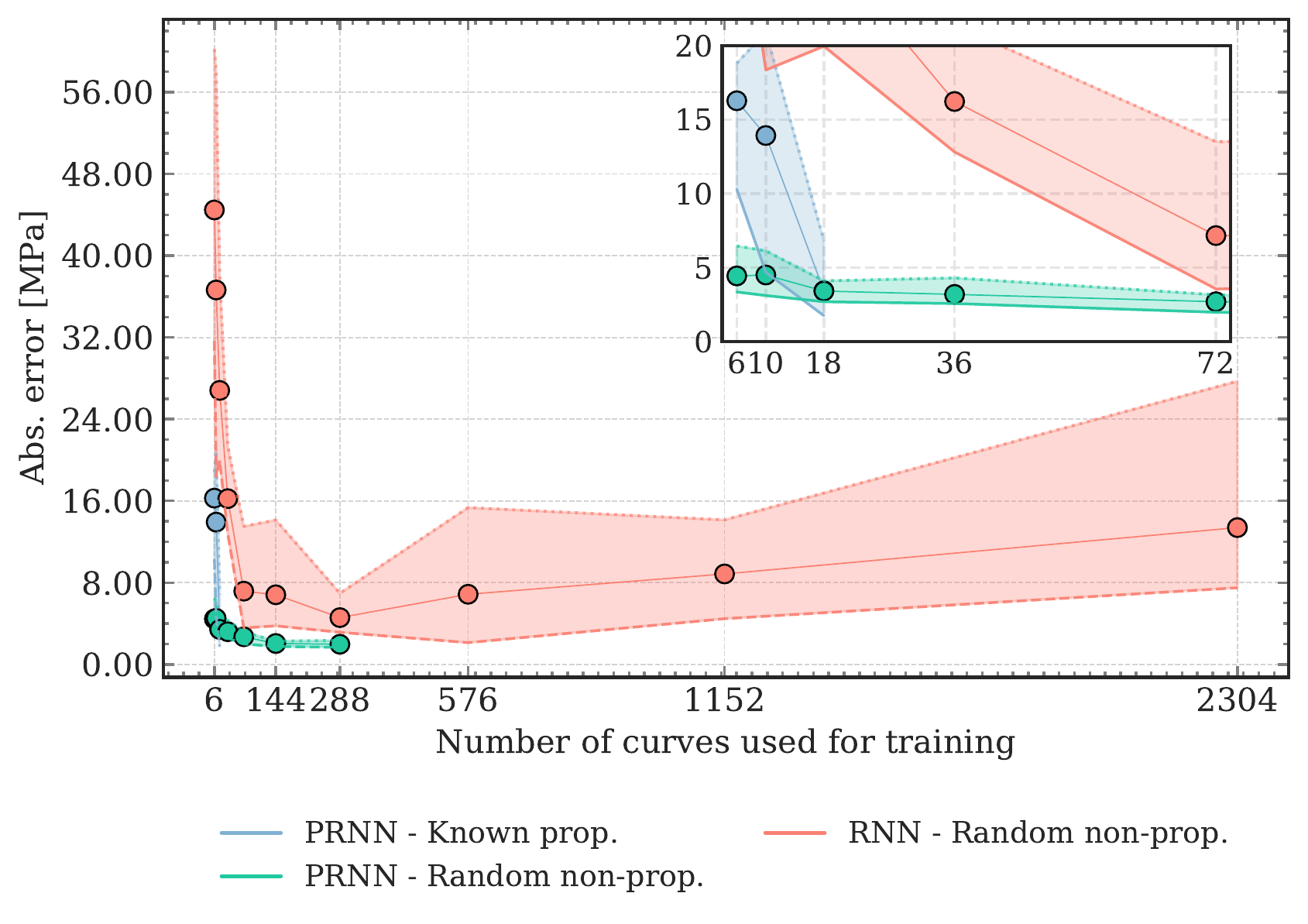}}
\hfill
\subfloat[Test set $\mathcal{T}_{\mathrm{III}}$ = \{100 Type III curves\}]{\label{fig:converrordiftrtestlvl3}\includegraphics[width=0.46\textwidth]{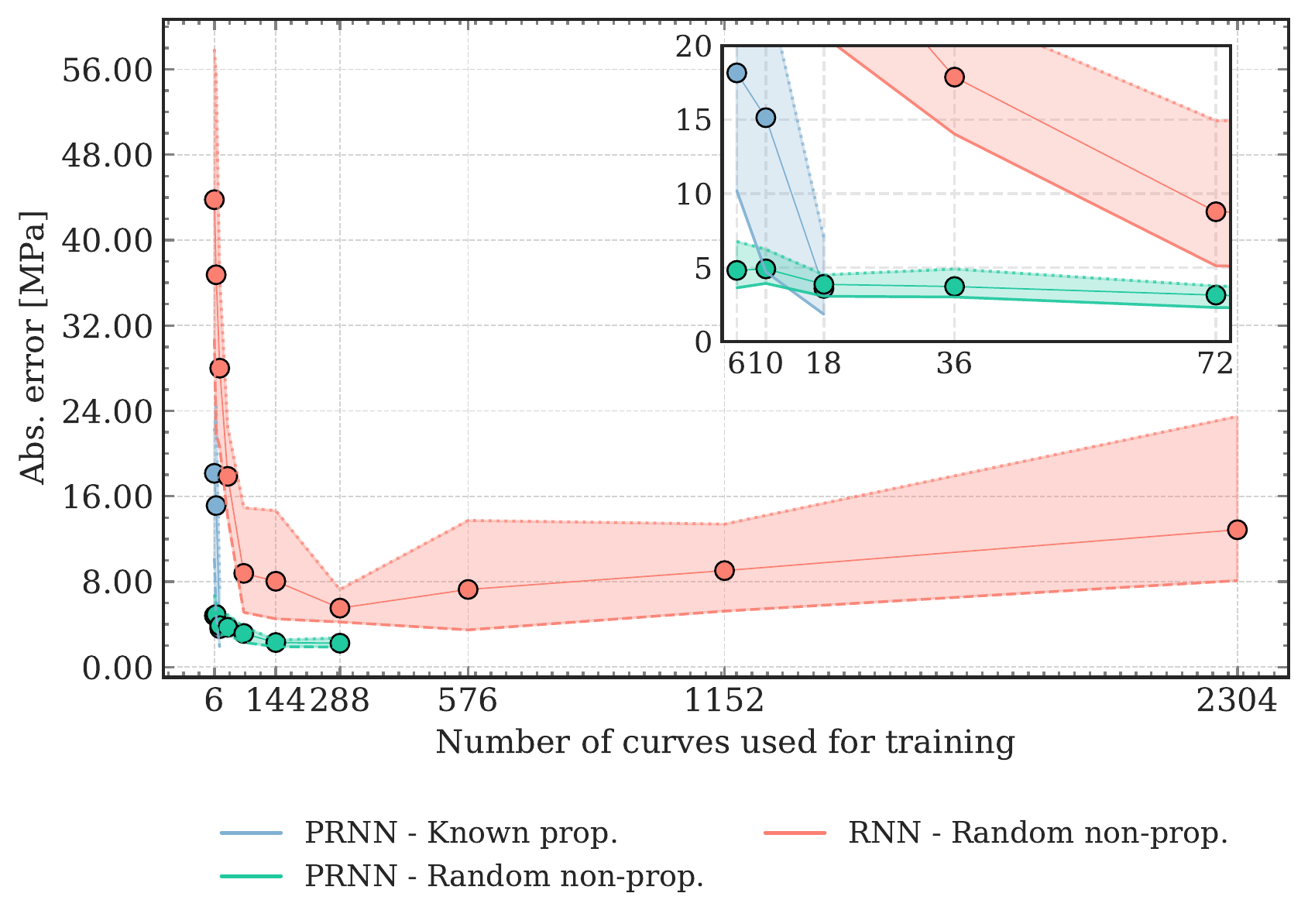}}
\hfill
\subfloat[Test set $\mathcal{T}_{\mathrm{IV}_\mathrm{a}}$ = \{100 Type IVa curves\}]{\label{fig:converrordiftrtestfar}\includegraphics[width=0.46\textwidth]{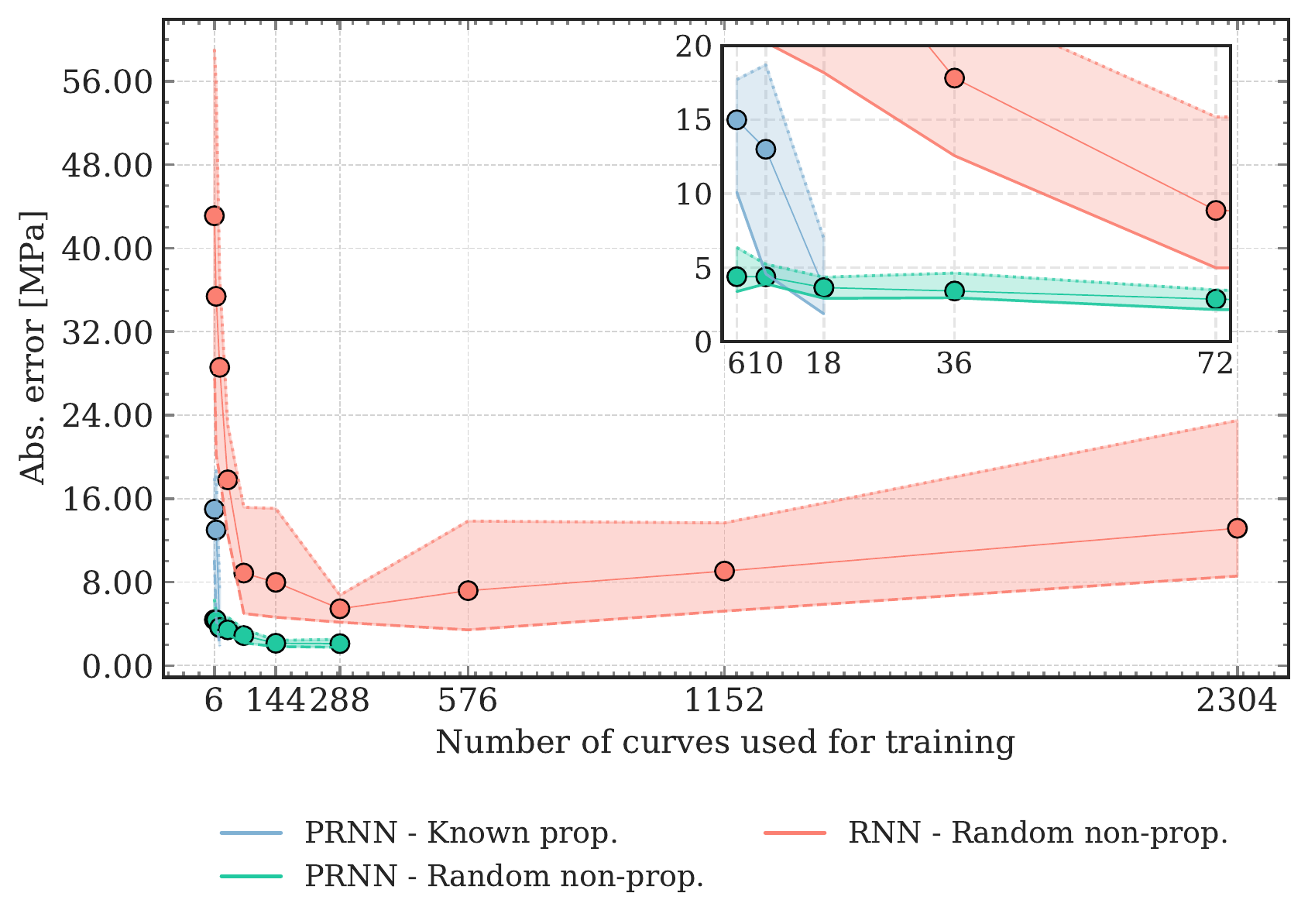}}
\hfill
\subfloat[Test set $\mathcal{T}_{\mathrm{IV}_\mathrm{b}}$ = \{100 Type IVb curves\}]{\label{fig:converrordiftrteststep}\includegraphics[width=0.46\textwidth]{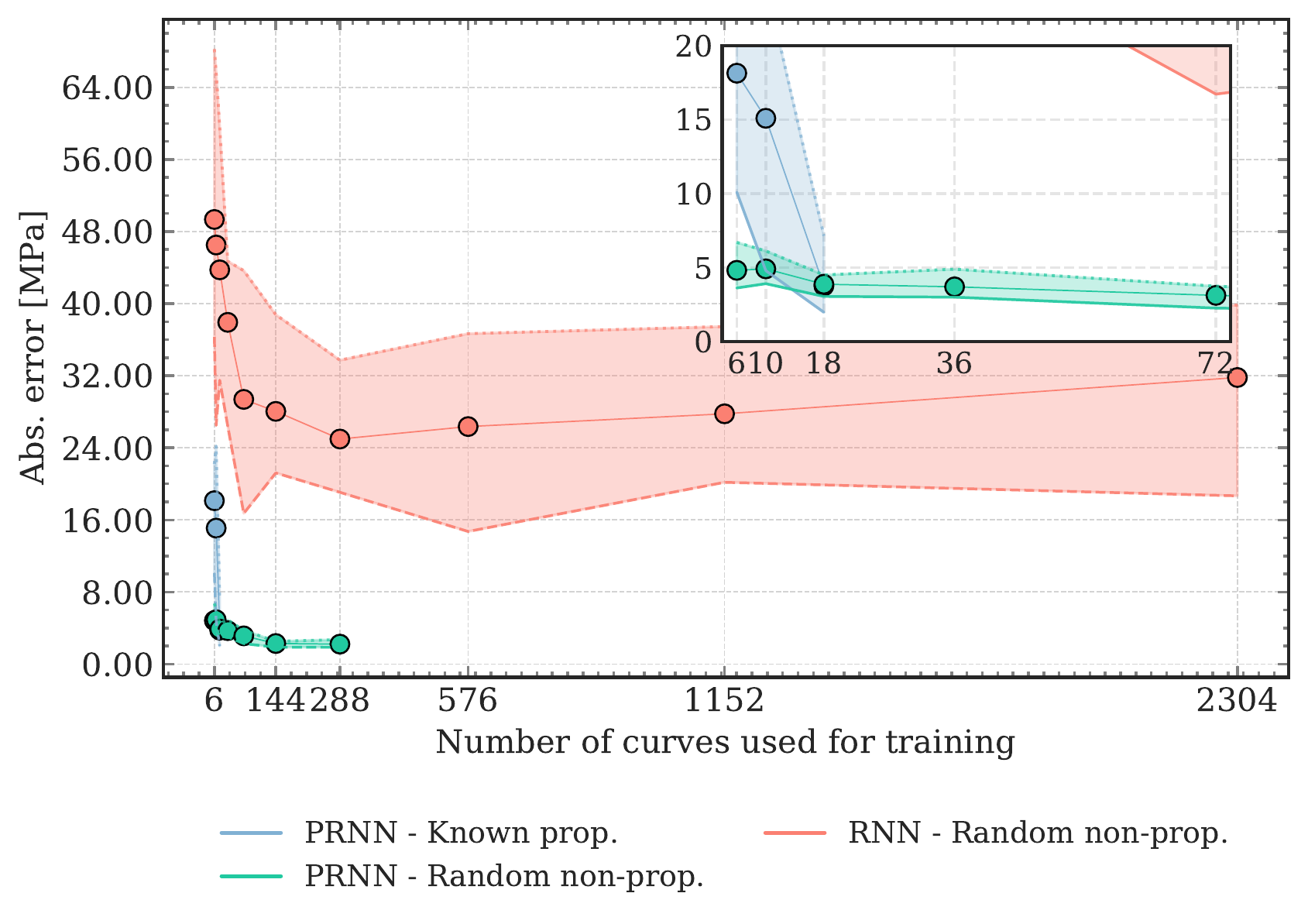}}
\hfill
\subfloat[Test set $\mathcal{T}_{\mathrm{V}}$ = \{100 Type V curves\}]{\label{fig:converrordiftrtestgp}\includegraphics[width=0.46\textwidth]{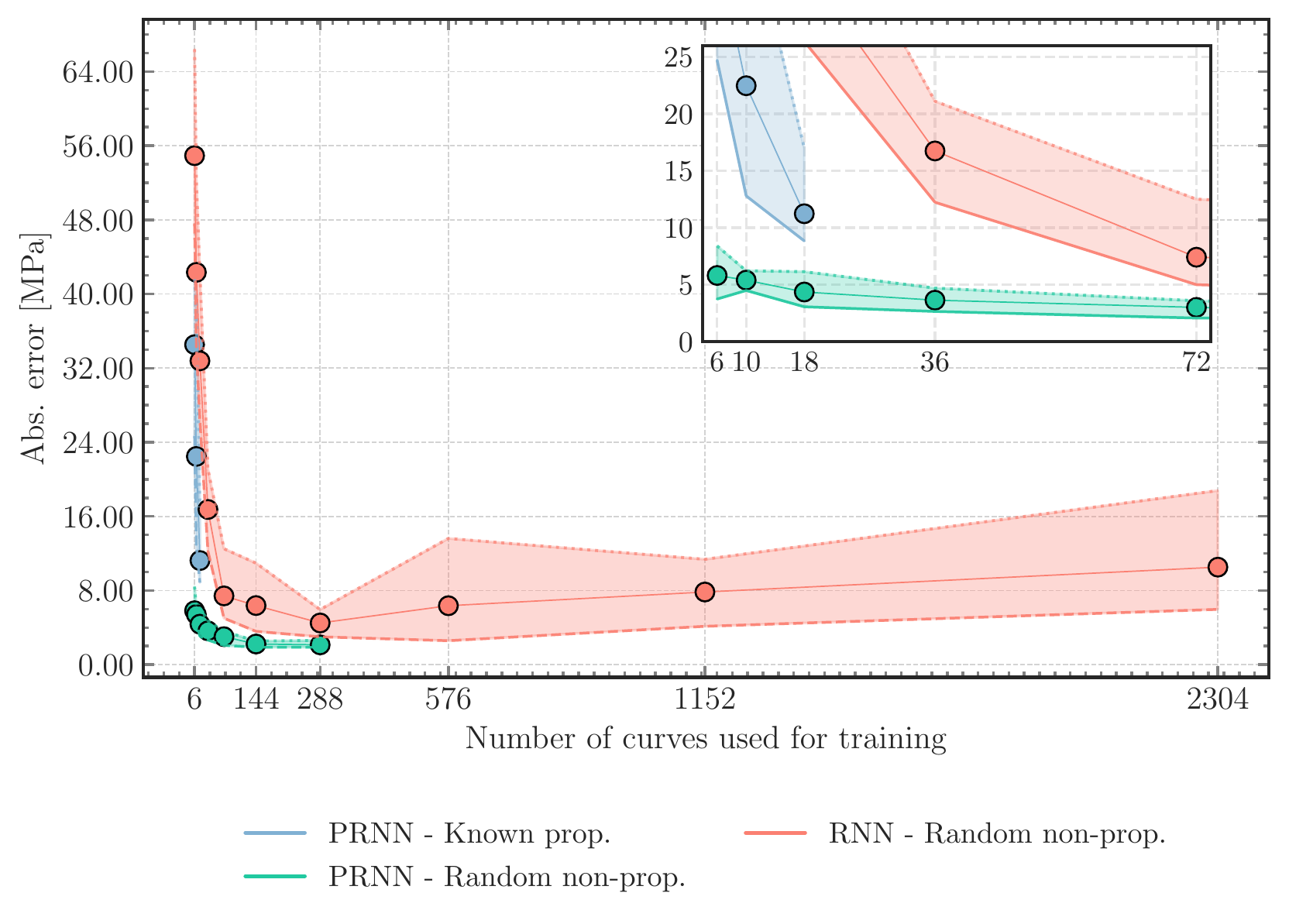}}
\caption{Absolute error of networks trained on different sampling strategies and different test sets}
\label{fig:convergencediftr}
\end{figure}

This time, all test sets discussed in Section \ref{subsec:trlvl2predlvl3}, \ref{subsec:trlvl2and3predlvl3} and \ref{subsec:trlvl2and3predunseen} are used again to assess the accuracy of the networks with the new sampling strategy. Figures \ref{fig:converrordiftrtestlvl2}-\ref{fig:converrordiftrtestgp} show the best, worst and average error for the cases studied so far in order. Based on this study, a few important insights are worth mentioning: after a certain point (around 576 curves), the RNNs reach an optimum level of accuracy and the addition of new curves no longer boosts predictions for proportional loading cases ($\mathcal{T}_{\mathrm{II}}$, $\mathcal{T}_{\mathrm{III}}$, $\mathcal{T}_{\mathrm{IVa}}$ and $\mathcal{T}_{\mathrm{IVb}}$). This is in line with the behavior observed in the previous sections, in which the RNNs would only perform well when trained with the same features as in the test set. And more importantly, changing the sampling strategy also showed to have limited effect on improving their performance. Granted, increasing even further the number of curves used for training as well as the complexity of the RNN might help in that task. However, the point stands that with the PRNN, this is not necessary. Note that for the same training set sizes, the PRNN with either the known or the random curves perfoms better than using the RNNs. 

Finally, when choosing between known and random curves for training the PRNN, the latter shows comparable errors with the first when predicting proportional loading, but is significantly more accurate (see detail in Fig. \ref{fig:converrordiftrtestgp}) for non-proportional loading. For that reason, the PRNN trained on Type V is chosen to illustrate the network's capacity in the following \fesqr \ examples. On average, the accuracy of the PRNN reaches a plateau around 36 curves. From that point on, the benefit of adding new data is limited. 

\section{\fesqr \ applications}
\label{sec:numex}

In this section, the PRNN trained with 36 non-proportional and non-monotonic Type V curves and lowest test error for test set $\mathcal{T}_{\mathrm{V}}$ in Section \ref{subsec:trlvlgppredlvlall} is employed as the constitutive model in two numerical examples. Results obtained with the PRNN as constitutive model are compared against results obtained with full \fesqr \ with the same micromodel that the network was trained to be a surrogate for. Both types of analysis are performed with an in-house Finite Element code using the open-source Jem/Jive C++ numerical analysis library \cite{NGUYENTHANH2020}. 

At the macroscale, an arc-length method with an adaptive-stepping scheme \cite{vanderMeer2012} is adopted to tackle potential convergence issues. This way, if a loading step does not converge with the given (full) step size, a reduction factor is applied to it until the loading step converges or until a maximum number of reductions in the initial step has been reached, terminating the analysis. Upon convergence, in the following step, the analysis is resumed with the full step size. In this work, each load step can be reduced by a factor of 0.4 for a maximum number of 5 times. 

All simulations, including the PRNN training, were executed on a single core of a Xeon E5-2630V4 processor on a cluster node with 128 GB RAM running CentOS 7. 

\subsection{Tapered bar}
\begin{figure}[!h]
\centering
\subfloat[Geometry and boundary conditions]{\label{fig:dogbonegeom}
\includegraphics[width=0.38\textwidth, valign = c]{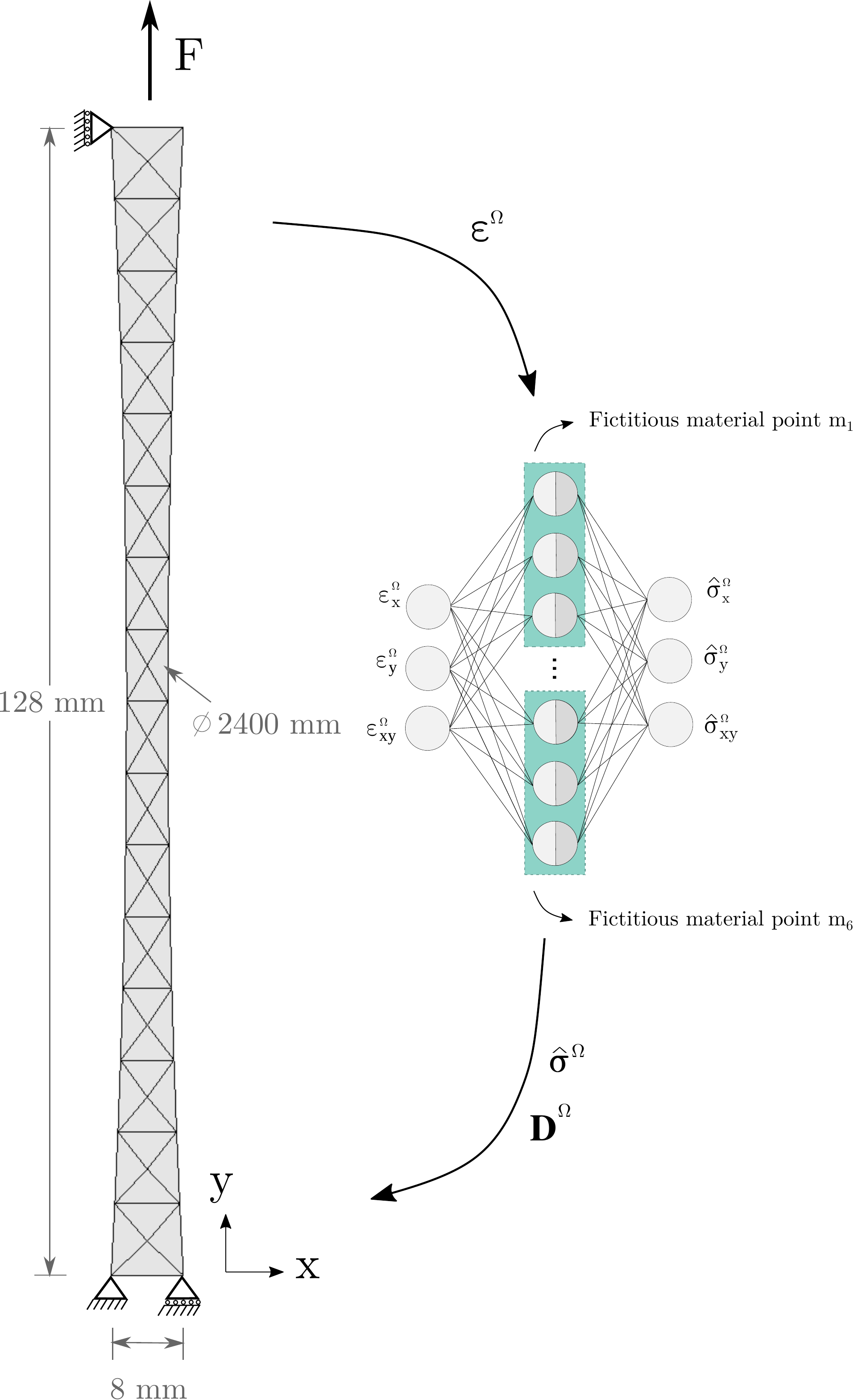}}
\hfill
\subfloat[Loading function]{\label{fig:dogboneshapefunc}\includegraphics[width=0.5\textwidth, valign = c]{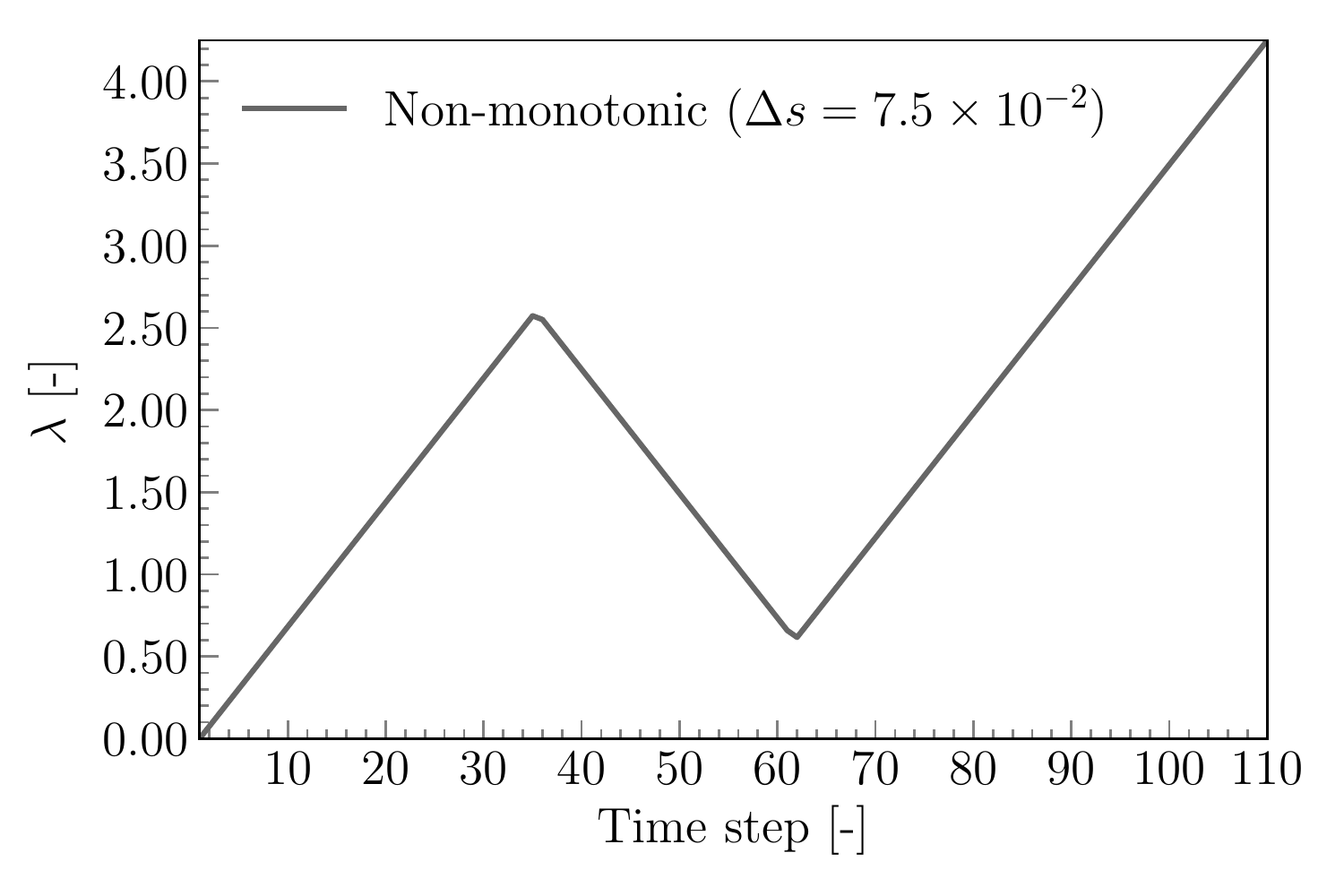}\vphantom{\includegraphics[width=0.38\textwidth,valign=c]{pictures/dogbonegeometry.pdf}}}
\caption{Tapered bar \fesqr \ example with (a) geometry and boundary conditions and (b) loading function}
\label{fig:dogbone} 
\end{figure}

The first example consists in a tapered composite specimen with length of 128 mm and height of 8 mm loaded in transverse tension. In this setting, the 36-fiber RVE model used to train the networks in the previous sections is embedded at each integration point of the macroscale. The geometry and the boundary conditions are shown in Fig. \ref{fig:dogbonegeom}. The \fesqr \ problem is solved for 110 load steps with unloading according to the function shown in Fig. \ref{fig:dogboneshapefunc}. At this point, the macroscopic response is already in the plastic regime. 

\begin{figure}[!ht]
\centering
\subfloat[Strain field at the end of the analysis using the PRNN]{\makebox[3.4\width]{\label{fig:dogbonedispfield}
\includegraphics[width=0.12\textwidth]{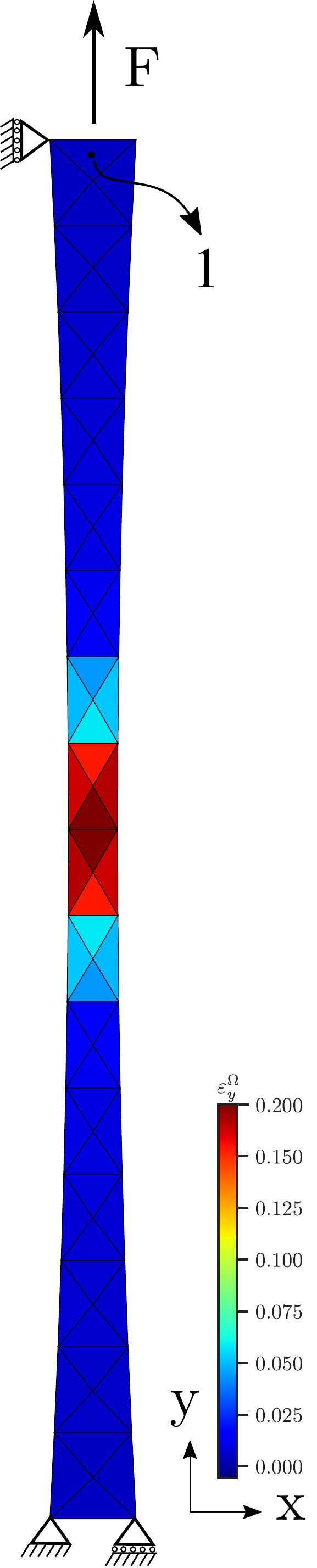}}}
\hfill
\subfloat[Detailed view of PRNN for integration point 1]{\makebox[1.1\width]{\label{fig:dogboneopennet}\includegraphics[width=0.51\textwidth]{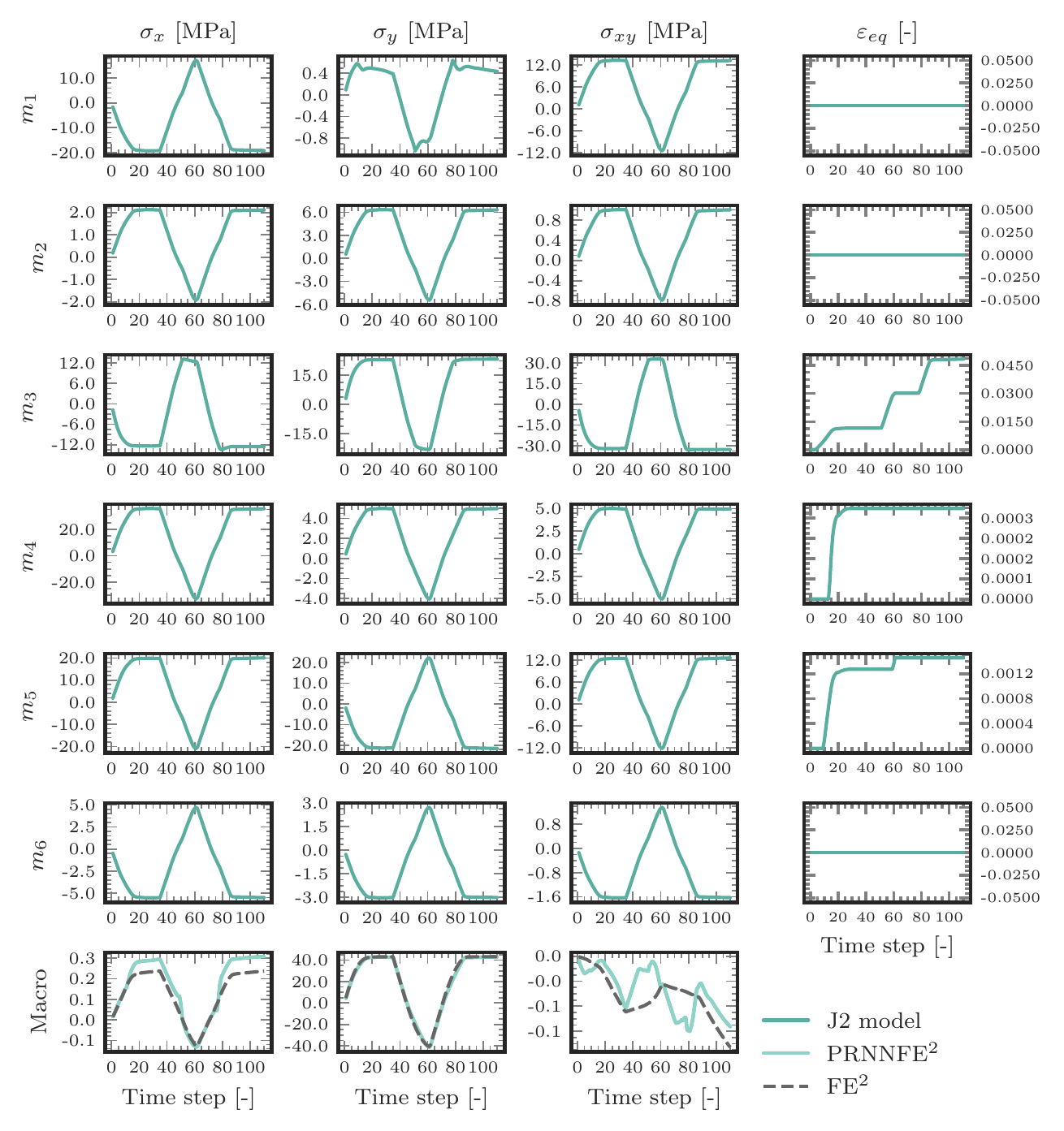}}}
\caption{Strain field using the PRNN on the left and detailed view of PRNN for a single macroscopic integration point on the right.}
\label{fig:dogboneresults} 
\end{figure}

The strain field at the end of the analysis is shown in Fig. \ref{fig:dogbonedispfield} along with the location of one macroscopic integration point. This point is used to illustrate the state of the PRNN throughout the time steps. Recall that the network used in this section consists of 18 units (\textit{i.e.} 6 fictitious material points). Thus, each row in Fig. \ref{fig:dogboneopennet} corresponds to a fictitious material point, each with its own stress path and internal variables (even though only one of them is plotted). 

In this case, each component of the macroscopic response (\textit{i.e.} homogenized stresses) is simply the linear combination of the local stresses of the 6 material models. It can be observed that the two macroscopic responses are in excellent agreement, with minor deviations in stress components with low magnitude. This is only visualized for a single point, but it is emphasized that in this multiscale problem, agreement in the evolution of the stresses in a single integration point indicates that the whole problem is solved accurately. Moreover, the equivalent plastic deformation of each material point is plotted in the last column. Note that despite the plastic response of the RVE after time step 20, three of the fictitious material points of the network ($m_1$, $m_2$ and $m_6$) remain in the elastic regime. 
\begin{figure}[!ht]
\centering
\label{fig:dogboneinfos} 
\subfloat[Load-displacement curve obtained with the full-order \fesqr \ approach and with PRNN for different macroscopic mesh discretizations]{\label{fig:dogbonecurve}\includegraphics[width=0.46\textwidth]{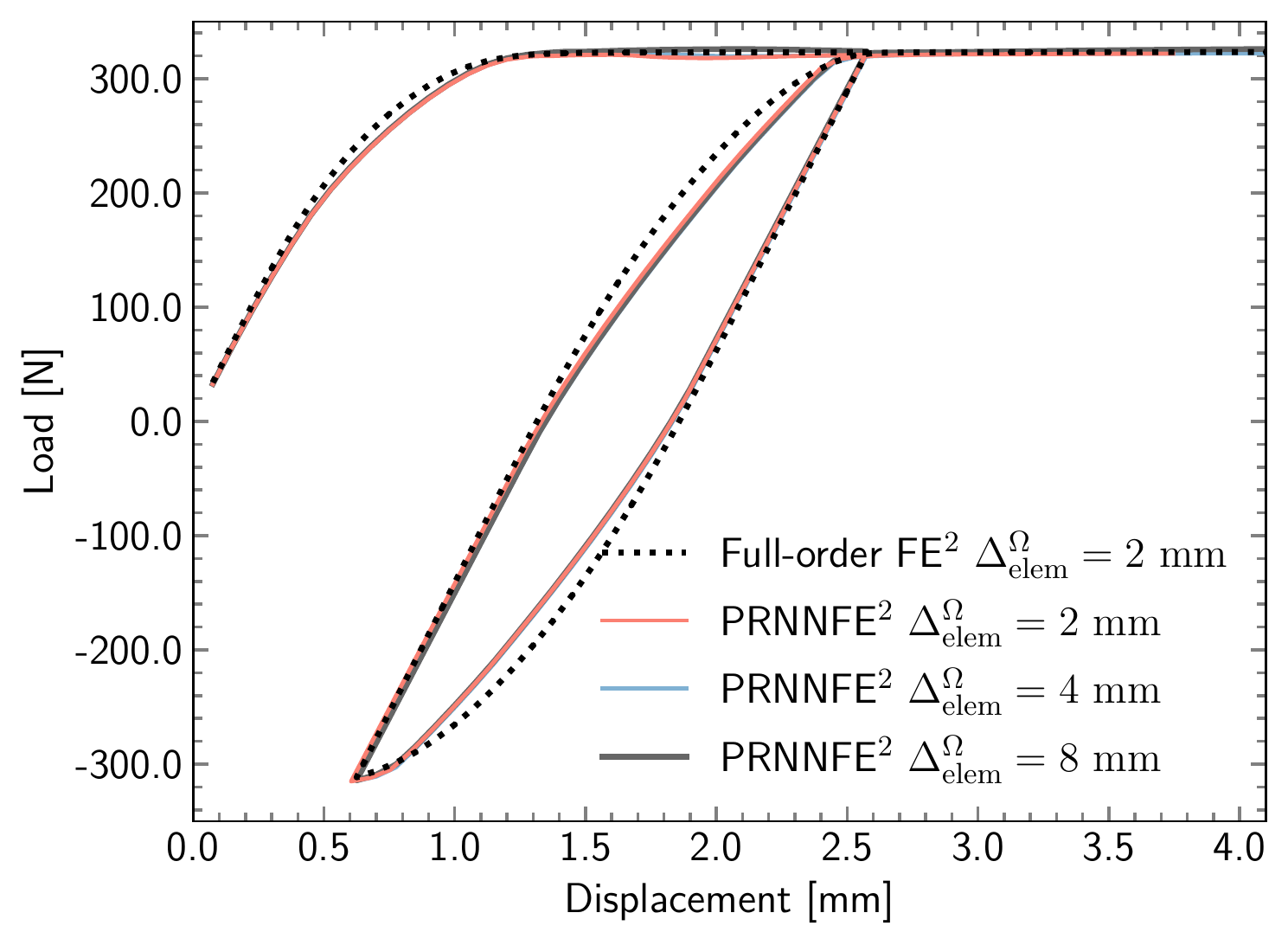}}
\hfill
\subfloat[Average RMSE at each time step of analysis with $\Delta_{\mathrm{elem}}^{\Omega} =$ 8 mm]{\label{fig:dogboneerror}\includegraphics[width=0.45\textwidth]{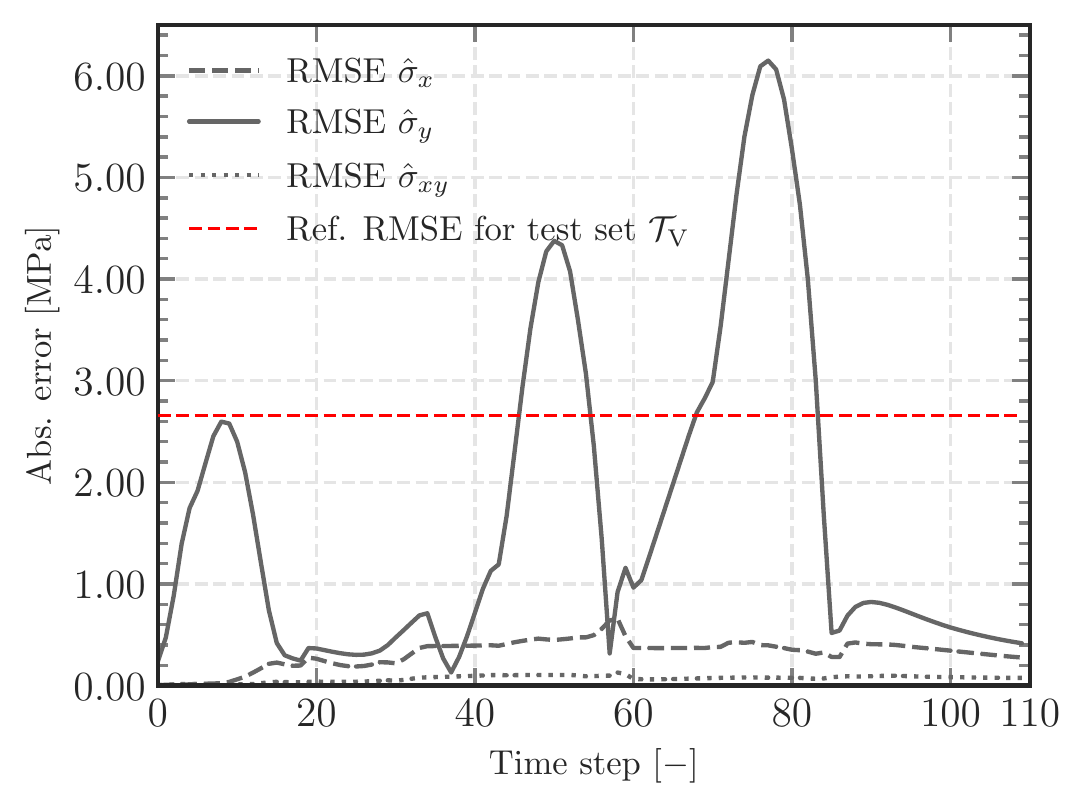}}
\caption{Tapered bar \fesqr \ example with (a) load-displacement curve with the full-order solution and best PRNN and (b) average error of PRNN's predictions at each time step of the analysis with $\Delta_{\mathrm{elem}}^{\Omega} = $ 8 mm.}
\end{figure}

The accuracy of the PRNN is further assessed by inspecting the load-displacement curve at the top edge of the bar. Fig. \ref{fig:dogbonecurve} shows the load-displacement curve using the full-order solution and the network’s response for different macroscopic mesh discretizations. For the refinement studied so far ($\Delta_{\mathrm{elem}}^{\Omega} = \ $ 8 mm), it is clear that the proposed network can capture accurately the entire nonlinear response, albeit with minor deviations when the tapered bar changes from tension to compression and then back again to tension. 

Next, the global accuracy of the method is verified. Since the surrogate model is not as accurate as the full-order model, a different equilibrium solution at a certain time step affects the equilibrium in the following loading steps, leading to accumulated error and diverging $\strainstresscurve$ paths. In this case, since no reduction in the step size was observed at any moment, the simple average error between the PRNN prediction and the full-order solution is calculated for each time step averaged over all integrations points at the macroscale, as illustrated in Fig. \ref{fig:dogboneerror}. For most of the simulation, the average absolute error in the predictions remains below 1 MPa with two peaks around 4 MPa and 6 MPa when the loading is reversed, following the trends observed in Fig. \ref{fig:dogbonecurve}. For reference, the lowest error of the network for test set $\mathcal{T}_{\mathrm{V}}$ is plotted.  
 \begin{table}[th]
\centering
\caption{\label{tab:dogbonecputimes} Computational cost for different mesh discretizations and efficiency of network in \fesqr \ approach.}
\begin{tabular}{lllll}
\hline\noalign{\smallskip}
& Macroscale element size ($\Delta_{\mathrm{elem}}^{\Omega}$) [mm] & 8 & 4 & 2 \\
& Number of elements at the macroscale & 64 & 134 & 454 \\ \noalign{\smallskip}\hline\noalign{\smallskip}
\multirow{3}{*}{Online}  & \fesqr \ wall-clock time  [s] & 21574 & 47644 & 178800 \\ 
& PRNN\fesqr \ wall-clock time [s] & 0.81 & 1.55 & 8.41 \\
& Speed-up$^a$ [-] & 26560 & 30746 & 21526   
          \\ \noalign{\smallskip} \cline{2-5} \noalign{\smallskip} 
\multirow{2}{*}{Offline} & Av. wall-clock time per curve (dataset gen.) [s] & 265$^b$ & N/A & N/A  \\
& Av. training time (excl. dataset gen.) [s] & 38045$^b$ & N/A & N/A \\
\noalign{\smallskip}\hline 
\multicolumn{4}{l}{\footnotesize $^a$Evaluated as \fesqr \ wall-clock time/PRNN\fesqr \ wall-clock time and averaged over 5 runs} \\
\multicolumn{4}{l}{\footnotesize $^b$One-off cost regardless of macroscopic mesh discretization} \\
\end{tabular}
\end{table}

For the purpose of assessing the efficiency of the network in accelerating the \fesqr \ simulations, three different levels of mesh refinement of the tapered bar are taken into account, being the coarses the one shown in Fig. \ref{fig:dogbone}. Table \ref{tab:dogbonecputimes} summarizes the wall-clock time spent in the analysis of the different discretizations, as well as the speed-up in comparison to the full-order solution. For the mesh used to illustrate this section ($\Delta_{\mathrm{elem}}^{\Omega} =$ 8 mm), replacing the solution of the BVP of the micromodel with the network led to a speed-up over 26000 with the accuracy reported in Fig. \ref{fig:dogboneerror}. Considering the offline costs, the training time is still lower than that of using the full-order solution with 134 macroscopic elements, which is a very modest number of elements for a multiscale problem.

Since the network is trained to replace the solution of the microscopic model, no additional training is required for the analysis of more complex cases where the macroscale problems require more elements and time steps. Hence, in general, higher speed-ups should be achieved with denser meshes. However, in this particular problem, this is not always the case. An increase from the coarsest to the intermediary mesh is observed, but no gain is achieved when refining even further. In that case, the reduction in performance  due to the higher number of iterations caused by the necessity of adaptively reducing the step size in order to ensure convergence. In constrast, the full-order solution was more numerically stable for this mesh density and the adaptive-stepping scheme was not triggered. 
\begin{figure}[!ht]
\centering
\subfloat[Macroscopic strains in $x$ and $y$]{\label{fig:dogbonehistogramxy}\includegraphics[width=0.33\textwidth]{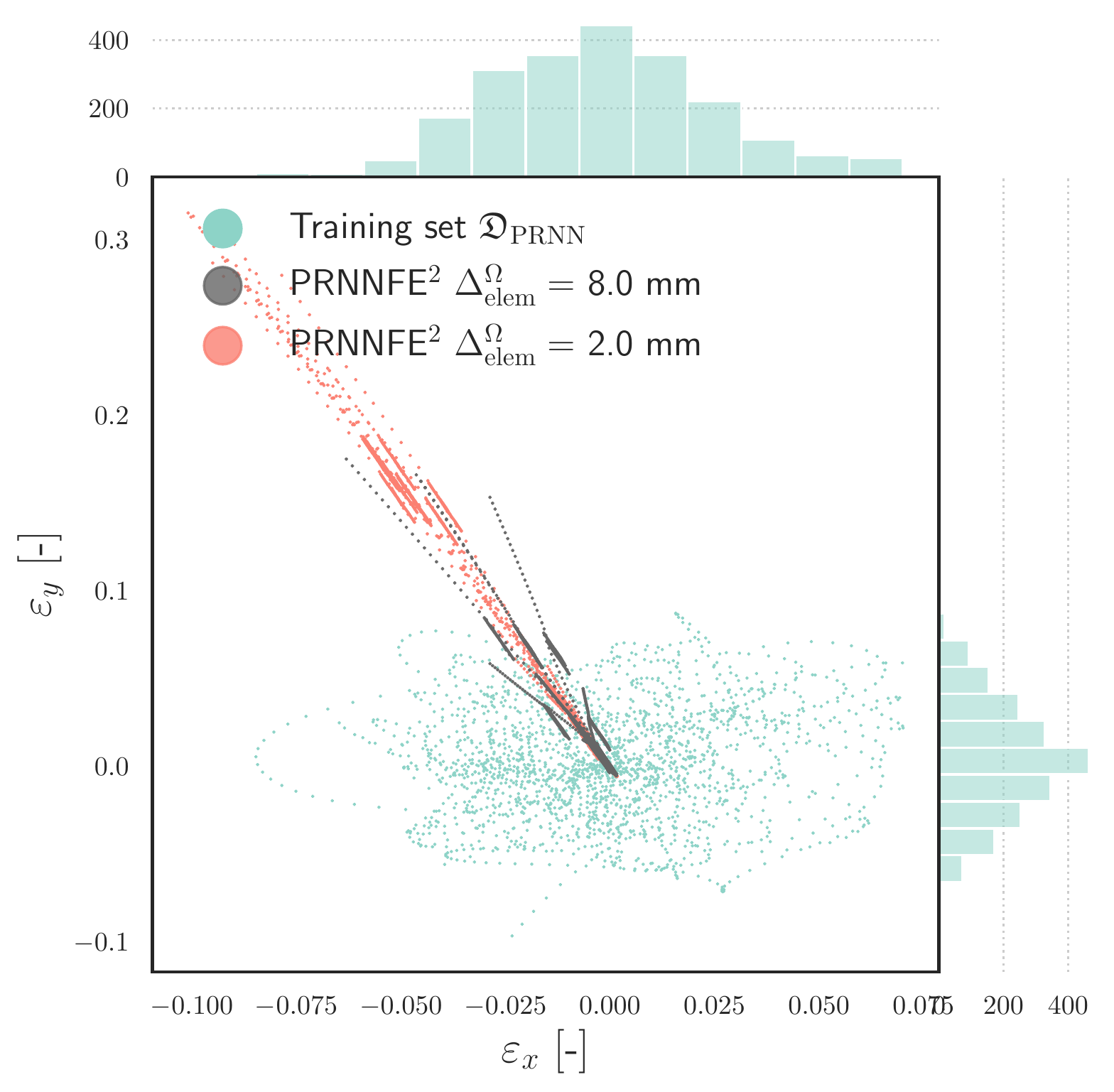}}
\subfloat[Macroscopic strains in $x$ and $xy$]{\label{fig:dogbonehistogramxxy}\includegraphics[width=0.33\textwidth]{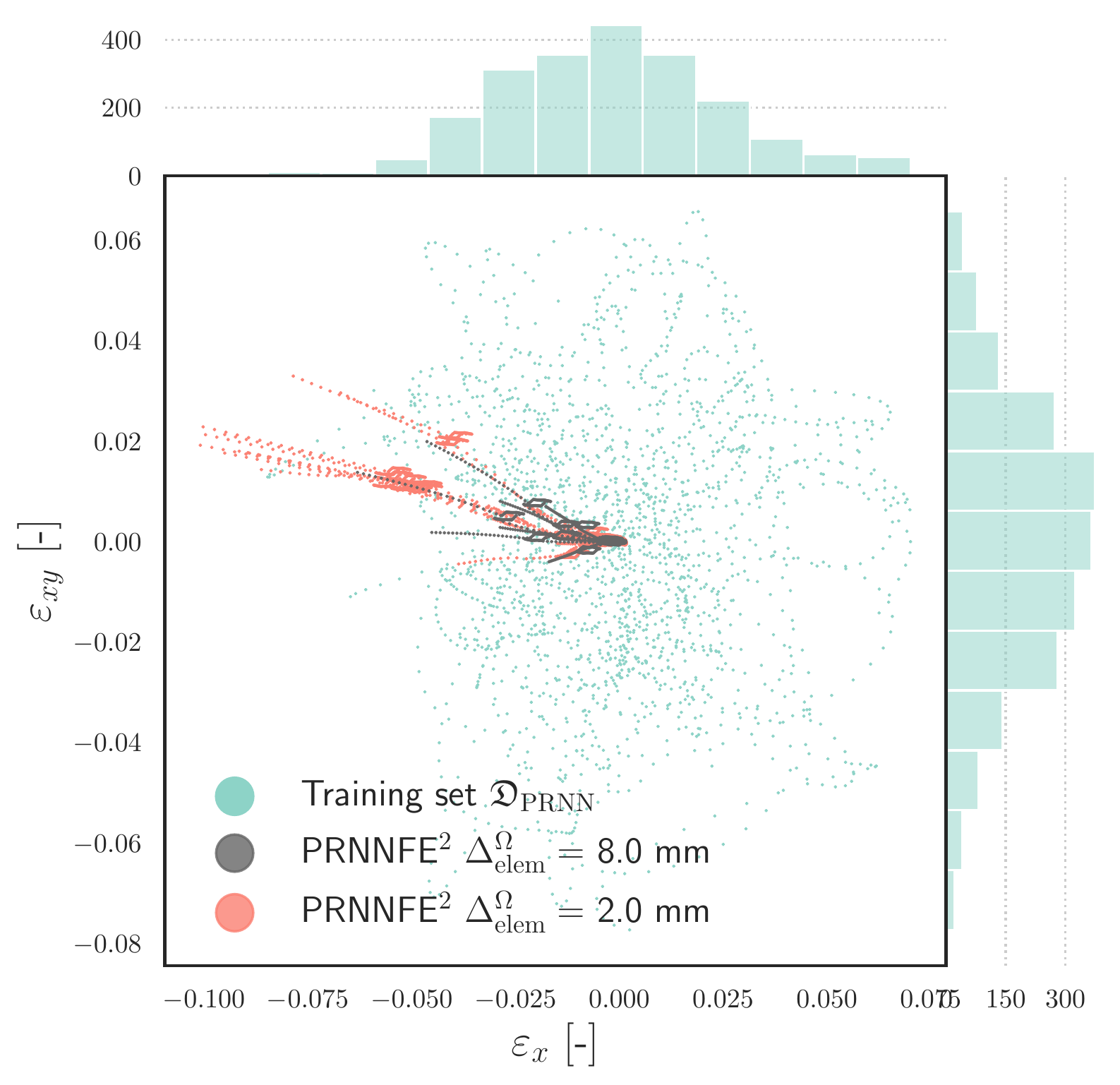}}
\subfloat[Macroscopic strains in $xy$ and $y$]{\label{fig:dogbonehistogramyxy}\includegraphics[width=0.33\textwidth]{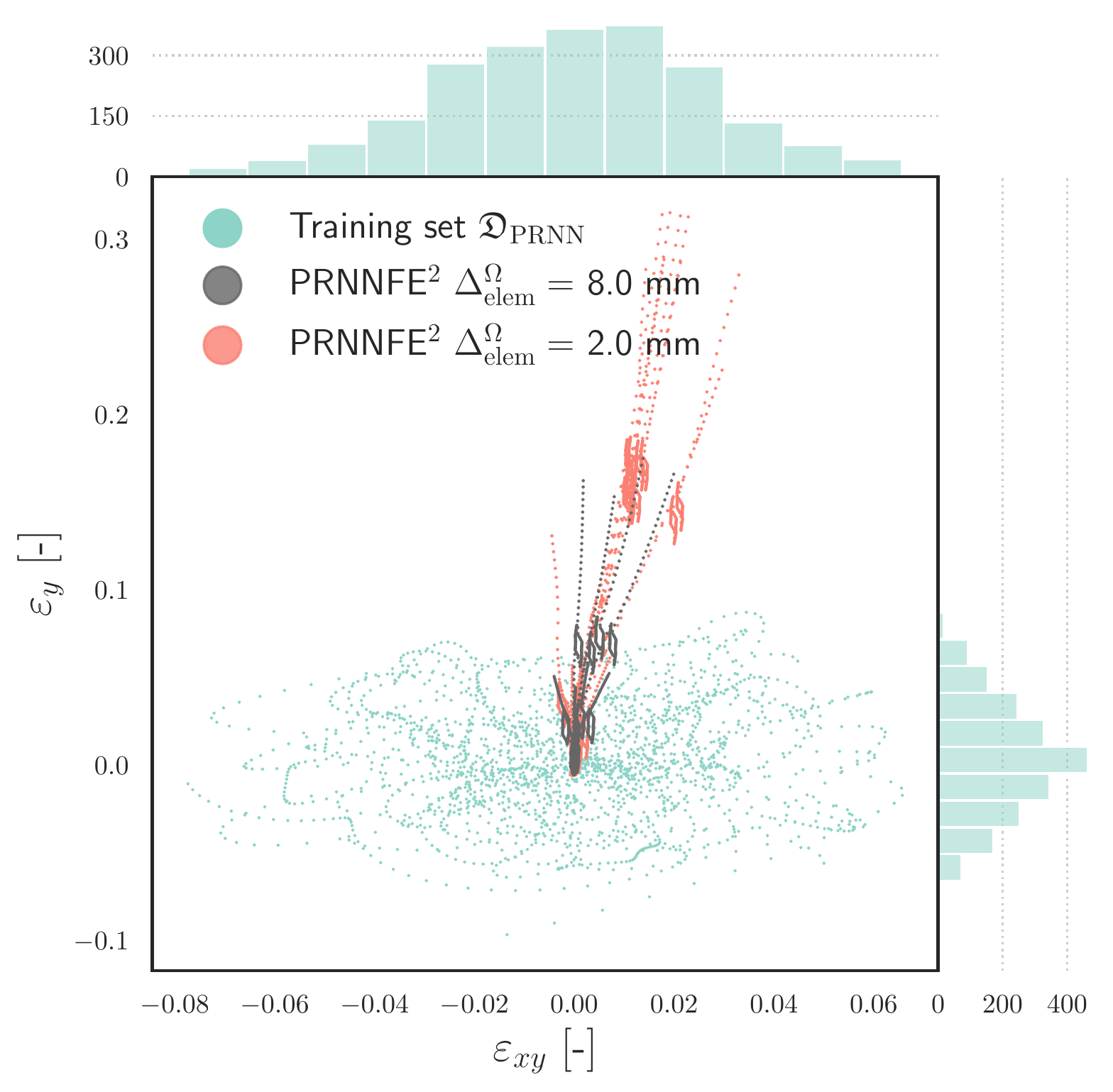}}
\caption{Joint distribution of strains from the training set of the best PRNN and strain distribution obtained by PRNNFE$^2$ for the tapered bar problem with different macroscopic mesh discretizations.}
\label{fig:dbhistogram} 
\end{figure}

As the mesh is refined and strain localization takes place (see the red region in Fig. \ref{fig:dogbonedispfield}), even higher strain levels are achieved, pushing the network to make predictions in unexplored regions during training, as illustrated in Fig. \ref{fig:dbhistogram}. Note that the network is already making far-reaching predictions in the coarsest discretization, although in a less extensive way. In the mesh with $\Delta_{\mathrm{elem}}^{\Omega} = \ $8 mm, the maximum strain does not exceed 0.2, while $\Delta_{\mathrm{elem}}^{\Omega} = \ $ 2 mm leads to strains higher than 0.3. In spite of these complicating aspects, it is worth mentioning this is still a significant speed-up. Moreover, a far less severe effect on the global accuracy is observed, as illustrated by the almost overlapping load-displacement curves in Fig. \ref{fig:dogbonecurve}. In that sense, the adaptive-stepping scheme plays an important role to help overcome convergence issues. 

\subsection{Plate with multiple holes}
\begin{figure}[!ht]
\centering
\includegraphics[width=0.76\textwidth]{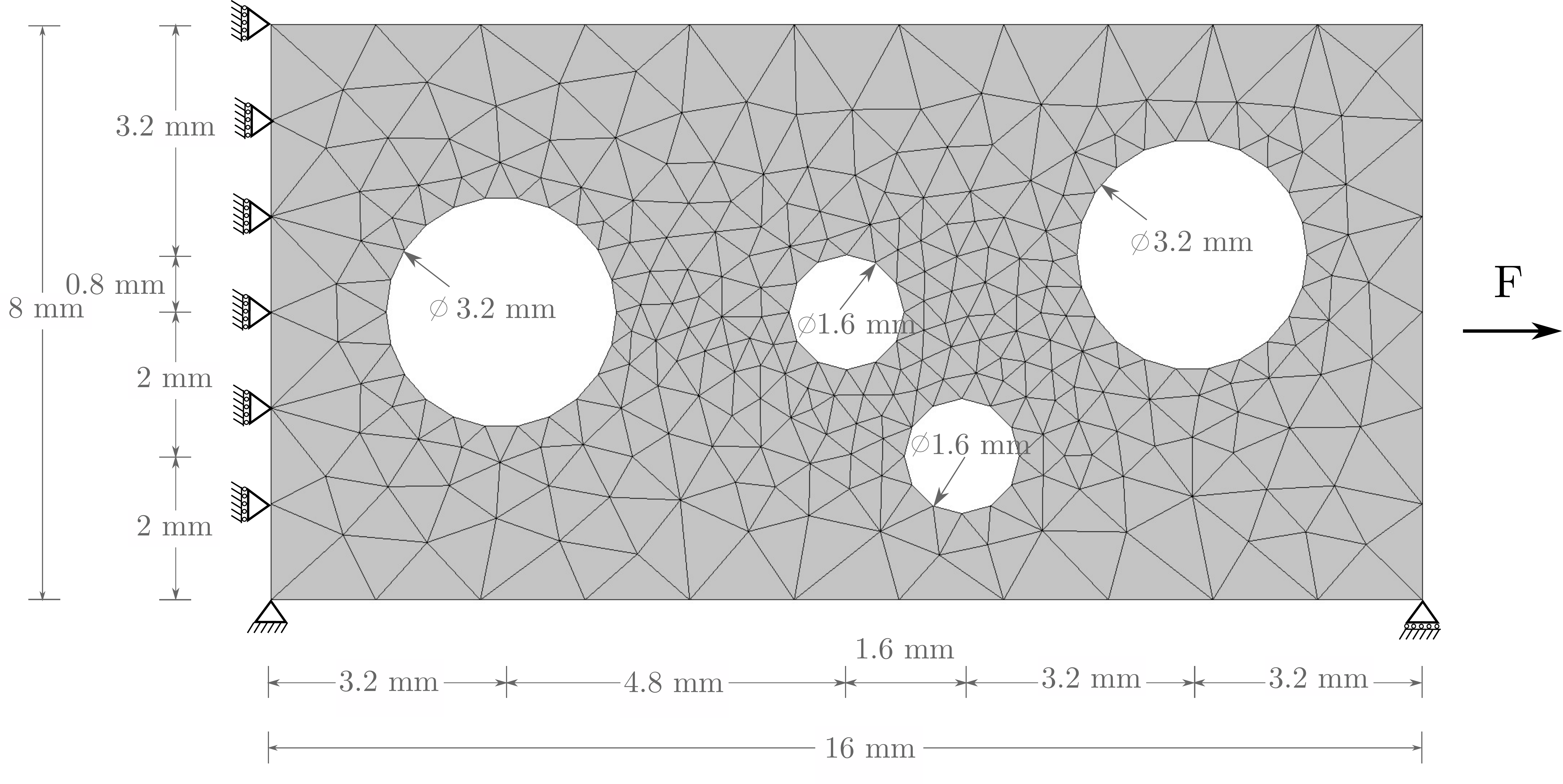}
\caption{Plate with cutouts: geometry and boundary and loading conditions.}
\label{fig:cutoutgeom} 
\end{figure}

\begin{figure}[!ht]
\centering
\label{fig:cutout} 
\subfloat[Load-displacement curve]{\label{fig:cutoutcurve}\includegraphics[width=0.45\textwidth]{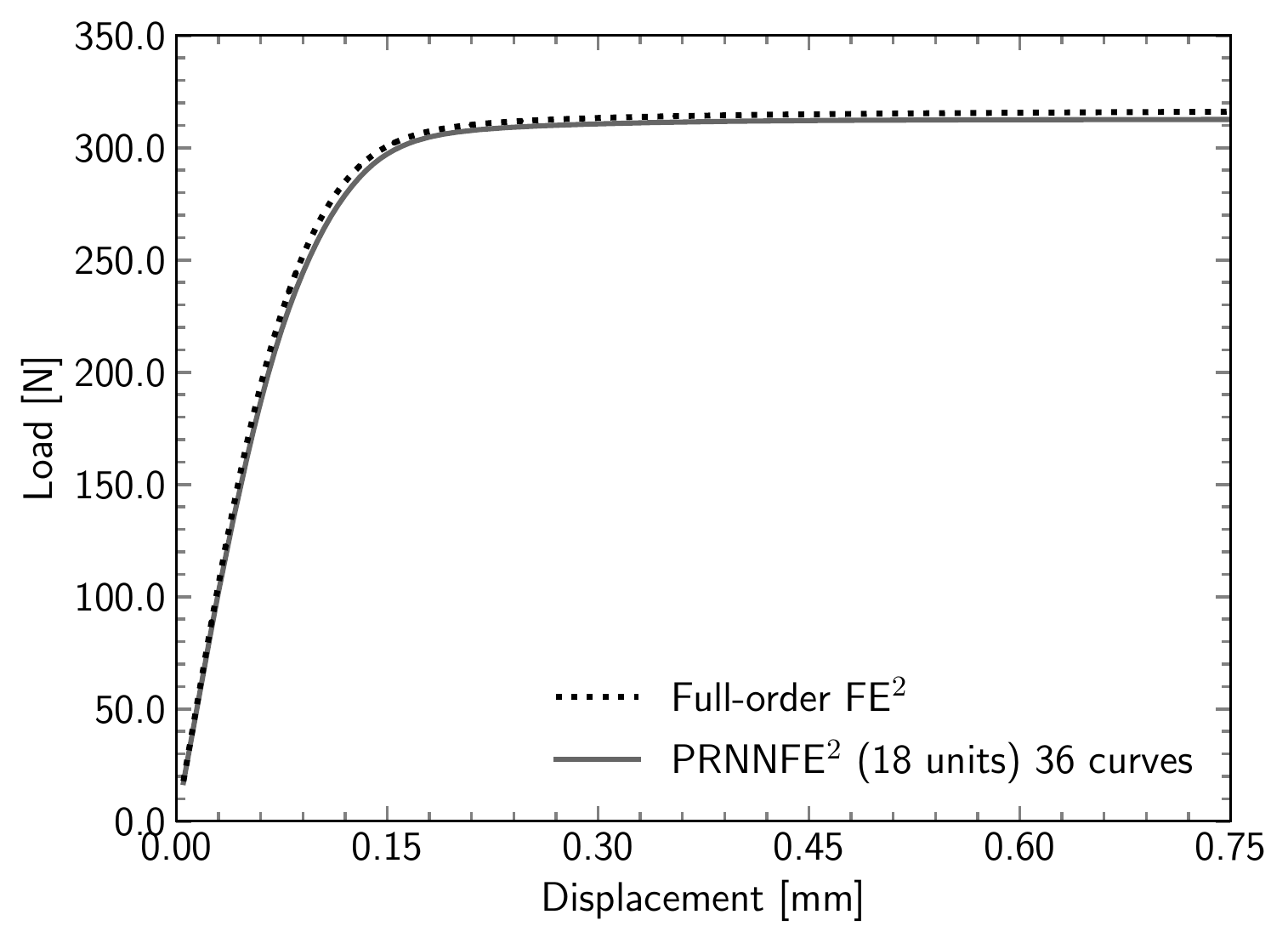}}
\hfill
\subfloat[Average RMSE]{\label{fig:cutouterror}\includegraphics[width=0.45\textwidth]{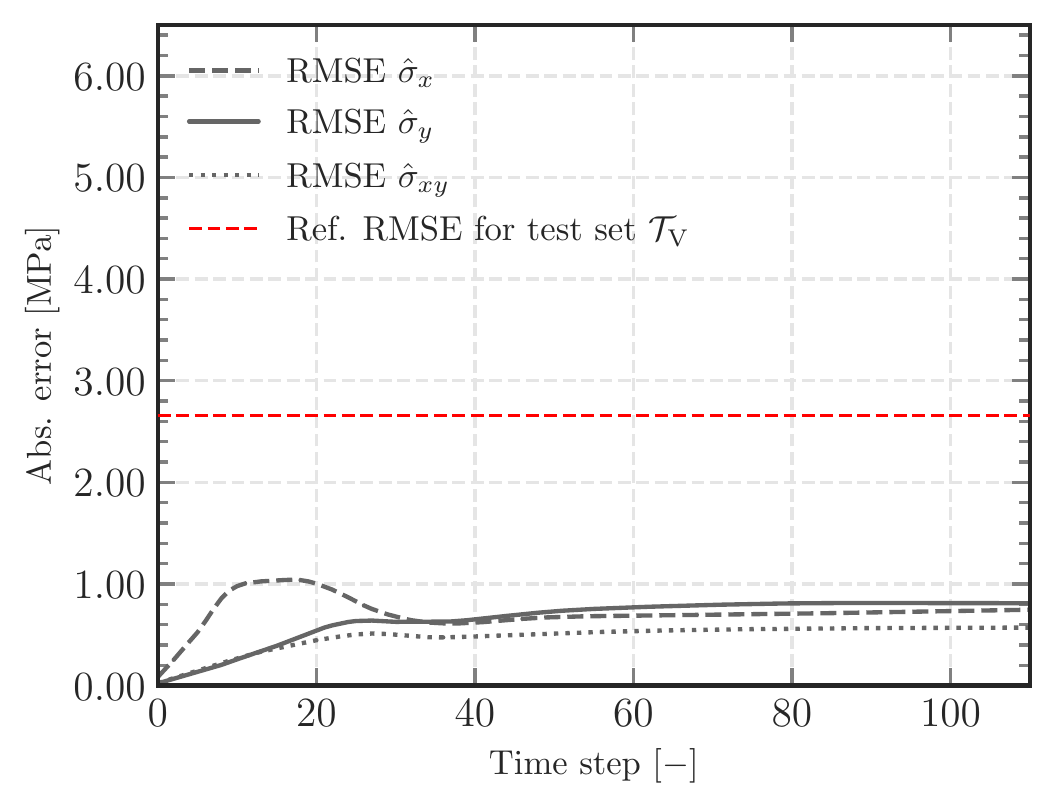}}
\caption{Plate with cutouts: (a) load-displacement curve and (b) average error of PRNN's predictions at each time step of the analysis.}
\end{figure}

As a final example, a composite plate with multiple cutouts with geometry and boundary and loading conditions as illustrated in Fig. \ref{fig:cutoutgeom} is studied. Again, an \fesqr \ approach is employed to solve the problem for the same microscopic model with which all the networks in Section \ref{sec:results} were trained for. This time, no unloading is imposed and 150 load steps with $\Delta$s $= 5.0 \times 10^{-3}$ are considered. 
The load-displacement curve at the right edge of the plate is plotted in Fig. \ref{fig:cutoutcurve} using both the full-order solution and the network. Again, good agreement is observed between the macroscopic responses. The slight inaccuracy between those are quantified in Fig. \ref{fig:cutouterror}, in which the average absolute error of the component with the highest magnitudes ($\sigma_x$) is  around 1 MPa for almost the entire simulation. 

\begin{figure}[ht!]
\centering
\subfloat[Displacement field at the end of the analysis with PRNN and selection of integration points for inspection]{\label{fig:cutoutlasttimestep}\includegraphics[width=0.58\textwidth, valign = c]{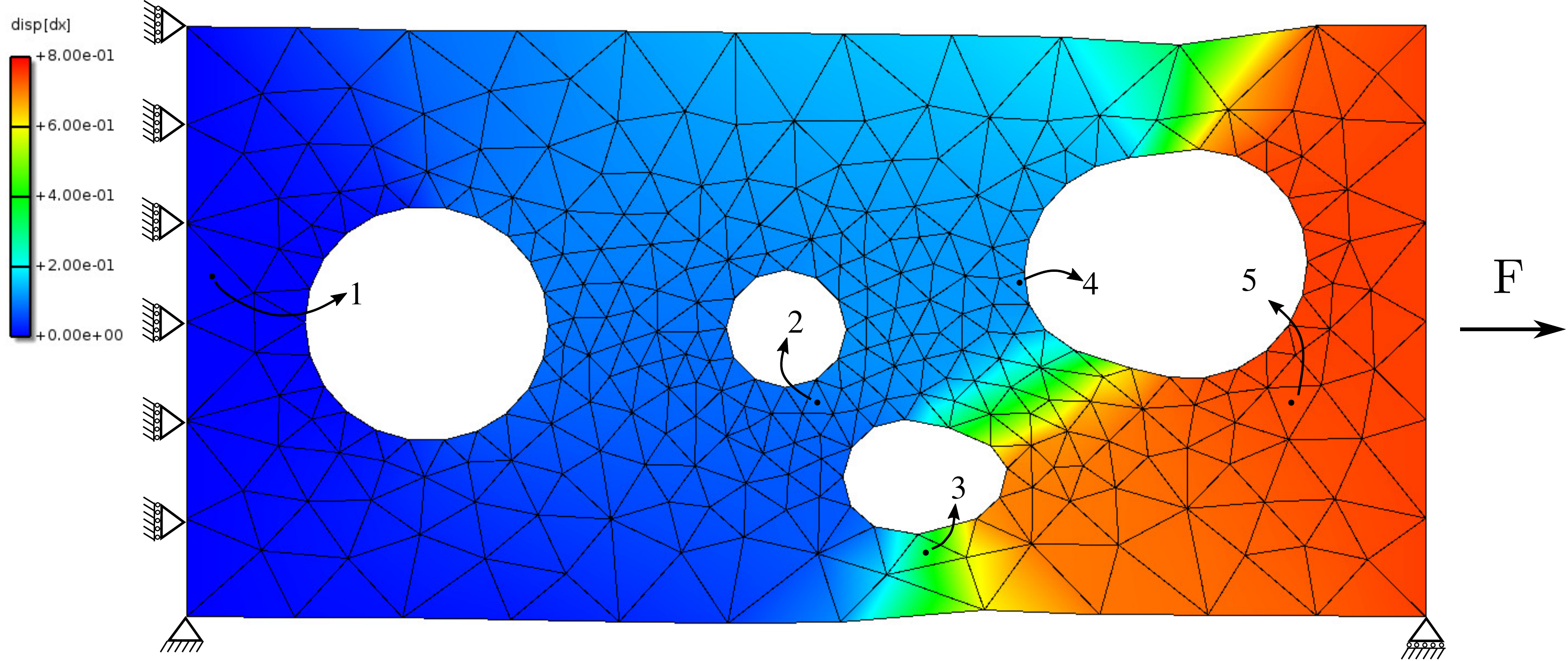}}
\hfill
\subfloat[Stress-time view of macroscopic integration points. Full-order solution in black and network's prediction in gray]{\label{fig:cutoutcurves}\includegraphics[width=0.40\textwidth, valign = c]{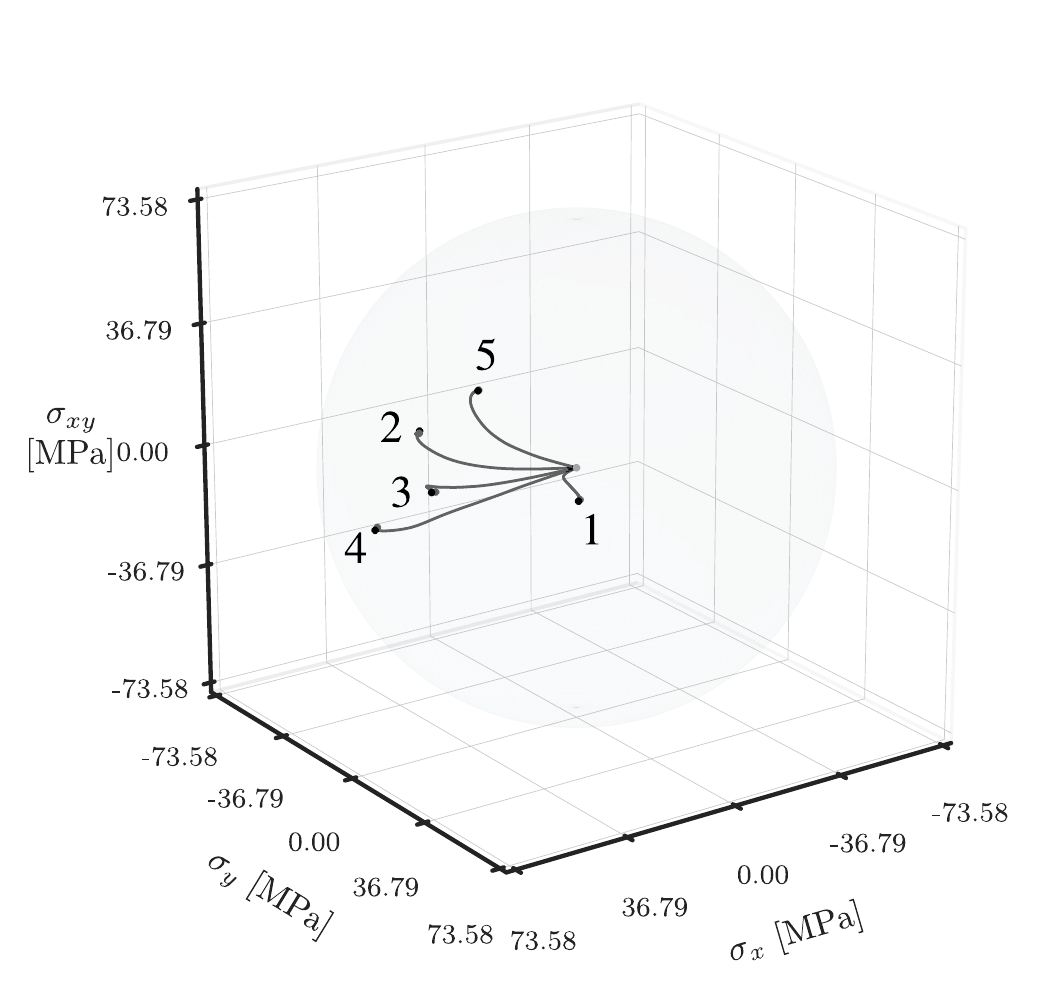}}
\caption{Plate with cutouts: (a) displacement field at the end of the analysis and (b) selected integration points shown in stress-time view.}
\label{fig:cutoutexample}
\end{figure}

The displacement field at the end of the analysis is shown in Fig. \ref{fig:cutoutlasttimestep}, where the location of five macroscopic integration points are marked for further inspection. The stress paths for each of these points are illustrated in Fig. \ref{fig:cutoutcurves}, where the full-order solution and the network prediction are plotted in black and gray, respectively. Note how the stress paths are non-proportional even for the relatively simple loading condition observed in the macroscale. The integration point on the edge of one of the cutouts, namely point 4, is also the one with the highest stress magnitude and the closest to a uniaxial state in the $x$ direction while the other points experience multiaxial loading more strongly. 

In terms of efficiency, the solution using the network is around 26705 times faster than the full-order solution, which took approximately 258780 s (around 72 hours). The order of magnitude in the speed-up is similar to that obtained in the tapered bar problem. Although no additional offline costs are incurred because the network has been trained before for the same microscopic model and it does not depend on the macroscopic problem at hand, it is worth stressing that the runtime of the full-order solution exceeds the sum of the online and offline costs of the PRNN. 

Finally, this example shows that the network can capture multiaxial stress states and non-proportional loading as obtained in \fesqr \ simulations accurately. No convergence issues were encountered in the PRNN\fesqr \ simulation which points to the smoothness of the predictions that is not always guaranteed with surrogate models (see \textit{e.g.} RNN curves in Fig. \ref{fig:samplecurvelevel3ongp}). 

\section{Concluding remarks}
\label{sec:conclusion}

In this paper, a novel network with embedded physics-based constitutive models is proposed as surrogate model for the behavior of path-dependent materials in \fesqr \ simulations. The central idea is to address common problems in modeling path-dependent materials using black-box models (\textit{e.g.} unique mapping between input and output and limited extrapolation abilities) by taking a step back and reintroducing physics into the network in a way that requires very little extra coding effort with respect to existing \fesqr \ frameworks. This is done by employing the same material models used for the microscopic level as part of one of the layers of the network.

To accommodate this non-standard neural layer the following changes with respect to standard neural network architectures are proposed. First, neurons are assembled in groups of the size of the number of strain/stress components of the problem. These are referred to as fictitious material points. Secondly, to take advantage of all the information coming from the physics-based material model, we store the updated internal variables used to fully describe the state of the fictitious material point in an auxiliary vector. With that, when new strain values are fed, the material point will start from the last stored internal variables. As a consequence, each subgroup follows a unique path without the need to increase the feature space with extra history variables. 

The properties and assumptions made by the physics-based constitutive model are inherited by the network and play a major role in reducing the amount of data required to mirror physical and complex behaviors such as elastic unloading/reloading. Here, the decomposition of the strain in elastic and plastic parts is an assumption built in the material model used to describe the nonlinear microscopic material phase and is also observed in the network when the local stresses of the fictitious material points are evaluated. This simple but highly-flexible arrangement allows the network to capture arbitrary unloading behaviors with only monotonic data, a stark contrast with other popular models such as RNNs. The PRNN inherits from \fesqr \ the idea that complex behavior of heterogeneous materials can be accurately described by letting simpler constitutive models that represent the microscopic constituents interact. The difference is that the interaction between the constituents is not based on micromechanics diretly but learned from data obtained from micromodel simulations.

Based on that, an extensive numerical comparison involving a state-of-the-art black-box model, namely a Bayesian Recurrent Neural Network (referred to in this work as RNN), was carried out in order to elucidate the abilities of the proposed network (referred to as PRNN). First, we trained both networks only on 18 monotonic curves with known directions and proportional loading in a similar fashion as done to calibrate classical mesomodels. Such strategy led to poor performance when trying to predict other random directions from the RNN, but good accuracy from our method (Fig. \ref{fig:convtrlvl2predlvl2}). Following that conclusion, the size of the RNN's training dataset was sequentially increased until the addition of more data did not result in significant gains in accuracy. At that stage, the PRNN performed with the same level of accuracy but with a factor of 16 times less data. 

Next, both networks were used to predict non-monotonic loading. For that scenario, the PRNN showed the same level of accuracy as before with the same minimal training dataset (Fig. \ref{fig:convtrlvl2predlvl3}). Such outstanding result is not observed in the RNNs, which again required a larger training set. This time, non-monotonic loading curves were added until the RNN's accuracy could no longer be significantly improved. As a result, a 32 times larger training dataset in comparison to the one used to train our network was necessary. Furthermore, while our network continues to perform well in all the scenarios tested so far, two other situations exposed the pitfalls of RNNs: (i) when trying to predict unloading in a different location than the one seen in training (Fig. \ref{fig:convlevel3onfar}) and (ii) when the step size was modified (Fig. \ref{fig:convlevel3onstep}). This is typically tackled by sampling different unloading behaviors with different step sizes, leading to the choice of arbitrarily long sequences. However, we showed that this is a trivial scenario for the PRNN. The network is only as sensitive to step size as the material models it includes. 

In the last test, both networks were used to predict non-proportional and non-monotonic paths and neither succeeded, although they failed at very different levels. While the lowest error of the RNN was around 32 MPa, the best PRNN led to an error around 9 MPa error (Fig. \ref{fig:convlevel3ongp}). Based on that, a second approach to generate the dataset was considered. Random strain paths were generated from Gaussian Processes priors, which produces non-proportional and non-monotonic loading as opposed to the proportional loading previously considered for training. This time, the size of the training dataset and the type of loading was also a variable for the PRNN. It was found that training the PRNN on random non-proportional and non-monotonic curves yields higher accuracy than training with known, proportional, and monotonic curves for all loading scenarios (Fig. \ref{fig:convergencediftr}). Although training with known directions is appealing, having a network that provides lower errors and consistent performance with random directions is also interesting. Ultimately, the PRNN consistently outperformed the RNN with 64 times less data.  

Finally, after ensuring the PRNN capacity in several challenging scenarios for black-box models, one of the networks trained on non-proportional and non-monotonic curves was chosen to surrogate the microscopic model in a set of two \fesqr \ examples. The first example concerned a tapered bar in transverse tension and was used to illustrate how the different material models in the PRNN behave for a single macroscopic integration point (Fig. \ref{fig:dogbone}). For different discretizations, speed-ups between 21000 and 31000 were obtained for the online phase (Table \ref{tab:dogbonecputimes}). Such substantial efficiency gain is explained from the dramatically reduced number of material model calls and the bypassing of solving the nonlinear microscopic system of equations for macroscopic stress evaluation for a single load step of each macroscopic integration point. 

In the last example, a similar order of magnitude of speed-up was observed and the accuracy of the PRNN was illustrated by comparing the $\strainstresscurve$ paths of different macroscopic integration points of a plate with multiple cutouts subjected to tension (Fig. \ref{fig:cutoutcurves}). For the analyzed cases, the time needed for a single \fesqr \ analysis exceeded the total offline and online time for the PRNN analysis, even though the selected problems had a very modest number of macroscopic elements. Moreover, performing subsequent macroscale analysis with the same material would not require any additional offline training and would therefore leverage the complete speed-up of four orders of magnitude.

\section*{Acknowledgements}
The authors acknowledge the TU Delft AI Initiative for their support through the SLIMM AI Lab. FM acknowledges financial support from the Dutch Research Council (NWO) under Vidi grant 16464.

\bibliographystyle{elsarticle-num-names} 
\bibliography{references.bib}

\begin{thebibliography}{38}
\expandafter\ifx\csname natexlab\endcsname\relax\def\natexlab#1{#1}\fi
\providecommand{\url}[1]{\texttt{#1}}
\providecommand{\href}[2]{#2}
\providecommand{\path}[1]{#1}
\providecommand{\DOIprefix}{doi:}
\providecommand{\ArXivprefix}{arXiv:}
\providecommand{\URLprefix}{URL: }
\providecommand{\Pubmedprefix}{pmid:}
\providecommand{\doi}[1]{\href{http://dx.doi.org/#1}{\path{#1}}}
\providecommand{\Pubmed}[1]{\href{pmid:#1}{\path{#1}}}
\providecommand{\bibinfo}[2]{#2}
\ifx\xfnm\relax \def\xfnm[#1]{\unskip,\space#1}\fi
\bibitem[{Goury et~al.(2016)Goury, Amsallem, Bordas, Liu, and
  Kerfriden}]{Gouryetal2016}
\bibinfo{author}{O.~Goury}, \bibinfo{author}{D.~Amsallem},
  \bibinfo{author}{S.~Bordas}, \bibinfo{author}{W.~Liu},
  \bibinfo{author}{P.~Kerfriden},
\newblock \bibinfo{title}{Automatised selection of load paths to construct
  reduced-order models in computational damage micromechanics: from
  dissipation-driven random selection to bayesian optimization},
\newblock \bibinfo{journal}{Computational Mechanics} \bibinfo{volume}{58}
  (\bibinfo{year}{2016}) \bibinfo{pages}{213--234}.
  \DOIprefix\doi{10.1007/s00466-016-1290-2}.
\bibitem[{Rocha et~al.(2021)Rocha, Kerfriden, and {van der Meer}}]{Rocha2021}
\bibinfo{author}{I.~B. C.~M. Rocha}, \bibinfo{author}{P.~Kerfriden},
  \bibinfo{author}{F.~P. {van der Meer}},
\newblock \bibinfo{title}{{On-the-fly construction of surrogate constitutive
  models for concurrent multiscale mechanical analysis through probabilistic
  machine learning}},
\newblock \bibinfo{journal}{Journal of Computational Physics: X}
  \bibinfo{volume}{9} (\bibinfo{year}{2021}) \bibinfo{pages}{1--30}.
  \DOIprefix\doi{10.1016/j.jcpx.2020.100083}.
  \href{http://arxiv.org/abs/2007.07749}{{\tt arXiv:2007.07749}}.
\bibitem[{Vlassis and Sun(2021)}]{VLASSIS2021}
\bibinfo{author}{N.~N. Vlassis}, \bibinfo{author}{W.~Sun},
\newblock \bibinfo{title}{Sobolev training of thermodynamic-informed neural
  networks for interpretable elasto-plasticity models with level set
  hardening},
\newblock \bibinfo{journal}{Computer Methods in Applied Mechanics and
  Engineering} \bibinfo{volume}{377} (\bibinfo{year}{2021})
  \bibinfo{pages}{113695}.
  \DOIprefix\doi{https://doi.org/10.1016/j.cma.2021.113695}.
\bibitem[{Lefik et~al.(2009)Lefik, Boso, and Schrefler}]{LEFIK20091785}
\bibinfo{author}{M.~Lefik}, \bibinfo{author}{D.~Boso},
  \bibinfo{author}{B.~Schrefler},
\newblock \bibinfo{title}{Artificial neural networks in numerical modelling of
  composites},
\newblock \bibinfo{journal}{Computer Methods in Applied Mechanics and
  Engineering} \bibinfo{volume}{198} (\bibinfo{year}{2009})
  \bibinfo{pages}{1785--1804}.
  \DOIprefix\doi{https://doi.org/10.1016/j.cma.2008.12.036},
  \bibinfo{note}{advances in Simulation-Based Engineering Sciences – Honoring
  J. Tinsley Oden}.
\bibitem[{Huang et~al.(2020)Huang, Fuhg, Wei{\ss}enfels, and
  Wriggers}]{Huang2020}
\bibinfo{author}{D.~Huang}, \bibinfo{author}{J.~N. Fuhg},
  \bibinfo{author}{C.~Wei{\ss}enfels}, \bibinfo{author}{P.~Wriggers},
\newblock \bibinfo{title}{{A machine learning based plasticity model using
  proper orthogonal decomposition}},
\newblock \bibinfo{journal}{Computer Methods in Applied Mechanics and
  Engineering} \bibinfo{volume}{365} (\bibinfo{year}{2020})
  \bibinfo{pages}{113008}. \DOIprefix\doi{10.1016/j.cma.2020.113008}.
  \href{http://arxiv.org/abs/2001.03438}{{\tt arXiv:2001.03438}}.
\bibitem[{Ghavamian and Simone(2019)}]{Ghavamian2019}
\bibinfo{author}{F.~Ghavamian}, \bibinfo{author}{A.~Simone},
\newblock \bibinfo{title}{{Accelerating multiscale finite element simulations
  of history-dependent materials using a recurrent neural network}},
\newblock \bibinfo{journal}{Computer Methods in Applied Mechanics and
  Engineering} \bibinfo{volume}{357} (\bibinfo{year}{2019})
  \bibinfo{pages}{112594}. \DOIprefix\doi{10.1016/j.cma.2019.112594}.
\bibitem[{Wu et~al.(2020)Wu, Nguyen, Kilingar, and Noels}]{WU2020113234}
\bibinfo{author}{L.~Wu}, \bibinfo{author}{V.~D. Nguyen}, \bibinfo{author}{N.~G.
  Kilingar}, \bibinfo{author}{L.~Noels},
\newblock \bibinfo{title}{A recurrent neural network-accelerated multi-scale
  model for elasto-plastic heterogeneous materials subjected to random cyclic
  and non-proportional loading paths},
\newblock \bibinfo{journal}{Computer Methods in Applied Mechanics and
  Engineering} \bibinfo{volume}{369} (\bibinfo{year}{2020})
  \bibinfo{pages}{113234}.
  \DOIprefix\doi{https://doi.org/10.1016/j.cma.2020.113234}.
\bibitem[{Mozaffar et~al.(2019)Mozaffar, Bostanabad, Chen, Ehmann, Cao, and
  Bessa}]{Mozaffar26414}
\bibinfo{author}{M.~Mozaffar}, \bibinfo{author}{R.~Bostanabad},
  \bibinfo{author}{W.~Chen}, \bibinfo{author}{K.~Ehmann},
  \bibinfo{author}{J.~Cao}, \bibinfo{author}{M.~A. Bessa},
\newblock \bibinfo{title}{Deep learning predicts path-dependent plasticity},
\newblock \bibinfo{journal}{Proceedings of the National Academy of Sciences}
  \bibinfo{volume}{116} (\bibinfo{year}{2019}) \bibinfo{pages}{26414--26420}.
  \DOIprefix\doi{10.1073/pnas.1911815116}.
  \href{http://arxiv.org/abs/https://www.pnas.org/content/116/52/26414.full.pdf}{{\tt
  arXiv:https://www.pnas.org/content/116/52/26414.full.pdf}}.
\bibitem[{Chen(2021)}]{Chen2021}
\bibinfo{author}{G.~Chen},
\newblock \bibinfo{title}{Recurrent neural networks ({RNNs}) learn the
  constitutive law of viscoelasticity},
\newblock \bibinfo{journal}{Computational Mechanics} \bibinfo{volume}{67}
  (\bibinfo{year}{2021}) \bibinfo{pages}{1009--1019}.
  \DOIprefix\doi{10.1007/s00466-021-01981-y}.
\bibitem[{Koeppe et~al.(2021)Koeppe, Bamer, Selzer, Nestler, and
  Markert}]{Koeppe2021}
\bibinfo{author}{A.~Koeppe}, \bibinfo{author}{F.~Bamer},
  \bibinfo{author}{M.~Selzer}, \bibinfo{author}{B.~Nestler},
  \bibinfo{author}{B.~Markert}, \bibinfo{title}{Explainable artificial
  intelligence for mechanics: physics-informing neural networks for
  constitutive models}, \bibinfo{year}{2021}.
  \href{http://arxiv.org/abs/2104.10683}{{\tt arXiv:2104.10683}}.
\bibitem[{Logarzo et~al.(2021)Logarzo, Capuano, and Rimoli}]{LOGARZO2021}
\bibinfo{author}{H.~J. Logarzo}, \bibinfo{author}{G.~Capuano},
  \bibinfo{author}{J.~J. Rimoli},
\newblock \bibinfo{title}{Smart constitutive laws: Inelastic homogenization
  through machine learning},
\newblock \bibinfo{journal}{Computer Methods in Applied Mechanics and
  Engineering} \bibinfo{volume}{373} (\bibinfo{year}{2021})
  \bibinfo{pages}{113482}.
  \DOIprefix\doi{https://doi.org/10.1016/j.cma.2020.113482}.
\bibitem[{Gorji et~al.(2020)Gorji, Mozaffar, Heidenreich, Cao, and
  Mohr}]{GORJI2020}
\bibinfo{author}{M.~B. Gorji}, \bibinfo{author}{M.~Mozaffar},
  \bibinfo{author}{J.~N. Heidenreich}, \bibinfo{author}{J.~Cao},
  \bibinfo{author}{D.~Mohr},
\newblock \bibinfo{title}{On the potential of recurrent neural networks for
  modeling path dependent plasticity},
\newblock \bibinfo{journal}{Journal of the Mechanics and Physics of Solids}
  \bibinfo{volume}{143} (\bibinfo{year}{2020}) \bibinfo{pages}{103972}.
  \DOIprefix\doi{https://doi.org/10.1016/j.jmps.2020.103972}.
\bibitem[{Raissi et~al.(2019)Raissi, Perdikaris, and Karniadakis}]{RAISSI2019}
\bibinfo{author}{M.~Raissi}, \bibinfo{author}{P.~Perdikaris},
  \bibinfo{author}{G.~Karniadakis},
\newblock \bibinfo{title}{Physics-informed neural networks: A deep learning
  framework for solving forward and inverse problems involving nonlinear
  partial differential equations},
\newblock \bibinfo{journal}{Journal of Computational Physics}
  \bibinfo{volume}{378} (\bibinfo{year}{2019}) \bibinfo{pages}{686--707}.
  \DOIprefix\doi{https://doi.org/10.1016/j.jcp.2018.10.045}.
\bibitem[{Haghighat et~al.(2021)Haghighat, Raissi, Moure, Gomez, and
  Juanes}]{HAGHIGHAT2021}
\bibinfo{author}{E.~Haghighat}, \bibinfo{author}{M.~Raissi},
  \bibinfo{author}{A.~Moure}, \bibinfo{author}{H.~Gomez},
  \bibinfo{author}{R.~Juanes},
\newblock \bibinfo{title}{A physics-informed deep learning framework for
  inversion and surrogate modeling in solid mechanics},
\newblock \bibinfo{journal}{Computer Methods in Applied Mechanics and
  Engineering} \bibinfo{volume}{379} (\bibinfo{year}{2021})
  \bibinfo{pages}{113741}.
  \DOIprefix\doi{https://doi.org/10.1016/j.cma.2021.113741}.
\bibitem[{Eghbalian et~al.(2022)Eghbalian, Pouragha, and
  Wan}]{Eghbalianetal2022}
\bibinfo{author}{M.~Eghbalian}, \bibinfo{author}{M.~Pouragha},
  \bibinfo{author}{R.~Wan}, \bibinfo{title}{A physics-informed deep neural
  network for surrogate modeling in classical elasto-plasticity},
  \bibinfo{year}{2022}. \DOIprefix\doi{10.48550/ARXIV.2204.12088}.
\bibitem[{Arora et~al.(2022)Arora, Kakkar, Dey, and Chakraborty}]{Arora2022}
\bibinfo{author}{R.~Arora}, \bibinfo{author}{P.~Kakkar},
  \bibinfo{author}{B.~Dey}, \bibinfo{author}{A.~Chakraborty},
  \bibinfo{title}{Physics-informed neural networks for modeling rate- and
  temperature-dependent plasticity}, \bibinfo{year}{2022}.
  \DOIprefix\doi{10.48550/ARXIV.2201.08363}.
\bibitem[{Masi et~al.(2021)Masi, Stefanou, Vannucci, and
  Maffi-Berthier}]{Masi2020}
\bibinfo{author}{F.~Masi}, \bibinfo{author}{I.~Stefanou},
  \bibinfo{author}{P.~Vannucci}, \bibinfo{author}{V.~Maffi-Berthier},
\newblock \bibinfo{title}{{Thermodynamics-based Artificial Neural Networks for
  constitutive modeling}},
\newblock \bibinfo{journal}{Journal of the Mechanics and Physics of Solids}
  \bibinfo{volume}{147} (\bibinfo{year}{2021}).
  \href{http://arxiv.org/abs/2005.12183}{{\tt arXiv:2005.12183}}.
\bibitem[{Masi and Stefanou(2022)}]{Masi2021}
\bibinfo{author}{F.~Masi}, \bibinfo{author}{I.~Stefanou},
\newblock \bibinfo{title}{Multiscale modeling of inelastic materials with
  thermodynamics-based artificial neural networks (tann)},
\newblock \bibinfo{journal}{Computer Methods in Applied Mechanics and
  Engineering} \bibinfo{volume}{398} (\bibinfo{year}{2022})
  \bibinfo{pages}{115190}.
  \DOIprefix\doi{https://doi.org/10.1016/j.cma.2022.115190}.
\bibitem[{Liu et~al.(2021)Liu, Tian, Tao, and Yu}]{LIU2021109152}
\bibinfo{author}{X.~Liu}, \bibinfo{author}{S.~Tian}, \bibinfo{author}{F.~Tao},
  \bibinfo{author}{W.~Yu},
\newblock \bibinfo{title}{A review of artificial neural networks in the
  constitutive modeling of composite materials},
\newblock \bibinfo{journal}{Composites Part B: Engineering}
  \bibinfo{volume}{224} (\bibinfo{year}{2021}) \bibinfo{pages}{109152}.
  \DOIprefix\doi{https://doi.org/10.1016/j.compositesb.2021.109152}.
\bibitem[{Liu et~al.(2019)Liu, Wu, and Koishi}]{Liu2019}
\bibinfo{author}{Z.~Liu}, \bibinfo{author}{C.~T. Wu},
  \bibinfo{author}{M.~Koishi},
\newblock \bibinfo{title}{{A deep material network for multiscale topology
  learning and accelerated nonlinear modeling of heterogeneous materials}},
\newblock \bibinfo{journal}{Computer Methods in Applied Mechanics and
  Engineering} \bibinfo{volume}{345} (\bibinfo{year}{2019})
  \bibinfo{pages}{1138--1168}. \DOIprefix\doi{10.1016/j.cma.2018.09.020}.
  \href{http://arxiv.org/abs/1807.09829}{{\tt arXiv:1807.09829}}.
\bibitem[{Liu(2021)}]{Liu2021}
\bibinfo{author}{Z.~Liu},
\newblock \bibinfo{title}{Cell division in deep material networks applied to
  multiscale strain localization modeling},
\newblock \bibinfo{journal}{CoRR} \bibinfo{volume}{abs/2101.07226}
  (\bibinfo{year}{2021}). \href{http://arxiv.org/abs/2101.07226}{{\tt
  arXiv:2101.07226}}.
\bibitem[{Kerfriden et~al.(2013)Kerfriden, Goury, Rabczuk, and
  Bordas}]{KERFRIDEN2013}
\bibinfo{author}{P.~Kerfriden}, \bibinfo{author}{O.~Goury},
  \bibinfo{author}{T.~Rabczuk}, \bibinfo{author}{S.~Bordas},
\newblock \bibinfo{title}{A partitioned model order reduction approach to
  rationalise computational expenses in nonlinear fracture mechanics},
\newblock \bibinfo{journal}{Computer Methods in Applied Mechanics and
  Engineering} \bibinfo{volume}{256} (\bibinfo{year}{2013})
  \bibinfo{pages}{169--188}.
  \DOIprefix\doi{https://doi.org/10.1016/j.cma.2012.12.004}.
\bibitem[{Bessa et~al.(2017)Bessa, Bostanabad, Liu, Hu, Apley, Brinson, Chen,
  and Liu}]{BESSA2017}
\bibinfo{author}{M.~Bessa}, \bibinfo{author}{R.~Bostanabad},
  \bibinfo{author}{Z.~Liu}, \bibinfo{author}{A.~Hu}, \bibinfo{author}{D.~W.
  Apley}, \bibinfo{author}{C.~Brinson}, \bibinfo{author}{W.~Chen},
  \bibinfo{author}{W.~Liu},
\newblock \bibinfo{title}{A framework for data-driven analysis of materials
  under uncertainty: Countering the curse of dimensionality},
\newblock \bibinfo{journal}{Computer Methods in Applied Mechanics and
  Engineering} \bibinfo{volume}{320} (\bibinfo{year}{2017})
  \bibinfo{pages}{633--667}.
  \DOIprefix\doi{https://doi.org/10.1016/j.cma.2017.03.037}.
\bibitem[{Oliver et~al.(2017)Oliver, Caicedo, Huespe, Hernández, and
  Roubin}]{OLIVER2017}
\bibinfo{author}{J.~Oliver}, \bibinfo{author}{M.~Caicedo},
  \bibinfo{author}{A.~Huespe}, \bibinfo{author}{J.~Hernández},
  \bibinfo{author}{E.~Roubin},
\newblock \bibinfo{title}{Reduced order modeling strategies for computational
  multiscale fracture},
\newblock \bibinfo{journal}{Computer Methods in Applied Mechanics and
  Engineering} \bibinfo{volume}{313} (\bibinfo{year}{2017})
  \bibinfo{pages}{560--595}.
  \DOIprefix\doi{https://doi.org/10.1016/j.cma.2016.09.039}.
\bibitem[{Fuhg et~al.(2021)Fuhg, Böhm, Bouklas, Fau, Wriggers, and
  Marino}]{Fuhg2021}
\bibinfo{author}{J.~N. Fuhg}, \bibinfo{author}{C.~Böhm},
  \bibinfo{author}{N.~Bouklas}, \bibinfo{author}{A.~Fau},
  \bibinfo{author}{P.~Wriggers}, \bibinfo{author}{M.~Marino},
\newblock \bibinfo{title}{Model-data-driven constitutive responses: Application
  to a multiscale computational framework},
\newblock \bibinfo{journal}{International Journal of Engineering Science}
  \bibinfo{volume}{167} (\bibinfo{year}{2021}) \bibinfo{pages}{103522}.
  \DOIprefix\doi{https://doi.org/10.1016/j.ijengsci.2021.103522}.
\bibitem[{Ghavamian et~al.(2017)Ghavamian, Tiso, and Simone}]{GHAVAMIAN2017}
\bibinfo{author}{F.~Ghavamian}, \bibinfo{author}{P.~Tiso},
  \bibinfo{author}{A.~Simone},
\newblock \bibinfo{title}{Pod–deim model order reduction for strain-softening
  viscoplasticity},
\newblock \bibinfo{journal}{Computer Methods in Applied Mechanics and
  Engineering} \bibinfo{volume}{317} (\bibinfo{year}{2017})
  \bibinfo{pages}{458--479}.
  \DOIprefix\doi{https://doi.org/10.1016/j.cma.2016.11.025}.
\bibitem[{Rocha et~al.(2018)Rocha, {van der Meer}, and Sluys}]{Rocha2018}
\bibinfo{author}{I.~B. C.~M. Rocha}, \bibinfo{author}{F.~P. {van der Meer}},
  \bibinfo{author}{L.~Sluys},
\newblock \bibinfo{title}{Efficient micromechanical analysis of
  fiber-reinforced composites subjected to cyclic loading through time
  homogenization and reduced-order modeling},
\newblock \bibinfo{journal}{Computer Methods in Applied Mechanics and
  Engineering} \bibinfo{volume}{345} (\bibinfo{year}{2018}).
  \DOIprefix\doi{10.1016/j.cma.2018.11.014}.
\bibitem[{Rocha et~al.(2020)Rocha, Kerfriden, and van~der Meer}]{Rocha2020}
\bibinfo{author}{I.~B. C.~M. Rocha}, \bibinfo{author}{P.~Kerfriden},
  \bibinfo{author}{F.~P. van~der Meer},
\newblock \bibinfo{title}{{Micromechanics-based surrogate models for the
  response of composites: A critical comparison between a classical mesoscale
  constitutive model, hyper-reduction and neural networks}},
\newblock \bibinfo{journal}{European Journal of Mechanics, A/Solids}
  \bibinfo{volume}{82} (\bibinfo{year}{2020}) \bibinfo{pages}{103995}.
  \DOIprefix\doi{10.1016/j.euromechsol.2020.103995}.
\bibitem[{Vijayaraghavan et~al.(2021)Vijayaraghavan, Wu, Noels, Bordas,
  Natarajan, and Beex}]{vijayaraghavan2021neuralnetwork}
\bibinfo{author}{S.~Vijayaraghavan}, \bibinfo{author}{L.~Wu},
  \bibinfo{author}{L.~Noels}, \bibinfo{author}{S.~P.~A. Bordas},
  \bibinfo{author}{S.~Natarajan}, \bibinfo{author}{L.~A.~A. Beex},
  \bibinfo{title}{Neural-network acceleration of projection-based
  model-order-reduction for finite plasticity: Application to rves},
  \bibinfo{year}{2021}. \href{http://arxiv.org/abs/2109.07747}{{\tt
  arXiv:2109.07747}}.
\bibitem[{Guo et~al.(2021)Guo, Rokoš, and Veroy}]{GUO2021}
\bibinfo{author}{T.~Guo}, \bibinfo{author}{O.~Rokoš},
  \bibinfo{author}{K.~Veroy},
\newblock \bibinfo{title}{Learning constitutive models from microstructural
  simulations via a non-intrusive reduced basis method},
\newblock \bibinfo{journal}{Computer Methods in Applied Mechanics and
  Engineering} \bibinfo{volume}{384} (\bibinfo{year}{2021})
  \bibinfo{pages}{113924}.
  \DOIprefix\doi{https://doi.org/10.1016/j.cma.2021.113924}.
\bibitem[{Nguyen et~al.(2012)Nguyen, Lloberas-Valls, Stroeven, and
  Sluys}]{NGUYEN2012}
\bibinfo{author}{V.~P. Nguyen}, \bibinfo{author}{O.~Lloberas-Valls},
  \bibinfo{author}{M.~Stroeven}, \bibinfo{author}{L.~J. Sluys},
\newblock \bibinfo{title}{{Computational homogenization for multiscale crack
  modeling. Implementational and computational aspects}},
\newblock \bibinfo{journal}{International Journal for Numerical Methods in
  Engineering} \bibinfo{volume}{89} (\bibinfo{year}{2012})
  \bibinfo{pages}{192--226}. \DOIprefix\doi{10.1002/nme.3237}.
\bibitem[{Kingma et~al.(2015)Kingma, Salimans, and Welling}]{Kingmaetal2015}
\bibinfo{author}{D.~P. Kingma}, \bibinfo{author}{T.~Salimans},
  \bibinfo{author}{M.~Welling}, \bibinfo{title}{Variational dropout and the
  local reparameterization trick}, \bibinfo{year}{2015}.
  \DOIprefix\doi{10.48550/ARXIV.1506.02557}.
\bibitem[{Kingma and Ba(2014)}]{Adametal2014}
\bibinfo{author}{D.~P. Kingma}, \bibinfo{author}{J.~Ba}, \bibinfo{title}{Adam:
  A method for stochastic optimization}, \bibinfo{year}{2014}.
  \DOIprefix\doi{10.48550/ARXIV.1412.6980}.
\bibitem[{Hernández et~al.(2017)Hernández, Caicedo, and
  Ferrer}]{HERNANDEZ2017}
\bibinfo{author}{J.~Hernández}, \bibinfo{author}{M.~Caicedo},
  \bibinfo{author}{A.~Ferrer},
\newblock \bibinfo{title}{Dimensional hyper-reduction of nonlinear finite
  element models via empirical cubature},
\newblock \bibinfo{journal}{Computer Methods in Applied Mechanics and
  Engineering} \bibinfo{volume}{313} (\bibinfo{year}{2017})
  \bibinfo{pages}{687--722}.
  \DOIprefix\doi{https://doi.org/10.1016/j.cma.2016.10.022}.
\bibitem[{Mori and Tanaka(1973)}]{MORI1973}
\bibinfo{author}{T.~Mori}, \bibinfo{author}{K.~Tanaka},
\newblock \bibinfo{title}{Average stress in matrix and average elastic energy
  of materials with misfitting inclusions},
\newblock \bibinfo{journal}{Acta Metallurgica} \bibinfo{volume}{21}
  (\bibinfo{year}{1973}) \bibinfo{pages}{571--574}.
  \DOIprefix\doi{https://doi.org/10.1016/0001-6160(73)90064-3}.
\bibitem[{Eshelby and Peierls(1957)}]{Eshelby1957}
\bibinfo{author}{J.~D. Eshelby}, \bibinfo{author}{R.~E. Peierls},
\newblock \bibinfo{title}{The determination of the elastic field of an
  ellipsoidal inclusion, and related problems},
\newblock \bibinfo{journal}{Proceedings of the Royal Society of London. Series
  A. Mathematical and Physical Sciences} \bibinfo{volume}{241}
  (\bibinfo{year}{1957}) \bibinfo{pages}{376--396}.
  \DOIprefix\doi{10.1098/rspa.1957.0133}.
  \href{http://arxiv.org/abs/https://royalsocietypublishing.org/doi/pdf/10.1098/rspa.1957.0133}{{\tt
  arXiv:https://royalsocietypublishing.org/doi/pdf/10.1098/rspa.1957.0133}}.
\bibitem[{Nguyen-Thanh et~al.(2020)Nguyen-Thanh, Nguyen, {de Vaucorbeil},
  {Kanti Mandal}, and Wu}]{NGUYENTHANH2020}
\bibinfo{author}{C.~Nguyen-Thanh}, \bibinfo{author}{V.~P. Nguyen},
  \bibinfo{author}{A.~{de Vaucorbeil}}, \bibinfo{author}{T.~{Kanti Mandal}},
  \bibinfo{author}{J.-Y. Wu},
\newblock \bibinfo{title}{Jive: An open source, research-oriented c++ library
  for solving partial differential equations},
\newblock \bibinfo{journal}{Advances in Engineering Software}
  \bibinfo{volume}{150} (\bibinfo{year}{2020}) \bibinfo{pages}{102925}.
  \DOIprefix\doi{https://doi.org/10.1016/j.advengsoft.2020.102925}.
\bibitem[{van~der Meer(2012)}]{vanderMeer2012}
\bibinfo{author}{F.~P. van~der Meer},
\newblock \bibinfo{title}{Mesolevel modeling of failure in composite laminates:
  Constitutive, kinematic and algorithmic aspects},
\newblock \bibinfo{journal}{Archives of Computational Methods in Engineering}
  \bibinfo{volume}{19} (\bibinfo{year}{2012}) \bibinfo{pages}{381--425}.
  \DOIprefix\doi{10.1007/s11831-012-9076-y}.

\end{thebibliography}

\end{document}